%% file: degenerate_diffusion-20200221v2-random_fields-b.tex
\documentclass[a4paper,12pt]{article}
\usepackage{dsfont}
\usepackage{amssymb}
\usepackage{latexsym}
\usepackage{amsmath}
\usepackage{color}
\usepackage{comment}
\usepackage{amsthm}
\usepackage[dvips]{graphicx}
\usepackage{layout} 
\usepackage{ulem}
\usepackage{enumerate}
\usepackage{bm}
\usepackage{bbm} 
\usepackage[bbgreekl]{mathbbol} 
\usepackage{authblk}
\usepackage{setspace} 
\newif\ifrs
\rstrue
\ifrs \usepackage{mathrsfs} \fi  
\newif\ifcol
\coltrue 

\newtheorem{theorem}{Theorem}[section]

\newtheorem{lemma}[theorem]{Lemma}

\newtheorem{proposition}[theorem]{Proposition}

\newtheorem{remark}[theorem]{Remark}

\numberwithin{equation}{section}
\newtheorem{theorem*}{Theorem}
\newtheorem{ass*}[theorem*]{Assumption}
\newtheorem{note*}[theorem*]{Note}
\newtheorem{lemma*}[theorem*]{Lemma}
\newtheorem{definition*}[theorem*]{Definition}
\newtheorem{proposition*}[theorem*]{Proposition}
\newtheorem{corollary*}[theorem*]{Corollary}
\newtheorem{remark*}[theorem*]{Remark}
\newtheorem{example*}[theorem*]{Example}
\numberwithin{equation}{section}
\newif\ifcol
\colfalse
\ifcol
\newcommand{\colorr}{\color[rgb]{0.8,0,0}}
\newcommand{\colorg}{\color[rgb]{0,0.5,0}}
\newcommand{\colorb}{\color[rgb]{0,0,0.8}}

\newcommand{\colorn}{\color[rgb]{1,1,1}}

\newcommand{\coloro}{\color[rgb]{1,0.4,0}}
\newcommand{\coloroy}{\color[rgb]{1,0.95,0}}

\else
\newcommand{\colorb}{\color{black}}

\newcommand{\colorr}{\color{black}}
\newcommand{\colorg}{\color{black}}

\newcommand{\colorn}{\color{black}}
\newcommand{\coloro}{\color{black}}

\newcommand{\coloroy}{\color{black}}
\fi
\newif\ifcol
\colfalse
\ifcol
\newcommand{\cred}{\color[rgb]{0.8,0,0}}

\newcommand{\cblue}{\color[rgb]{0,0,0.8}}
\else
\newcommand{\cred}{\color{black}}
\newcommand{\cblue}{\color{black}}
\fi
\excludecomment{en-text}
\includecomment{jp-text}
\includecomment{comment}
\input nakamacro291103++.tex

\begin{document}

\title{Adaptive and non-adaptive estimation for degenerate diffusion processes
\footnote{
This work was in part supported by 
Japan Science and Technology Agency CREST JPMJCR14D7; 
Japan Society for the Promotion of Science Grants-in-Aid for Scientific Research 
No. 17H01702 (Scientific Research);  
and by a Cooperative Research Program of the Institute of Statistical Mathematics. 
}
}
\author[1]{Arnaud Gloter}
\author[2,3]{Nakahiro Yoshida}
\affil[1]{Laboratoire de Math\'ematiques et Mod\'elisation d'Evry, Universit\'e d'Evry
\footnote{
Laboratoire de Math\'ematiques et Mod\'elisation d'Evry, CNRS, Univ Evry, 
Universit\'e Paris-Saclay, 91037, Evry, France. e-mail: arnaud.gloter@univ-evry.fr}
        }
\affil[2]{Graduate School of Mathematical Sciences, University of Tokyo
\footnote{Graduate School of Mathematical Sciences, University of Tokyo: 3-8-1 Komaba, Meguro-ku, Tokyo 153-8914, Japan. e-mail: nakahiro@ms.u-tokyo.ac.jp}
        }
\affil[3]{CREST, Japan Science and Technology Agency
        }
\maketitle
\ \\
{\it Summary} 
We discuss parametric estimation of a degenerate diffusion system 
from time-discrete observations. 
The first component of the degenerate diffusion system 
has a parameter $\theta_1$ in a non-degenerate diffusion coefficient and 
a parameter $\theta_2$ in the drift term. 
The second component has a drift term parameterized by $\theta_3$ and no diffusion term. 
Asymptotic normality is proved in three different situations for 
an adaptive estimator for $\theta_3$ with some initial estimators for $(\theta_1,\theta_2)$, 
an adaptive one-step estimator for $(\theta_1,\theta_2,\theta_3)$ with some initial estimators 
for {\cred them}, 
and a joint quasi-maximum likelihood estimator for $(\theta_1,\theta_2,\theta_3)$ 
without any initial estimator. 
Our estimators incorporate information of the increments of both components. 
Thanks to this construction, 
the asymptotic variance of the estimators for $\theta_1$ is smaller than the standard one 
based only on the first component. 
The convergence of the estimators for $\theta_3$ is much faster than the other parameters. 
The resulting asymptotic variance is smaller than that of an estimator only using the increments of 
the second component. 

\section{Introduction}\label{202002201325}
In this article, we will discuss {\cred parametric estimation} for 
a hypo-elliptic diffusion process. 
More precisely, given a stochastic basis $(\Omega,\calf,\F,P)$ with a right-continuous filtration $\F=(\calf_t)_{t\in\bbR_+}$, 
$\bbR_+=[0,\infty)$, suppose that an $\F$-adapted 
process $Z_t=(X_t,Y_t)$ satisfies the stochastic differential equation
\bea\label{201905262213}
\left\{\begin{array}{ccl}
dX_t &=& A(Z_t,\theta_2)dt+B(Z_t,\theta_1)dw_t
\y
dY_t&=&H(Z_t,\theta_3)dt
\end{array}\right.
\eea
Here $A:\bbR^{\sfd_Z}\times\overline{\Theta}_2\to\bbR^{\sfd_X}$, 
$B:\bbR^{\sfd_Z}\times\overline{\Theta}_1\to\bbR^{\sfd_X}\otimes\bbR^\sfr$, 
$H:\bbR^{\sfd_Z}\times\overline{\Theta}_3\to\bbR^{\sfd_Y}$, and 
$w=(w_t)_{t\in\bbR_+}$ is an $\sfr$-dimensional $\F$-Wiener process. 
The spaces $\Theta_i$ ($i=1,2,3$) are the unknown parameter spaces of 
{\cred the components of} 
$\theta=(\theta_1,\theta_2,\theta_3)$ to be estimated 
from the data $(Z_{\tj})_{j=0,1,...,n}$, {\cblue where} $\tj=jh$, $h=h_n$ satisfying 
$h\to0$, $nh\to\infty$ and $nh^2\to0$ as $n\to\infty$.

Estimation theory 
has been well developed for {\cred diffusion processes}. 
Even focusing {\cblue on} parametric estimation for ergodic diffusions, there is huge amount of studies: 
Kutoyants \cite{Kutoyants1984,Kutoyants2004,Kutoyants1997}, 
Prakasa Rao \cite{PrakasaRao1983, PrakasaRao1988},
Yoshida \cite{Yoshida1992b,yoshida2011polynomial},
Bibby and S{\o}rensen \cite{BibbySoerensen1995},
Kessler \cite{Kessler1997},
K{\"u}chler and S{\o}rensen \cite{kuchler1997exponential},
{\cblue Genon--Catalot et al. \cite{genoncatalot1999parameter},}
Gloter \cite{gloter2000discrete,gloter2001parameter,gloter2007efficient}, 
Sakamoto and Yoshida \cite{SakamotoYoshida2009}, 
Uchida \cite{uchida2010contrast}, 
Uchida and Yoshida \cite{UchidaYoshida2001,uchida2011estimation, uchida2012adaptive}, 
Kamatani and Uchida \cite{KamataniUchida2014}, 
De Gregorio and Iacus \cite{de2012adaptive}, 
{\cblue Genon--Catalot and Lar\'edo \cite{genoncatalot2016estimation},}
Suzuki and Yoshida \cite{suzuki2020penalized} 
among many others. 
{\cred Nakakita and Uchida \cite{nakakita2019inference} and 
Nakakita et al. \cite{nakakita2020quasi} studied estimation under measurement error; 
related are 
Gloter and Jacod \cite{gloter2001diffusions1,gloter2001diffusions}. 
}
\noindent
{\cblue Non parametric estimation for the coefficients of an ergodic diffusion has also been widely studied : Dalayan and Kutoyants \cite{dalalyan2002asymptotically}, Kutoyants \cite{Kutoyants2004}, Dalalyan \cite{dalalyan2005sharp}, Dalalyan and Rei{ss} \cite{dalalyan2006asymptotic, dalalyan2007asymptotic}, Comte and Genon--Catalot \cite{comte2006penalized}, Comte et al. \cite{comte2007penalized}, Schmisser \cite{schmisser2013penalized}, to name a few}.
Historically attentions were paid to inference for non-degenerate cases.  

Recently there is a growing interest in hypo-elliptic diffusions, that appear in various applied fields. 
Examples of the hypo-elliptic diffusion include 
the harmonic oscillator, the Van der Pol oscillator and 
the FitzHugh-Nagumo neuronal model; see e.g. 
Le\'on and Samson \cite{leon2018hypoelliptic}. 
%
For parametric estimation of {\cblue hypo-elliptic} diffusions, we refer the reader to 
Gloter \cite{gloter2006parameter} for a discretely observed integrated diffusion process, 
{\cred and}
Samson and Thieullen \cite{samson2012contrast} for a 
contrast estimator. 
Comte et al. \cite{comte2017adaptive} gave adaptive estimation under partial observation. 
Recently, 
Ditlevsen and Samson \cite{ditlevsen2019hypoelliptic} studied 
filtering and inference for {\cblue hypo-elliptic} diffusions from complete and partial observations. 
When the observations are discrete and complete, 
they showed asymptotic normality of their estimators under the assumption that 
the true value of some of parameters are known. 
{\cblue Melnykova \cite{melnykova2019parametric} studied 
	the estimation problem for the model \eqref{201905262213}, comparing contrast functions and least square estimates. The 
	contrast functions we propose in this paper are different from the one
	 in \cite{melnykova2019parametric}. 
}	
	

In this paper, we will present several estimation schemes. 
Since we assume discrete-time observations of $Z{\cred=(Z_t)_{t\in\bbR+}}$, 
quasi-likelihood estimation for $\theta_1$ and $\theta_2$ {\cred is} known; 
only difference from the standard diffusion case is the existence of 
the covariate $Y{\cred=(Y_t)_{t\in\bbR+}}$ 
in the equation of $X{\cred=(X_t)_{t\in\bbR+}}$ but it causes no theoretical difficulty. 
We will give an exposition for construction of those standard estimators 
in Sections \ref{201905291607} and \ref{201906041938} for selfcontainedness. 
Thus, our first approach in Section \ref{201906041935} 
is toward estimation of $\theta_3$ with initial estimators 
for $\theta_1$ and $\theta_2$. 
The idea for construction of the quasi-likelihood function in the elliptic case was 
based on the local Gaussian approximation of the transition density. 
Then it is natural to approximate the distribution of the increments of $Y$ 
by {\cred that of} the principal Gaussian variable in the expansion of {\cred the increment}. 
However, this method causes deficiency, as we will observe there; {\cred see} 
Remark \ref{201906181743} on p.\pageref{201906181743}. 
We present a more efficient method by incorporating an additional Gaussian part from $X$. 
The rate of convergence attained by the estimator for $\theta_3$ is $n^{-1/2}h^{1/2}$ and it is 
much faster than the rate $(nh)^{-1/2}$ for $\theta_2$ and $n^{-1/2}$ for $\theta_1$. 
Section \ref{202001141623} treats some adaptive estimators 
using suitable initial estimators for $(\theta_1,\theta_2,\theta_3)$, and 
shows joint {\cred asymptotic normality}. 
Then it should be remarked that the asymptotic variance of our estimator $\hat{\theta}_1$ 
for $\theta_1$ has improved {\cred that of} the ordinary volatility parameter estimator, 
e.g. $\hat{\theta}_1^0$ recalled in Section \ref{201905291607}, 
that would be 
asymptotically optimal if the system consisted only of $X$. 
In Section \ref{202001141632}, 
we consider a non-adaptive joint quasi-maximum likelihood estimator. 
This method does not require initial estimators. 
From computational point of view, adaptive methods often have merits 
by reducing dimension of parameters, 
but 
the non-adaptive method is still theoretically interesting. 
{\cred
Section \ref{202001141611} collects the assumptions under which we will work. 
Section \ref{202001141618} offers several basic estimates to the increments of $Z$.  
}
To investigate efficiency of the presented estimators, {\cred we need} 
the LAN property of the exact likelihood function of the hypo-elliptic diffusion. 
We will discuss this problem elsewhere.

\section{Assumptions}\label{202001141611}
We assume that $\Theta_i$ ($i=1,2,3$) are bounded open domain in $\bbR^{\sfp_i}$, respectively, and 
${\colorg\Theta=\prod_{i=1}^3}\Theta_i$ has a good boundary so that Sobolev's embedding inequality holds, that is, 
there exists a positive constant $C_{{\colorg\Theta}}$ such that 
{\colorg
\bea\label{202001131246}
\sup_{\theta\in\Theta}|f_i(\theta_i)|&\leq& C_{\Theta}\sum_{k=0}^1 
\|\partial_\theta^kf\|_{L^p(\Theta)}
\eea
for all $f\in C^1(\Theta)$ and $p>\sum_{i=1}^3\sfp_i$. 
If $\Theta$ has a Lipschitz boundary, then this condition is satisfied. 
Obviously, the embedding inequality (\ref{202001131246}) is valid for functions depending only on a part of components of $\theta$.} 

In this paper, we will propose an estimator for $\theta$ and show its consistency and asymptotic normality. 

Given a finite-dimensional real vector space ${\sf E}$, 
denote by $C^{a,b}_p(\bbR^{\sfd_Z}\times\Theta_i;{\sf E})$ the set of functions $f:\bbR^{\sfd_Z}\times\Theta_i\to{\sf E}$ 
such that $f$ is continuously differentiable $a$ times in $z\in\bbR^{\sfd_Z}$ and $b$ times in $\theta_i\in\Theta$ 
in any order and $f$ and all such derivatives are continuously extended to $\bbR^{\sfd_Z}\times\overline{\Theta}_i$, 
moreover, they are of at most polynomial growth in $z\in\bbR^{\sfd_Z}$ uniformly in $\theta\in\Theta$. 
Let ${\colorr C}=BB^\star$, $\star$ denoting the matrix transpose. 
We suppose that the process $(Z_t)_{t\in\bbR_+}$ that generates the data satisfies 
the stochastic differential equation (\ref{201905262213}) for a true value 
$\theta^*=(\theta_1^*,\theta_2^*,\theta_3^*)\in\Theta_1\times\Theta_2\times\Theta_3$. 

\bd
\im[[A1\!\!]]
{\bf (i)} 
$A\in C^{i_A,j_A}_p(\bbR^{\sfd_Z}\times\Theta_2;\bbR^{\sfd_X})$ and 
$B\in C^{i_B,j_B}_p(\bbR^{\sfd_Z}\times\Theta_1;\bbR^{\sfd_X}\otimes\bbR^\sfr)$. 

\bd\im[(ii)] $H\in C^{i_H,j_H}_p(\bbR^{\sfd_Z}\times\Theta_3;\bbR^{\sfd_Y})$. 
\ed
\ed
We will denote 
$F_x$ for $\partial_xF$, $F_y$ for $\partial_yF$, and $F_i$ for $\partial_{\theta_i}F$. 

\bd
\im[[A2\!\!]] {\bf (i)} 
$\sup_{t\in\bbR_+}\|Z_t\|_p<\infty$ 
for every $p>1$. 
\bd\im[(ii)] There exists a probability measure $\nu$ on $\bbR^{\sfd_Z}$ such that 
\beas 
\frac{1}{T}\int_0^T f(Z_t)\>dt &\to^p& \int f(z)\nu(dz)\qquad(T\to\infty)
\eeas
for any {\cred bounded continuous} 
function $f:\bbR^{\sfd_Z}\to\bbR$. 

\im[(iii)] {\coloro The function $\theta_1\mapsto C(Z_t,\theta_1)^{-1}$ is continuous 
{\colorg on $\overline{\Theta}_1$}
a.s., and} 
\beas
\sup_{\theta_1\in\Theta_1}\sup_{t\in\bbR_+}\|\det C(Z_t,\theta_1)^{-1}\|_p&<&\infty
\eeas 
for every $p>1$. 

\im[(iv)] For the $\bbR^{{\colorr\sfd_Y}}\otimes\bbR^{{\colorr\sfd_Y}}$ valued function 
$V(z,\theta_1,\theta_3)=H_x(z,\theta_3)C(z,\theta_1)H_x(z,\theta_3)^\star$, 
{\coloro the function $(\theta_1,\theta_3)\mapsto V(Z_t,\theta_1,\theta_3)^{-1}$ is continuous 
{\colorg on $\overline{\Theta}_1\times\overline{\Theta}_3$} a.s., 
and }
\beas 
\sup_{(\theta_1,\theta_3)\in\Theta_1\times\Theta_3}\sup_{t\in\bbR_+}\|\det V(Z_t,\theta_1,\theta_3)^{-1}\|_p&<&\infty
\eeas 
for every $p>1$. 
\ed
\ed
\begin{remark}\label{201906021604}\rm 
{\bf (a)} 
It follows from $[A2]$ that 
the convergence in $[A2]$ (ii) holds for any {\cred continuous} 
function $f$ of at most polynomial growth. 
\begin{en-text}
By the compactness of $\overline{\Theta}_1\times\overline{\Theta}_3$ and the continuity of the function, 
\beas 
\inf_{(\theta_1,\theta_3)\in\overline{\Theta}_1\times\overline{\Theta}_3}\det V(Z_t,\theta_1,\theta_3)
&=&
\liminf_{m\to\infty}\det V(Z_t,\theta_1^m,\theta_3^m)
\eeas
for a sequence $(\theta_1^m,\theta_3^m)_{m\in\bbN}$ that is dense in $\overline{\Theta}_1\times\overline{\Theta}_3$. 
Then by Fatou's lemma, under $[A2]$ (iv), 
\beas 
P\bigg[\inf_{(\theta_1,\theta_3)\in\overline{\Theta}_1\times\overline{\Theta}_3}\det V(Z_t,\theta_1,\theta_3)>0\bigg]=1
\eeas 
for all $t\in\bbZ_+$, 
and hence  
the function $(\theta_1,\theta_3)\mapsto V(Z_t,\theta_1,\theta_3)^{-1}$ is continuous a.s. 
(Otherwise, we could construct an estimating function avoiding the null set of $V(Z_\tjm,\theta_1,\theta_3)$, as a matter of fact, if we do not mind information loss. ) 
In the same way, the function $\theta_1\mapsto C(Z_t,\theta_1)^{-1}$ is continuous a.s. by $[A2]$ (iii). 
\end{en-text}
\bd\im[(b)] 
{\coloro We implicitly assume the existence of $C(Z_T,\theta_1)^{-1}$ and $V(Z_t,\theta_1,\theta_3)^{-1}$ 
in (iii) and (iv) of $[A2]$. }
\im[(c)] 
Fatou's lemma implies 
\beas 
\int |z|^p\nu(dz)+\sup_{\theta_1\in\overline{\Theta}_1}\int \big(\det C(z,\theta_1)\big)^{-p}\nu(dz)
+\sup_{(\theta_1,\theta_3)\in\overline{\Theta}_1\times\overline{\Theta}_3}
\int \big(\det V(z,\theta_1,\theta_3)\big)^{-p}\nu(dz)<\infty
\eeas 
for any $p>0$. 
\ed
\end{remark}

\begin{en-text}
\bi
\im $\sup_{\theta\in\Theta}\|R_j(\theta)\|_p=O(h^2)$ for every $p>1$. 

\im Conditionally on $\calf_\tjm$, the distribution of $\Delta_jY$ is approximated by 
the normal distribution\\ $N\big(hF_n(Z_\tjm,\theta_3),h^3K(Z_\tjm,(\theta_1,\theta_3)\big)$. 

\im This suggests the use of the quasi-log likelihood function 
\beas 
\bbH^{[3]}_n(\theta_3)
 &=& 
 -\half\sum_{j=1}^n \bigg\{
 \frac{\big(\Delta_jY-hF_n(Z_\tjm,\theta_3)\big)^2}{h^3K(Z_\tjm,\hat{\theta}_1,\theta_3)}
 +\log\big(h^3K(Z_\tjm,\hat{\theta}_1,\theta_3)\big)\bigg\}
 \eeas
where 
$\hat{\theta}_1=\hat{\theta}_{1,n}$ and $\hat{\theta}_2^0=\hat{\theta}_{2,n}$ implicitly in $F_n$ 
are consistent estimators for $\theta_1$ and $\theta_2$ respectively based on the data $(Z_\tj)_{j=0,1,...,n}$. 
\im $B$ is non-degenerate. 
\im Observe $(Z_\tj)_{j=0,1,...,n}$, $\tj=t_j^n=jh$, $h=h_n=\tj-\tjm$
\im Assume $h\to0$, $nh\to\infty$, $nh^2\to0$. 
\im Approximation to $\Delta_jY=Y_\tj-Y_\tjm$
\ei
\end{en-text}

Let 
\beas 
\bbY^{(1)}(\theta_1) 
&=& 
-\half\int \bigg\{
\text{Tr}\big(C(z,\theta_1)^{-1}C(z,\theta_1^*)\big)-\sfd_X
+\log\frac{\det C(z,\theta_1)}{\det C(z,\theta_1^*)}\bigg\}\nu(dz).
\eeas
Since 
$|\log x|\leq x+x^{-1}$ for $x>0$, $\bbY^{(1)}(\theta_1)$ is a continuous function on $\overline{\Theta}_1$ 
well defined 
under {\colorg $[A1]$ and }$[A2]$. 
Let 
\beas 
\bbY^{(J,1)}(\theta_1)
&=& 
-\half\int \bigg\{\text{Tr}\big(C(z,\theta_1)^{-1}C(z,\theta_1^*)\big)
+\text{Tr}\big(V(z,\theta_1,\theta_3^*)^{-1}V(z,\theta_1^*,\theta_3^*)\big)
-\sfd_Z
\\&&
+\log \frac{\det C(z,\theta_1)\det V(\theta_1,\theta_3^*)}{\det C(z,\theta_1^*)\det V(\theta_1^*,\theta_3^*)}
\bigg\}
\nu(dz)
\eeas
Let 
\bea\label{202002171821}
\bbY^{(2)}(\theta_2)
&=&
-\half \int 
C(z,\theta_1^*)^{{\cred -1}}\big[\big(A(z,\theta_2)-A(z,\theta_2^*)\big)^{\otimes2}\big]\nu(dz)
\eea
Let 
\beas 
\bbY^{(3)}(\theta_3) &=&
-\int 6V(z,\theta_1^*,\theta_3)^{-1}
\big[\big({\colorr H(z,\theta_3)-H(z,\theta_3^*)}\big)^{\otimes2}\big]\nu(dz). 
\eeas
The random field $\bbY^{(3)}$ is well defined under {\colorg $[A1]$ and }$[A2]$. 
{\coloro Let 
\beas 
\bbY^{(J,3)}(\theta_1,\theta_3) &=&
-\int 6V(z,\theta_1,\theta_3)^{-1}
\big[\big(H(z,\theta_3)-H(z,\theta_3^*)\big)^{\otimes2}\big]\nu(dz). 
\eeas
}

We will assume all or some of the following identifiability conditions
\bd
\im[[A3\!\!]] 
{\bf (i)} There exists a positive constant $\chi_1$ such that 
\beas 
\bbY^{(1)}(\theta_1) &\leq& -\chi_1|\theta_1-\theta_1^*|^2\qquad(\theta_1\in\Theta_1). 
\eeas
\bd 
\im[(i$'$)] There exists a positive constant $\chi_1'$ such that 
\beas 
\bbY^{(J,1)}(\theta_1)&\leq& -\chi_1'|\theta_1-\theta_1^*|^2\qquad(\theta_1\in\Theta_1). 
\eeas

\im[(ii)] There exists a positive constant $\chi_2$ such that 
\beas 
\bbY^{(2)}(\theta_2) &\leq& -\chi_2|\theta_2-\theta_2^*|^2\qquad(\theta_2\in\Theta_2). 
\eeas

\im[(iii)]
There exists a positive constant $\chi_3$ such that 
\beas 
\bbY^{(3)}(\theta_3) &\leq& -\chi_3|\theta_3-\theta_3^*|^2\qquad(\theta_3\in\Theta_3). 
\eeas

{\coloro
\im[(iii$'$)]
There exists a positive constant $\chi_3$ such that 
\beas 
\bbY^{(J,3)}(\theta_1,\theta_3) &\leq& -\chi_3|\theta_3-\theta_3^*|^2
\qquad(\theta_1\in\Theta_1,\>\theta_3\in\Theta_3). 
\eeas
}

\ed
\ed

\section{Basic estimation of the increments}\label{202001141618}
We denote $U^{\otimes k}$ for $U\otimes \cdots\otimes U$ ($k$-times) for a tensor $U$. 
For tensors  
$S^1=(S^1_{i_{1,1},...,i_{1,d_1};j_{1,1},...,j_{1,k_1}})$, ..., 
$S^m=(S^m_{i_{m,1},...,i_{m,d_m};j_{m,1},...,j_{m,k_m}})$ 
and and a tensor $T=(T^{i_{1,1},...,i_{1,d_1},...,i_{m,1},...,i_{m,d_m}})$, we write 
\beas 
T[S^1,...,S^m]&=&T[S^1\otimes\cdots\otimes S^m]
\\&=&
\bigg(
\sum_{i_{1,1},...,i_{1,d_1},...,i_{m,1},...,i_{m,d_m}}
T^{i_{1,1},...,i_{1,d_1},...,i_{m,1},...,i_{m,d_m}}
S^1_{i_{1,1},...,i_{1,d_1};j_{1,1},...,j_{1,k_1}}
\\&&\hspace{3cm}\cdots S^m_{i_{m,1},...,i_{m,d_m};j_{m,1},...,j_{m,k_m}}\bigg)
_{j_{1,1},...,j_{1,k_1},...,j_{m,1},...,j_{m,k_m}}.
\eeas
This notation will be applied for a tensor-valued tensor $T$ as well.

We have 
\bea\label{201906021458}
h^{-1/2}\Delta_jX
&=& 
h^{-1/2}\int_\tjm^\tj B(Z_t,\theta_1^*)dw_t
+h^{-1/2}\int_\tjm^\tj A(Z_t,\theta_2^*)dt
\nn\\&=& 
h^{-1/2}B(Z_\tjm,\theta_1^*)\Delta_jw
+r^{(\ref{201905280538})}_j
\eea
where 
\bea\label{201905280538}
r^{(\ref{201905280538})}_j 
&=&
h^{-1/2}\int_\tjm^\tj (B(Z_t,\theta_1^*)-B(Z_\tjm,\theta_1^*))dw_t
+h^{-1/2}\int_\tjm^\tj A(Z_t,\theta_2^*)dt
\eea

\begin{lemma}\label{201906041850}
{\bf (a)} Under $[A1]$ with $(i_A,j_A,i_B,j_B,i_H,j_H)=(0,0,0,0,0,0)$ and $[A2]$ $(i)$, 
\bea\label{201905280550}
\sup_{s,t\in\bbR_+,\>|s-t|\leq \Delta} \|Z_s-Z_t\|_p=O(\Delta^{1/2})\quad(\Delta\down0)
\eea
for every $p>1$. 
\bd\im[(b)] 
Under $[A1]$ with $(i_A,j_A,i_B,j_B,i_H,j_H)=(0,0,1,0,0,0)$ and $[A2]$ $(i)$, 
$r^{(\ref{201905280538})}_j=\overline{O}_{L^\inftym}(h^{1/2})$, i.e., 
\beas 
\sup_n\sup_j\|r^{(\ref{201905280538})}_j\|_p&=&O(h^{1/2})
\eeas 
for every $p>1$. 
\ed
\end{lemma}
\proof 
(a) is trivial. For (b), the first term on the right-hand side of (\ref{201905280538}) can be estimated by 
the Burkholder-Davis-Gundy inequality, Taylor's formula for 
$B(Z_t,\theta_1^*)-B(Z_\tjm,\theta_1^*)$ and by (\ref{201905280550}). 
\qed\halflineskip

We have 
\beas 
h^{-1/2}\Delta_jX
&=& 
h^{-1/2}\int_\tjm^\tj B(Z_t,\theta_1^*)dw_t
+h^{-1/2}\int_\tjm^\tj A(Z_t,\theta_2^*)dt
\\&=& 
h^{-1/2}\int_\tjm^\tj B(Z_t,\theta_1^*)dw_t
+
h^{1/2}A(Z_\tjm,\theta_2^*)
+r^{(\ref{201905280600})}_j
\eeas
where
\bea\label{201905280600}
r^{(\ref{201905280600})}_j
&=&
h^{-1/2}\int_\tjm^\tj \big(A(Z_t,\theta_2^*)-A(Z_\tjm,\theta_2^*)\big)dt
\eea
Then 

\begin{lemma}\label{201906041900}
$r^{(\ref{201905280600})}_j=\overline{O}_{L^\inftym}(h)$, i.e., 
\bea\label{201905280613}
\sup_n\sup_j\big\|r^{(\ref{201905280600})}_j\big\|_p=O(h)
\eea 
for every $p>1$ 
if $[A1]$ for $(i_A,j_A,i_B,j_B,i_H,j_H)=(1,0,0,0,0,0)$ and $[A2]$ $(i)$ hold. 
\end{lemma}
\proof 
Thanks to (\ref{201905280550}). 
\qed\halflineskip

Let 
\beas
L_H(z,\theta_1,\theta_2,\theta_3)
&=&
H_x(z,\theta_3)[A(z,\theta_2)]+\half H_{xx}(z,\theta_3)[C(z,\theta_1)]+H_y(z,\theta_3)[H(z,\theta_3)].
\eeas
Define the $\bbR^{\sfd_Y}$-valued function $G_n(z,\theta_1,\theta_2,\theta_3)$ by 
\beas 
G_n(z,\theta_1,\theta_2,\theta_3) 
&=& 
H(z,\theta_3)+\frac{h}{2} L_H\big(z,\theta_1,\theta_2,\theta_3\big).
\eeas
Write 
\beas 
\zeta_j &=& \sqrt{{\colorr{3}}}\int_\tjm^\tj\int_\tjm^t dw_sdt
\eeas
Then 
$
E\big[\zeta_j^{\otimes2}\big] = h^3I_{\sfr} 
$
for the $\sfr$-dimensional identity matrix $I_{\sfr}$. 

We have 
\bea\label{201905291809}&&
\Delta_jY -hG_n(Z_\tjm,\theta_1,\theta_2,\theta_3)
\nn\\&=&
\Delta_jY-hH(Z_\tjm,\theta_3)-\frac{h^2}{2} L_H\big(Z_\tjm,\theta_1,\theta_2,\theta_3\big)
\nn\\&=&
hH(Z_\tjm,\theta_3^*)-hH(Z_\tjm,\theta_3)
\nn\\&&
+H_x(Z_\tjm,\theta_3^*)B(Z_\tjm,\theta_1^*)\int_\tjm^\tj\int_\tjm^t dw_sdt
\nn\\&&
+\int_\tjm^\tj\int_\tjm^t \big\{H_x(Z_s,\theta_3^*)B(Z_s,\theta_1^*)-H_x(Z_\tjm,\theta_3^*)B(Z_\tjm,\theta_1^*)\big\}dw_sdt
\nn\\&&
+\int_\tjm^\tj\int_\tjm^t \big(L_H(Z_s,\theta_1^*,\theta_2^*,\theta_3^*)-L_H(Z_\tjm,\theta_1,\theta_2,\theta_3)\big)dsdt
\nn\\&=&
\big\{hH(Z_\tjm,\theta_3^*)-hH(Z_\tjm,\theta_3)\big\}
+\kappa(Z_\tjm,\theta_1^*,\theta_3^*)\zeta_j+\rho_j(\theta_1,\theta_2,\theta_3)
\eea
where 
\beas 
\kappa(Z_\tjm,\theta_1^*,\theta_3^*)
&=&
3^{-1/2}H_x(Z_\tjm,\theta_3^*)B(Z_\tjm,\theta_1^*)
\eeas
and
\bea\label{201906021454}
\rho_j(\theta_1,\theta_2,\theta_3)
&=&
\int_\tjm^\tj\int_\tjm^t \big\{H_x(Z_s,\theta_3^*)B(Z_s,\theta_1^*)-H_x(Z_\tjm,\theta_3^*)B(Z_\tjm,\theta_1^*)\big\}dw_sdt
\nn\\&&
+\int_\tjm^\tj\int_\tjm^t \big(L_H(Z_s,\theta_1^*,\theta_2^*,\theta_3^*)-L_H(Z_\tjm,\theta_1,\theta_2,\theta_3)\big)dsdt.
\eea

Let 
\bea\label{201906030041}
\cald_j(\theta_1,\theta_2,\theta_3) &=&
\left(\begin{array}{c}
h^{-1/2}\big(\Delta_jX-{\colorr h}A(Z_\tjm,\theta_2)\big)\y
h^{-3/2}\big(\Delta_jY-{\colorr h}G_n(Z_\tjm,\theta_1,\theta_2,\theta_3)\big)
\end{array}\right).
\eea

\begin{lemma}\label{201906021438}
Suppose that $[A1]$ with $(i_A,j_A,i_B,j_B,i_H,j_H)=(1,0,1,0,3,0)$ and $[A2]$ $(i)$ are satisfied. 
Then 
\bd
\im[(a)] 
$\ds
\sup_n\sup_j\big\|\rho_j(\theta_1^*,\theta_2^*,\theta_3^*)\big\|_p\yeq O(h^2)
$
for every $p>1$. 
\im[(b)] 
$\ds \sup_n\sup_j\big\|\cald_j(\theta_1^*,\theta_2^*,\theta_3^*)\big\|_p\yl\infty$ for every $p>1$. 
\ed
\end{lemma}
\proof 
It is possible to show (a) by 
(\ref{201906021454}) 
and using the estimate (\ref{201905280550}) with the help of Taylor's formula. 
Additionally to the representation (\ref{201905291809}), by 
using (\ref{201906021458}) and (\ref{201905280538}), 
we obtain (b). \qed\halflineskip

\begin{en-text}
\beas 
h^{-1/2}\Delta_jX
&=& 
h^{-1/2}\int_\tjm^\tj B(Z_t,\theta_1^*)dw_t
+h^{-1/2}\int_\tjm^\tj A(Z_t,\theta_2^*)dt
\\&=& 
h^{-1/2}B(Z_\tjm,\theta_1^*)\Delta_jw
+h^{-1/2}\int_\tjm^\tj (B(Z_t,\theta_1^*)-B(Z_\tjm,\theta_1^*))dw_t
\\&&
+h^{-1/2}\int_\tjm^\tj A(Z_t,\theta_2^*)dt
\\&=& 
h^{-1/2}B(Z_\tjm,\theta_1^*)\Delta_jw
+h^{-1/2}\int_\tjm^\tj \int_\tjm^t B_x(Z_s,\theta_1^*)B(Z_s,\theta_1^*)dw_sdw_t
\\&&
+h^{-1/2}\int_\tjm^\tj \int_\tjm^t L_B(Z_s,\theta_1^*,\theta_2^*,\theta_3^*)dsdw_t
+h^{-1/2}\int_\tjm^\tj A(Z_t,\theta_2^*)dt, 
\eeas
\end{en-text}

We denote by $(B_xB)(z,\theta_2)$ the tensor defined by
${\cred (B_xB)}
(z,\theta_2)[u_1\otimes u_2]=
B_x(z,\theta_2)[{\cred u_2},B(z,\theta_2)[{\cred u_1}]]
$ for $u_1,u_2\in\bbR^\sfr$.  
Moreover, we write $dw_sdw_t$ for $dw_s\otimes dw_t$, and 
$(B_xB)(Z_\tjm,\theta_2^*)\int_\tjm^\tj \int_\tjm^tdw_sdw_t$ 
for $(B_xB)(Z_\tjm,\theta_2^*)\big[\int_\tjm^\tj \int_\tjm^tdw_sdw_t\big]$. 
We will apply this rule in similar situations. 
Let 
\bea\label{201996160909}
L_B(z,\theta_1,\theta_2,\theta_3)
&=&
B_x(z,\theta_1)[A(z,\theta_2)]+\half B_{xx}(z,\theta_3)[C(z,\theta_1)]+B_y(z,\theta_3)[H(z,\theta_3)].
\eea

\begin{lemma}\label{201906031518}
Suppose that $[A1]$ with $(i_A,j_A,i_B,j_B,i_H,j_H)=(1,{\colorr1},2,0,0,0)$ and $[A2]$ $(i)$ are satisfied. 
Then 
\bea\label{201906030207} 
h^{-1/2}\big(\Delta_jX-hA(Z_\tjm,\theta_2)\big)
&=& 
\xi^{(\ref{201906030126})}_j+\xi^{(\ref{201906030150})}_j
+r^{(\ref{201906030137})}_j(\theta_2)
\eea
where 
\bea\label{201906030126}
\xi^{(\ref{201906030126})}_j
&=&
h^{-1/2}B(Z_\tjm,\theta_1^*)\Delta_jw,
\eea
\bea\label{201906030150}
\xi^{(\ref{201906030150})}_j
&=&
h^{-1/2}(B_xB)(Z_\tjm,\theta_1^*)\int_\tjm^\tj \int_\tjm^tdw_sdw_t,
\eea
and 
\bea\label{201906030137}
r^{(\ref{201906030137})}_j(\theta_2)
&=&
h^{-1/2}\int_\tjm^\tj \int_\tjm^t((B_xB)(Z_s,\theta_1^*)-(B_xB)(Z_\tjm,\theta_1^*))dw_sdw_t
\nn\\&&
+h^{-1/2}\int_\tjm^\tj\int_\tjm^tL_B(Z_s,\theta_1^*,\theta_2^*,\theta_3^*)dsdw_t
\nn\\&&
+h^{-1/2}\int_\tjm^\tj \big(A(Z_t,\theta_2^*)-A(Z_\tjm,\theta_2)\big)dt.
\eea
Moreover, 
\bea\label{201906030210}
\sup_n\sup_j\|r^{(\ref{201906030137})}_j(\theta_2^*)\|_p&=&O(h)
\eea
for every $p>1$, and 
{\cred
\bea\label{201906031424}
\big|r^{(\ref{201906030137})}_j(\theta_2)\big|
&\leq&
r_{n,j}^{(\ref{201906031425})}\big\{h^{1/2}\big|\theta_2-\theta_2^*\big|+h\big\}
\eea
}
with 
some random variables $r_{n,j}^{(\ref{201906031425})}$ satisfying 
\bea\label{201906031425}
\sup_n\sup_j\big\|r_{n,j}^{(\ref{201906031425})}\big\|_p
&<&
\infty
\eea
for every $p>1$. 
\end{lemma}
\proof 
The decomposition (\ref{201906030207}) is obtained by It\^o's formula. 
The estimate (\ref{201906030210}) is verified by (\ref{201905280550}) 
since $\partial_z(B_xB)$ and $\partial_zA$ are bound by a polynomial in $z$ 
uniformly in $\theta$. 
The estimate (\ref{201906031424}) uses $\partial_2A$ for $\theta_2$ near $\theta_2^*$ 
as well as $\partial_zA$ evaluated at $\theta_2^*$: 
{\cred
\beas
\big|r^{(\ref{201906030137})}_j(\theta_2)\big|1_{\{|\theta_2-\theta_2^*|<r\}}
&\leq&
r_{n,j}^{(\ref{201906031425})}\big\{h^{1/2}\big|\theta_2-\theta_2^*\big|+h\big\}
1_{\{|\theta_2-\theta_2^*|<r\}}
\eeas
with some positive constant $r$ and  
some random variables $r_{n,j}^{(\ref{201906031425})}$ satisfying (\ref{201906031425}).
}
{\colorg The small number $r$ was taken to ensure convexity of the vicinity of $\theta_2^*$.}
{\cred For $\theta_2$ such that $|\theta_2-\theta_2^*|\geq r$, the estimate (\ref{201906031424}) is valid 
by enlarging $r_{n,j}^{(\ref{201906031425})}$ if necessary. }
\qed\halflineskip
\begin{en-text}
\beas &&
h^{-1/2}\big(\Delta_jX-hA(Z_\tjm,\theta_2)\big)
\\&=& 
h^{-1/2}\int_\tjm^\tj B(Z_t,\theta_2^*)dw_t
+h^{-1/2}\int_\tjm^\tj A(Z_t,\theta_2^*)dt-h^{1/2}A(Z_\tjm,\theta_2)
\\&=& 
h^{-1/2}B(Z_\tjm,\theta_2^*)\Delta_jw
+h^{-1/2}\int_\tjm^\tj (B(Z_t,\theta_2^*)-B(Z_\tjm,\theta_2^*))dw_t
\\&&
+h^{-1/2}\int_\tjm^\tj \big(A(Z_t,\theta_2^*)-A(Z_\tjm,\theta_2)\big)dt
\\&=& 
h^{-1/2}B(Z_\tjm,\theta_2^*)\Delta_jw
+h^{-1/2}B_x(Z_\tjm,\theta_2^*)B(Z_\tjm,\theta_2^*)\int_\tjm^\tj \int_\tjm^tdw_sdw_t
+r^{(\ref{201906030137})}_j(\theta_2)
\eeas
$\sup_n\sup_j\|r_j(\theta_2^*)\|_p=O(h)$
\end{en-text}

\begin{lemma}\label{201906031516}
{\bf (a)} 
Suppose that $[A1]$ with $(i_A,j_A,i_B,j_B,i_H,j_H)=(1,1,2,1,3,0)$ and $[A2]$ $(i)$ are satisfied. 
Then 
\bea\label{201906030236}
\Delta_jY -hG_n(Z_\tjm,\theta_1,\theta_2,\theta_3^*)
&=&
\xi_j^{(\ref{201906030334})}+\xi_j^{(\ref{201906030345})}
+h^{3/2}r_j^{(\ref{201906030335})}(\theta_1,\theta_2)
+h^{3/2}r_j^{(\ref{201906030336})}(\theta_1,\theta_2)
\nn\\&&
\eea
where
\bea\label{201906030334}
\xi_j^{(\ref{201906030334})}
&=&
\kappa(Z_\tjm,\theta_1^*,\theta_3^*)\zeta_j,
\eea
\bea\label{201906030345}
\xi_j^{(\ref{201906030345})}
&=&
((H_xB)_xB)(Z_\tjm,\theta_1^*,\theta_3^*)\int_\tjm^\tj\int_\tjm^t\int_\tjm^sdw_rdw_sdt,
\eea
\bea\label{201906030335}
r_j^{(\ref{201906030335})}(\theta_1,\theta_2)
&=&
h^{-3/2}\int_\tjm^\tj\int_\tjm^t\int_\tjm^s
\big\{((H_xB)_xB)(Z_r,\theta_1^*,\theta_3^*)
\nn\\&&
-((H_xB)_xB)(Z_\tjm,\theta_1^*,\theta_3^*)\big\}
dw_rdw_sdt
\nn\\&&
+h^{-3/2}\int_\tjm^\tj\int_\tjm^t\int_\tjm^sL_{H_xB}(Z_r,\theta_1^*,\theta_2^*,\theta_3^*)drdw_sdt
\nn\\&&
+h^{-3/2}\int_\tjm^\tj\int_\tjm^t \big(L_H(Z_s,\theta_1,\theta_2,\theta_3^*)-L_H(Z_\tjm,\theta_1,\theta_2,\theta_3^*)\big)dsdt
\eea
with 
\beas 
L_{H_xB}(z,\theta_1,\theta_2,\theta_3)
&=&
(H_xB)_x(z,\theta_1,\theta_3)[A(z,\theta_2)]+\half(H_xB)_{xx}(z,\theta_1,\theta_3)[C(z,\theta_1)]
\\&&
+(H_xB)_y(z,\theta_1,\theta_3)[H(z,\theta_3)],
\eeas
and
\bea\label{201906030336}
r_j^{(\ref{201906030336})}(\theta_1,\theta_2)
&=&
h^{-3/2}\int_\tjm^\tj\int_\tjm^t \big(L_H(Z_s,\theta_1^*,\theta_2^*,\theta_3^*)-L_H(Z_s,\theta_1,\theta_2,\theta_3^*)\big)dsdt.
\eea
Moreover, 
\bea\label{201906030446}
\sup_n\sup_j\bigg\|
\sup_{(\theta_1,\theta_2)\in\overline{\Theta}_1\times\overline{\Theta}_2}
\big|r_j^{(\ref{201906030335})}(\theta_1,\theta_2)\big|\bigg\|_p
&=&
O(h)
\eea
for every $p>1$, and 
{\cred 
\bea\label{201906030450}
\big|r_j^{(\ref{201906030336})}(\theta_1,\theta_2)\big|
&\leq&
h^{1/2}r_{n,j}^{(\ref{201906030451})}\big\{\big|\theta_1-\theta_1^*\big|
+\big|\theta_2-\theta_2^*\big|\big\}
\nn\\&&
\eea
for all $(\theta_1,\theta_2)\in\overline{\Theta}_1\times\overline{\Theta}_2$
}
with 
some random variables $r_{n,j}^{(\ref{201906030451})}$ satisfying 
\bea\label{201906030451}
\sup_n\sup_j\big\|r_{n,j}^{(\ref{201906030451})}\big\|_p
&<&
\infty
\eea
for every $p>1$. 

\bd\im[(b)] 
Suppose that $[A1]$ with $(i_A,j_A,i_B,j_B,i_H,j_H)=({\colorg1,1,2},1,2,{\colorg0})
$ 
and $[A2]$ $(i)$ are satisfied. 
Then there exist random variables $r_{n,j}^{(\ref{201906031442})}$ and a number $\rho$ such that 
\beas 
\sup_{\theta_3\in\overline{\Theta}_3}\big|\cald_j(\theta_1,\theta_2,\theta_3)-\cald_j(\theta_1^*,\theta_2^*,\theta_3)\big|
&\leq&
h^{1/2}
r_{n,j}^{(\ref{201906031442})}\big\{\big|\theta_1-\theta_1^*\big|+\big|\theta_2-\theta_2^*\big|\big\}
\eeas
for all $(\theta_1,\theta_2)\in B((\theta_1^*,\theta_2^*),\rho)$ 
and that 
\bea\label{201906031442}
\sup_n\sup_j\big\|r_{n,j}^{(\ref{201906031442})}\big\|_p
&<&
\infty
\eea
for every $p>1$. 
\ed
\begin{en-text}
\im[(c)] (erase this later)
Suppose that $[A1]$ with $(i_A,j_A,i_B,j_B,i_H,j_H)=(0,0,0,0,2,0)$ and $[A2]$ $(i)$ are satisfied. 
Then 
\beas 
\sup_n\sup_j\bigg\|\sup_{(\theta_1,\theta_2,\theta_3)\in\overline{\Theta}_1\times\overline{\Theta}_2\times\overline{\Theta}_3}
\big|\cald_j(\theta_1,\theta_2,\theta_3)-\cald_j(\theta_1^*,\theta_2^*,\theta_3)\big|\bigg\|_p
&=&
O(h^{1/2})
\eeas
for every $p>1$. 
\end{en-text}
\end{lemma}
\proof 
By (\ref{201905291809}), we have 
\bea\label{201906030250}
\Delta_jY -hG_n(Z_\tjm,\theta_1,\theta_2,\theta_3^*)
&=&
\xi_j^{(\ref{201906030334})}+\rho_j(\theta_1,\theta_2,\theta_3^*)
\eea
and 
\beas
\rho_j(\theta_1,\theta_2,\theta_3^*)
&=&
\int_\tjm^\tj\int_\tjm^t \big\{H_x(Z_s,\theta_3^*)B(Z_s,\theta_1^*)-H_x(Z_\tjm,\theta_3^*)B(Z_\tjm,\theta_1^*)\big\}dw_sdt
\\&&
+\int_\tjm^\tj\int_\tjm^t \big(L_H(Z_s,\theta_1^*,\theta_2^*,\theta_3^*)-L_H(Z_\tjm,\theta_1,\theta_2,\theta_3^*)\big)dsdt.
\eeas
\begin{en-text}
\\&=&
((H_xB)_xB)(Z_\tjm,\theta_1^*,\theta_3^*)\int_\tjm^\tj\int_\tjm^t\int_\tjm^sdw_rdw_sdt
+h^{3/2}\tilde{\rho}_j+h^{3/2}\hat{\rho}_j
\eeas
where 
$\sup_n\sup_j\|\tilde{\rho}_j\|_p=O(h)$ and 
\beas
\hat{\rho}_j
&=&
h^{-3/2}
\int_\tjm^\tj\int_\tjm^t \big(L_H(Z_s,\theta_1^*,\theta_2^*,\theta_3^*)-L_H({\colorr Z_s},\hat{\theta}_1,\theta_2^*,\theta_3^*)\big)dsdt.
\eeas

$|\hat{\rho}_j|\leq h^{1/2}\hat{r}_j|\hat{\theta}_1-\theta_1^*|$ with 
$\sup_n\sup_j\|\hat{r}_j\|_p<\infty$
\end{en-text}
Then the decomposition (\ref{201906030236}) is obvious. 
The first and third terms on the right-hand side of (\ref{201906030335}) can be estimated with 
Taylor's formula and (\ref{201905280550}), and the second term is easy to estimate. 
Thus, we obtain (\ref{201906030446}). 
Since $\partial_{(\theta_1,\theta_2)}L_H(z,\theta_1,\theta_2,\theta_3^*)$ is bound by a polynomial in $z$ uniformly in ${\cred (\theta_1,\theta_2)}$, 
there exist random variables $r_{n,j}^{(\ref{201906030451})}$ that satisfy (\ref{201906030450}) and 
(\ref{201906030451}). 
{\cred $[$ 
First show (\ref{201906030450}) on the set 
$\{|(\theta_1,\theta_2)-(\theta_1^*,\theta_2^*)|<r\}$, next see this estimate is valid 
on $\big(\overline{\Theta}_1\times\overline{\Theta}_2\big)\setminus\{|(\theta_1,\theta_2)-(\theta_1^*,\theta_2^*)|<r\}$ by redefining $r_{n,j}^{(\ref{201906030451})}$ if necessary. 
}
We obtained (a). 
The assertion (b) is easy to verify with (\ref{201905291809}), (\ref{201906021454}) 
and Lemma \ref{201906031518}. 
\qed\halflineskip

\begin{lemma}\label{201906040435}
Suppose that $[A1]$ with $(i_A,j_A,i_B,j_B,i_H,j_H)=(0,0,0,0,2,1)$ and $[A2]$ $(i)$ are satisfied. 
Then 
\beas 
\sup_{(\theta_1,\theta_2)\in\overline{\Theta}_1\times\overline{\Theta}_2}
\big|\cald_j(\theta_1,\theta_2,\theta_3)-\cald_j(\theta_1,\theta_2,\theta_3')\big|
&\leq& 
h^{-1/2}r_{n,j}^{(\ref{201906040420})}\big|\theta_3-\theta_3'\big|
\quad({\cred \theta_3,\theta_3'\in\overline{\Theta}_3})
\eeas
for 
some random variables {\cred $r_{n,j}^{(\ref{201906040420})}$} such that 
\bea\label{201906040420}
\sup_n\sup_j\big\|r_{n,j}^{(\ref{201906040420})}\big\|_p
&<&
\infty
\eea
for every $p>1$. 
\end{lemma}
\proof
\beas &&
\cald_j(\theta_1,\theta_2,\theta_3)-\cald_j(\theta_1,\theta_2,\theta_3')
\\&=&
\left(\begin{array}{c}0\\ 
h^{-1/2}\big(H(Z_\tjm,\theta_3')-H(Z_\tjm,\theta_3)\big)
+\frac{h^{\colorr{1/2}}}{2}\big(L_H(Z_\tjm,\theta_1,\theta_2,\theta_3')-L_H(Z_\tjm,\theta_1,\theta_2,\theta_3)\big)
\end{array}\right)
\eeas
Therefore the lemma is obvious. 
{\cred 
Apply the Taylor formula for the argument $\theta_3$ if $\theta_3$ and $\theta_3'$ are close, otherwise 
and if necessary, 
redifine $r_{n,j}^{(\ref{201906040420})}$. }
\qed\halflineskip

\section{An adaptive estimator for $\theta_3$
}\label{201906041935}
We will work with some initial estimators $\hat{\theta}_1^0$ for $\theta_1^0$ 
and $\hat{\theta}_2^0$ for $\theta_2$. 
The following standard convergence rates, in part or fully, will be assumed for these estimators:
\bd
\im[[A4\!\!]] 
{\bf (i)} 
$\ds \hat{\theta}_1^0-\theta_1^* \yeq O_p(n^{-1/2})$ as $n\to\infty$
\bd
\im[(ii)] 
$\ds \hat{\theta}_2^0-\theta_2^*\yeq O_p(n^{-1/2}h^{-1/2})$  as $n\to\infty$
\ed\ed

Sections \ref{201905291607} and \ref{201906041938} 
recall certain standard estimators for $\theta_1$ and $\theta_2$, respectively.  
The expansions (\ref{201906021458}) and (\ref{201905291809}) 
with Lemma {\colorg\ref{201906031516}} 
suggest 
two approaches for estimating $\theta_3$. 
The first approach is based on the likelihood of $h^{-3/2}\big(\Delta_jY-{\colorr h}G_n(Z_\tjm,\theta_1,\theta_2,\theta_3)\big)$ only. 
The second one uses the likelihood corresponding to $\cald_j(\theta_1,\theta_2,\theta_3)$. 
However, it is possible to show that the first approach gives less optimal asymptotic variance; 
see Remark \ref{201906181743}. 
So, we will treat the second approach here. 

\subsection{{\coloro Adaptive quasi-likelihood function for $\theta_3$}}
Let 
\beas
S(z,\theta_1,\theta_3)&=&
\left(\begin{array}{ccc}
C(z,\theta_1)&& 2^{-1}C(z,\theta_1)H_x(z,\theta_3)^\star
\\
2^{-1}H_x(z,\theta_3)C(z,\theta_1)&& 3^{-1}H_x(z,\theta_3)C(z,\theta_1)H_x(z,\theta_3)^\star
\end{array}\right)
\eeas
Then 
\bea\label{201906040958} &&
S(z,\theta_1,\theta_3)^{-1}
\nn\\&=&
\left(\begin{array}{ccc}
C(z,\theta_1)^{-1}+3H_x(z,\theta_3)^\star V(z,\theta_1,\theta_3)^{-1}H_x(z,\theta_3)
&&-6H_x(z,\theta_3)^\star V(z,\theta_1,\theta_3)^{-1}\y
-6V(z,\theta_1,\theta_3)^{-1}H_x(z,\theta_3)&&12V(z,\theta_1,\theta_3)^{-1}
\end{array}\right). 
\nn\\&&
\eea
Recall that 
\beas 
V(z,\theta_1,\theta_3) &=& H_x(z,\theta_3)C(z,\theta_1)H_x(z,\theta_3)^\star.
\eeas
Let 
\beas
\hat{S}(z,\theta_3)&=&S(z,\hat{\theta}_1^0,\theta_3).
\eeas

We define a log quasi-likelihood function by 
\bea\label{201906041919} 
\bbH^{(3)}_n(\theta_3)
 &=& 
 -\half\sum_{j=1}^n \bigg\{
 \hat{S}(Z_\tjm,\theta_3)^{-1}\big[\cald_j(\hat{\theta}_1^0,\hat{\theta}_2^0,\theta_3)^{\otimes2}\big]
 +\log\det \hat{S}(Z_\tjm,\theta_3)\bigg\}.
\eea
Let $\hat{\theta}_3^0$ be a quasi-maximum likelihood estimator (QMLE) for $\theta_3$
 for $\bbH^{(3)}_n$, that is, $\hat{\theta}_3^0$ is a 
 $\overline{\Theta}_3$-valued  measurable mapping satisfying 
 \beas 
 \bbH^{(3)}_n(\hat{\theta}_3^0)
 &=&
 \max_{\theta_3\in\overline{\Theta}_3}\bbH^{(3)}_n(\theta_3). 
 \eeas
 The QMLE 
 $\hat{\theta}_3^0$ for $ \bbH^{(3)}_n$ 
 depends on $n$ as it does on the data $(Z_\tj)_{j=0,1,...,n}$; 
 $\hat{\theta}_1^0$ in the function $\hat{S}$ also depends on $(Z_\tj)_{j=0,1,...,n}$.

We introduce the following random fields 
depending on $n$. 
\beas
\Psi_1(\theta_1,\theta_2,\theta_3,\theta_1',\theta_2',\theta_3')
&=&
\sum_{j=1}^nS(Z_\tjm,\theta_1,\theta_3)^{-1}
\left[\cald_j(\theta_1',\theta_2',\theta_3'),\>
\left(\begin{array}{c}0\\
2^{-1}\partial_1 L_H(Z_\tjm,\theta_1,\theta_2,\theta_3)
\end{array}\right)\right]
\\&=&
\sum_{j=1}^nS(Z_\tjm,\theta_1,\theta_3)^{-1}
\left[\cald_j(\theta_1',\theta_2',\theta_3'),\>
\left(\begin{array}{c}0\\
2^{-1}H_{xx}(z,\theta_3)[\partial_1C(Z_\tjm,\theta_1)]
\end{array}\right)\right],
\eeas
\beas
\Psi_2(\theta_1,\theta_2,\theta_3,\theta_1',\theta_2',\theta_3')
&=&
\sum_{j=1}^nS(Z_\tjm,\theta_1,\theta_3)^{-1}
\left[\cald_j(\theta_1',\theta_2',\theta_3'),\>
\left(\begin{array}{c}\partial_2A(Z_\tjm,\theta_1,\theta_2)\\
2^{-1}\partial_2 L_H(Z_\tjm,\theta_1,\theta_2,\theta_3)
\end{array}\right)\right]
\\&=&
\sum_{j=1}^nS(Z_\tjm,\theta_1,\theta_3)^{-1}
\left[\cald_j(\theta_1',\theta_2',\theta_3'),\>
\left(\begin{array}{c}\partial_2A(Z_\tjm,\theta_1,\theta_2)\\
2^{-1}H_x(z,\theta_3)[\partial_2A(Z_\tjm,\theta_2)]
\end{array}\right)\right]
\eeas
\beas
\widetilde{\Psi}_2(\theta_1,\theta_2,\theta_3,\theta_1',\theta_2',\theta_3')
&=&
\sum_{j=1}^nS(Z_\tjm,\theta_1,\theta_3)^{-1}
\left[\widetilde{\cald}_j(\theta_1',\theta_2',\theta_3'),\>
\left(\begin{array}{c}\partial_2A(Z_\tjm,\theta_2)\\
2^{-1}\partial_2 L_H(Z_\tjm,\theta_1,\theta_2,\theta_3)
\end{array}\right)\right]
\\&=&
\sum_{j=1}^nS(Z_\tjm,\theta_1,\theta_3)^{-1}
\left[\widetilde{\cald}_j(\theta_1',\theta_2',\theta_3'),\>
\left(\begin{array}{c}\partial_2A(Z_\tjm,\theta_2)\\
2^{-1}H_x(z,\theta_3)[\partial_2A(Z_\tjm,\theta_2)]
\end{array}\right)\right],
\eeas
where 
\beas 
\widetilde{\cald}_j(\theta_1^*,\theta_2^*,\theta_3^*)
&=&
\left(\begin{array}{c}
{\cred\xi^{(\ref{201906030126})}_j+\xi^{(\ref{201906030150})}_j}\y
h^{-3/2}\big(\xi_j^{(\ref{201906030334})}+\xi_j^{(\ref{201906030345})}\big)
\end{array}\right),
\eeas
and
\beas
\Psi_3(\theta_1,\theta_2,\theta_3,\theta_1',\theta_2',\theta_3')
&=&
\sum_{j=1}^nS(Z_\tjm,\theta_1,\theta_3)^{-1}
\bigg[\cald_j(\theta_1',\theta_2',\theta_3')
\\&&\qquad\otimes
\left(\begin{array}{c}0\\
\partial_3 H(Z_\tjm,\theta_3)
+2^{-1}h\partial_3 L_H(Z_\tjm,\theta_1,\theta_2,\theta_3)
\end{array}\right)\bigg]
\eeas
\beas
\widetilde{\Psi}_3(\theta_1,\theta_2,\theta_3,\theta_1',\theta_2',\theta_3')
&=&
\sum_{j=1}^nS(Z_\tjm,\theta_1,\theta_3)^{-1}
\bigg[\widetilde{\cald}_j(\theta_1',\theta_2',\theta_3')
\\&&\qquad\otimes
\left(\begin{array}{c}0\\
\partial_3 H(Z_\tjm,\theta_3)
+2^{-1}h\partial_3 L_H(Z_\tjm,\theta_1,\theta_2,\theta_3)
\end{array}\right)\bigg]
\eeas
\beas
\Psi_{3,1}(\theta_1,\theta_3,\theta_1',\theta_2',\theta_3')
&=&
\sum_{j=1}^nS(Z_\tjm,\theta_1,\theta_3)^{-1}
\bigg[\cald_j(\theta_1',\theta_2',\theta_3'),\>
\left(\begin{array}{c}0\\
\partial_3 H(Z_\tjm,\theta_3)
\end{array}\right)\bigg]
\eeas
\beas
\widetilde{\Psi}_{3,1}(\theta_1,\theta_3,\theta_1',\theta_2',\theta_3')
&=&
\sum_{j=1}^nS(Z_\tjm,\theta_1,\theta_3)^{-1}
\bigg[\widetilde{\cald}_j(\theta_1',\theta_2',\theta_3'),\>
\left(\begin{array}{c}0\\
\partial_3 H(Z_\tjm,\theta_3)
\end{array}\right)\bigg]
\eeas
\beas
\Psi_{3,2}(\theta_1,\theta_2,\theta_3,\theta_1',\theta_2',\theta_3')
&=&
\sum_{j=1}^nS(Z_\tjm,\theta_1,\theta_3)^{-1}
\bigg[\cald_j(\theta_1',\theta_2',\theta_3'),\>
\left(\begin{array}{c}0\\
2^{-1}h\partial_3 L_H(Z_\tjm,\theta_1,\theta_2,\theta_3)
\end{array}\right)\bigg]
\eeas
\beas
\Psi_{3,3}(\theta_1,\theta_3,\theta_1',\theta_2',\theta_3')
&=&
\half\sum_{j=1}^n\big(S^{-1}(\partial_3S)S^{-1}\big)(Z_\tjm,\theta_1,\theta_3)
\big[\cald_j(\theta_1',\theta_2',\theta_3')^{\otimes2}-S(Z_\tjm,\theta_1,\theta_3)\big]
\eeas
\beas
\Psi_{33,1}(\theta_1,\theta_2,\theta_3)
&=&
-\sum_{j=1}^nS(Z_\tjm,\theta_1,\theta_3)^{-1}
\left[
\left(\begin{array}{c}0\\
\partial_3H(Z_\tjm,\theta_3)+2^{-1}h\partial_3L_H(Z_\tjm,\theta_1,\theta_2,\theta_3)
\end{array}\right)^{\otimes2}
\right]
\eeas
\beas
\Psi_{33,2}(\theta_1,\theta_2,\theta_3,\theta_1',\theta_2',\theta_3')
&=&
\sum_{j=1}^nS(Z_\tjm,\theta_1,\theta_3)^{-1}
\bigg[\cald_j(\theta_1',\theta_2',\theta_3')
\\&&\qquad\otimes
\left(\begin{array}{c}0\\
\partial_3^2 H(Z_\tjm,\theta_3)
+2^{-1}h\partial_3^2 L_H(Z_\tjm,\theta_1,\theta_2,\theta_3)
\end{array}\right)\bigg]
\eeas
\beas
\Psi_{33,3}(\theta_1,\theta_3)
&=&
-\half\sum_{j=1}^n
\big\{\big(S^{-1}(\partial_3S)S^{-1}\big)(Z_\tjm,\theta_1,\theta_3)
\big[\partial_3S(Z_\tjm,\theta_1,\theta_3)\big]\big\}
\eeas
\beas
\Psi_{33,4}(\theta_1,\theta_2,\theta_3,\theta_1',\theta_2',\theta_3')
&=&
- {\colorr2}\sum_{j=1}^n
S^{-1}(\partial_3S)S^{-1}(Z_\tjm,\theta_1,\theta_3)
\left[\cald_j(\theta_1',\theta_2',\theta_3')\right.
\\&&\hspace{3cm}\otimes
\left.
\left(\begin{array}{c}0\\
\partial_3H(Z_\tjm,\theta_3)+2^{-1}h\partial_3L_H(Z_\tjm,\theta_1,\theta_2,\theta_3)
\end{array}\right)
\right]
\eeas
\beas
\Psi_{33,5}(\theta_1,\theta_3,\theta_1',\theta_2',\theta_3')
&=&
\half \sum_{j=1}^n\partial_3\big\{
\big(S^{-1}(\partial_3S)S^{-1}\big)(Z_\tjm,\theta_1,\theta_3)\big\}
\big[\cald_j(\theta_1',\theta_2',\theta_3')^{\otimes2}-S(Z_\tjm,\theta_1,\theta_3)\big].
\eeas
%


\subsection{Consistency of $\hat{\theta}_3^0$}

\begin{lemma}\label{201906021729}
Suppose that $[A1]$ with $(i_A,j_A,i_B,j_B,i_H,j_H)=(0,0,0,1,1,1)$ and $[A2]$ $(i)$, $(iii)$ and $(iv)$ are fulfilled. 
Then 
\beas
\sup_{t\in\bbR_+}\bigg\|\sup_{(\theta_1,\theta_3)\in\overline{\Theta}_1\times\overline{\Theta}_3}
\big\{\big|S(Z_t,\theta_1,\theta_3)\big|+
\det S(Z_t,\theta_1,\theta_3)^{-1}+\big|S(Z_t,\theta_1,\theta_3)^{-1}\big|
\big\}\bigg\|_p
&<&
\infty
\eeas
for every $p>1$ 
\end{lemma}
\proof
{\colorg By $[A2]$ (iii) and (iv),}
$\det S(Z_\tjm,\theta_1,\theta_3)^{-1}$ as well as $S(Z_\tjm,\theta_1,\theta_3)$ 
is continuous on $\overline{\Theta}_1\times\overline{\Theta}_3$ a.s., 
and continuously differentiable on $\Theta_1\times\Theta_3$. 
Moreover we see
\beas 
\sup_{t\in\bbR_+}\sum_{i=0,1}\sup_{(\theta_1,\theta_3)\in\Theta_1\times\Theta_3}
\big\|\partial_{(\theta_1,\theta_3)}^i\big(\det S(Z_t,\theta_1,\theta_3)^{-1}\big)\big\|_p
&<&
\infty
\eeas
for every $p>1$ {\colorg from (\ref{201906040958})}. 
This implies that 
\beas
\sup_{t\in\bbR_+}\bigg\|\sup_{(\theta_1,\theta_3)\in\Theta_1\times\Theta_3}\big(\det S(Z_t,\theta_1,\theta_3)^{-1}\big)
\bigg\|_p
&<&
\infty
\eeas
for every $p>1$ 
by Sobolev's inequality. The inequality 
\beas
\sup_{t\in\bbR_+}\bigg\|\sup_{(\theta_1,\theta_3)\in\Theta_1\times\Theta_3}\big|S(Z_t,\theta_1,\theta_3)
\big|\bigg\|_p
&<&
\infty
\eeas
for every $p>1$ 
is rather easy to show. 
\qed\halflineskip

Let 
\beas 
\bbY^{(3)}_n(\theta_3)
&=&
n^{-1}h\big\{\bbH^{(3)}_n(\theta_3)-\bbH^{(3)}_n(\theta_3^*)\big\}.
\eeas
\begin{theorem}\label{201906021825}
Suppose that $[A1]$ with $(i_A,j_A,i_B,j_B,i_H,j_H)=(1,{\colorg1},{\colorg2},{\colorr1},3,1)
$ and $[A2]$ are satisfied. 
Then   
\bea\label{201906021845}
\sup_{\theta_3\in\overline{\Theta}_3}
\big|\bbY_n^{(3)}(\theta_3)
-\bbY^{(3)}(\theta_3)
\big|
&\to^p&
0
\eea
as $n\to\infty$, if $\hat{\theta}_1^0\to^p\theta_1^*$ and $\hat{\theta}_2^0\to^p\theta_2^*$. 
Moreover, $\hat{\theta}_3^0\to^p\theta_3^*$ 
if $[A3]$ $(iii)$ is additionally satisfied. 
\end{theorem}
\begin{en-text}
{\cred {\bf Modify this paragraph later.} 
Remark that if we assume $[A3]$ $(i)$ in the first part of Theorem \ref{201906021825} 
(in place of the convergence of $\hat{\theta}_1^0$), then 
the convergence $\hat{\theta}_1^0\to^p\theta_1^*$ 
for the QMLE $\hat{\theta}_1$ 
follows from 
Theorem \ref{201905291601} (a). \koko about $\hat{\theta}_2^0$\y
}
\end{en-text}
\noindent
{\it Proof of Theorem \ref{201906021825}.} 
We have 
\beas 
\bbY_n^{(3)}(\theta_3)
&=& 
n^{-1}h^{1/2}\sum_{j=1}^n 
\hat{S}(Z_\tjm,\theta_3)^{-1}\big[h^{1/2}\delta_j(\hat{\theta}_1^0,\hat{\theta}_2^0,\theta_3),
\cald_j(\hat{\theta}_1^0,\hat{\theta}_2^0,\theta_3^*)\big]
\\&&
-\frac{1}{2n}\sum_{j=1}^n  \hat{S}(Z_\tjm,\theta_3)^{-1}
\big[\big(h^{1/2}\delta_j(\hat{\theta}_1^0,\hat{\theta}_2^0,\theta_3)\big)^{\otimes2}\big]
+n^{-1}hR^{(\ref{201906021417})}_n(\theta_3)
\eeas
where 
\beas 
\delta_j(\theta_1,\theta_2,\theta_3) 
&=& 
-\cald_j(\theta_1,\theta_2,\theta_3)+\cald_j(\theta_1,\theta_2,\theta_3^*)
\eeas
and
\bea\label{201906021417}
R^{(\ref{201906021417})}_n(\theta_3)
&=&
-\half\sum_{j=1}^n\big(
\hat{S}(Z_\tjm,\theta_3)^{-1}-\hat{S}(Z_\tjm,\theta_3^*)^{-1}\big)
\big[\cald_j(\hat{\theta}_1^0,\hat{\theta}_2^0,\theta_3^*)^{\otimes2}\big]
\nn\\&&
-\half\sum_{j=1}^n \log\frac{\det \hat{S}(Z_\tjm,\theta_3)}{\det \hat{S}(Z_\tjm,\theta_3^*)}
\eea

By Lemma \ref{201906021438} (b), Lemma \ref{201906031516} (b) 
and Lemma \ref{201906021729}, we obtain  
\beas\label{201906021619} 
n^{-1}h\sup_{\theta_3\in\Theta_3}\big|R^{(\ref{201906021417})}_n(\theta_3)
\big|
&=&
{\colorr O_p(h).}
\eeas
By definition, 
\beas 
h^{1/2}\delta_j(\hat{\theta}_1^0,\hat{\theta}_2^0,\theta_3) 
&=& 
\left(\begin{array}{c}
0
\y
H(Z_\tjm,\theta_3)-H(Z_\tjm,\theta_3^*)
+\frac{h}{2}\big(L_H(Z_\tjm,\hat{\theta}_1^{{\cred0}},\hat{\theta}_2^0,\theta_3)-L_H(Z_\tjm,\hat{\theta}_1^{{\cred0}},\hat{\theta}_2^{{\cred0}},\theta_3^*)\big)
\end{array}\right)
\eeas
Since the functions $A(z,\theta_2)$, $H(z,\theta_3)$ and $L_H(z,\theta_1,\theta_2,\theta_3)$ are 
dominated by a polynomial in $z$ uniformly in $\theta$, 
by using the above formula, 
it is easy to show 
\bea\label{201906021842} 
\sup_{\theta_3\in\overline{\Theta}_3}
\left|\bbY^{(3)}_n(\theta_3)
-\bbY^{(\ref{201906021814})}_n(\hat{\theta}_1^0,\theta_3)
\right|
&=&
{\colorr O_p(h^{1/2})}
\eea
for 
\bea\label{201906021814}
\bbY^{(\ref{201906021814})}_n(\theta_1,\theta_3)
&=&
-\frac{1}{2n}\sum_{j=1}^nS(Z_\tjm,\theta_1,\theta_3)^{-1}
\left[\left(\begin{array}{c}0\\
H(Z_\tjm\,\theta_3)-H(Z_\tjm\,\theta_3^*)
\end{array}\right)^{\otimes2}\right].
\eea
The derivative $\partial_1S_x(z,\theta_1,\theta_3)$ is dominated by a polynomial in $z$ 
uniformly in $\theta$. Therefore 
\bea\label{201906021843}
\sup_{\theta_3\in\overline{\Theta}_3}\big|
\bbY^{(\ref{201906021814})}_n(\hat{\theta}_1^0,\theta_3)-\bbY^{(\ref{201906021814})}_n(\theta_1^*,\theta_3)\big|
\to^p0. 
\eea
Finally, the estimate (\ref{201905280550}) gives 
\bea\label{201906021844} 
\sup_{\theta_3\in\overline{\Theta}_3}\left|
\bbY^{(\ref{201906021814})}_n(\theta_1^*,\theta_3)
+\frac{1}{2nh}\int_0^{nh}S(Z_t,\theta_1^*,\theta_3)^{-1}
\left[\left(\begin{array}{c}0\\
H(Z_t,\theta_3)-H(Z_t,\theta_3^*)
\end{array}\right)^{\otimes2}\right]dt
\right|
\to^p0. 
\eea
Now 
(\ref{201906021845}) follows from (\ref{201906021842}), (\ref{201906021843}), (\ref{201906021844}) 
and $[A2]$ (ii) since ${\colorg\partial_3}^iH(z,\theta_1,\theta_3)$ $(i=0,1)$ are dominated by a polynomial in $z$ 
uniformly in $\theta_3$. 
Then the convergence $\hat{\theta}_3^0\to^p\theta_3$ as $n\to\infty$ is obvious under Condition $[A3]$ (iii). 
\qed\halflineskip

\subsection{Asymptotic normality of $\hat{\theta}_3^0$} 
\begin{en-text}
We will prove 
\beas 
\hat{\theta}_3^0-\theta_3^* &=& \xout{{\colorr O_p(n^{-1/2}h^{1/2})}}{\colorr O_p(n^{-1/2})}
\eeas
as a preliminary estimate of the error. 
\end{en-text}

Let 
\beas 
M^{(3)}_n
&=&
n^{-1/2}\sum_{j=1}^nS(Z_\tjm,\theta_1^*,\theta_3^*)^{-1}
\left[
\left(\begin{array}{c}
h^{-1/2}B(Z_\tjm,\theta_2^*)\Delta_jw\y
h^{-3/2}\kappa(Z_\tjm,\theta_1^*,\theta_3^*)\zeta_j
\end{array}\right),\>
\left(\begin{array}{c}0\\
\partial_3H(Z_\tjm,\theta_3^*)
\end{array}\right)
\right]. 
\eeas
Let 
\bea\label{201906190133} 
\Gamma_{33}
&=&
\int S(z,\theta_1^*,\theta_3^*)^{-1}
\left[
\left(\begin{array}{c}0\\
\partial_3H(z,\theta_3^*)
\end{array}\right)^{\otimes2}
\right]
\nu(dz)
\nn\\&=&
\int 12V(z,\theta_1^*,\theta_3^*)^{-1}\big[\big(\partial_3H(z,\theta_3^{{\colorg*}})\big)^{\otimes2}\big]\nu(dz)
\nn\\&=&
\int12\partial_3H(z,\theta_3^*)^\star V(z,\theta_1^*,\theta_3^*)^{-1}\partial_3H(z,\theta_3^*)\nu(dz).
\eea

\begin{lemma}\label{201906031714}
Suppose that $[A1]$ with $(i_A,j_A,i_B,j_B,i_H,j_H)=(1,1,2,1,3,1)$, $[A2]$, 
{\colorr $[A3]$ $(iii)$} and $[A4]$ are satisfied. 
Then 
{\colorr 
\beas 
n^{-1/2}h^{1/2}\>\partial_3\bbH_n^{(3)}(\theta_3^*)
-M^{(3)}_n
&=&
o_p(1)
\eeas
}
as $n\to\infty$. 
\end{lemma}
\proof 
From $[A3]$ (iii), $\Gamma_{33}$ is non-singular. 
From (\ref{201906041919}) and (\ref{201906030041}), we have 
\bea\label{201906030105}
n^{-1/2}h^{1/2}\>\partial_3\bbH_n^{(3)}(\theta_3^*)
&=&
R^{(\ref{201906030100})}_n(\hat{\theta}_1^0,\hat{\theta}_2^0)+R^{(\ref{201906022139})}_n(\hat{\theta}_1^0,\hat{\theta}_2^0)
+R^{(\ref{201906030102})}_n(\hat{\theta}_1^0,\hat{\theta}_2^0)
\eea
where 
\bea\label{201906030100}
R^{(\ref{201906030100})}_n(\hat{\theta}_1^0,\hat{\theta}_2^0)
&=&
n^{-1/2}\Psi_{3,1}(\hat{\theta}_1^0,\theta_3^*,\hat{\theta}_1^0,\hat{\theta}_2^0,\theta_3^*),
\eea
\bea\label{201906022139}
R^{(\ref{201906022139})}_n(\hat{\theta}_1^0,\hat{\theta}_2^0)
&=&
n^{-1/2}\Psi_{3,2}(\hat{\theta}_1^0,\hat{\theta}_2^0,\theta_3^*,\hat{\theta}_1^0,\hat{\theta}_2^0,\theta_3^*)
\nn\\&&
\eea
and 
\bea\label{201906030102}
R^{(\ref{201906030102})}_n(\hat{\theta}_1^0,\hat{\theta}_2^0)
&=&
n^{-1/2}h^{1/2}\Psi_{3,3}(\hat{\theta}_1^0,\theta_3^*,\hat{\theta}_1^0,\hat{\theta}_2^0,
\theta_3^*).
\eea

We have 
\beas 
\cald_j({\colorr\hat{\theta}_1^0,\hat{\theta}_2^0},\theta_3^*)
-\cald_j({\colorr\hat{\theta}_1^0,\theta_2^*},\theta_3^*)
&=& 
-h^{1/2}
\left[\begin{array}{c}
\big(A(Z_\tjm,\hat{\theta}_2^0)-A(Z_\tjm,\theta_2^*)\big)\y
2^{-1}H_x(Z_\tjm,\theta_3^*)\big[A(Z_\tjm,\hat{\theta}_2^0)-A(Z_\tjm,\theta_2^*)\big]
\end{array}\right], 
\eeas
and so only by algebraic computation we obtain 
\bea\label{201906091302} 
\hat{S}(Z_\tjm,\theta_3^*)^{-1}
\left[\cald_j({\colorr\hat{\theta}_1^0,\hat{\theta}_2^0},\theta_3^*)-\cald_j({\colorr\hat{\theta}_1^0,\theta_2^*},\theta_3^*),\>
\left(\begin{array}{c}0\\
\partial_3H(Z_\tjm,\theta_3^*)
\end{array}\right)
\right]
&=&
0. 
\eea
Applying Lemma 
{\colorg \ref{201906031516} (b)}  
under $[A4]$, and next using 
the results in Lemmas \ref{201906031518} and \ref{201906031516}, 
we see 
\bea\label{201906031505}
R^{(\ref{201906030100})}_n(\hat{\theta}_1^0,\hat{\theta}_2^0)
&=&
{\colorg n^{-1/2}}
\Psi_{3,1}(\hat{\theta}_1^0,\theta_3^*,\theta_1^*,\theta_2^*,\theta_3^*)
+O_p({\colorr h^{{\colorr1/2}}})
\nn\\&=&
{\colorg n^{-1/2}}
\tilde{\Psi}_{3,1}(\hat{\theta}_1^0,\theta_3^*,\theta_1^*,\theta_2^*,\theta_3^*)+o_p(1)
\nn\\&&
\eea
since $(nh^2)^{1/2}=o(1)$. 
{\colorr 
Consider the random field 
\bea\label{201906091137}
\Phi_n^{(\ref{201906091137})}(u_1)
&=&
{\colorg n^{-1/2}\big\{}
\widetilde{\Psi}_{3,1}(\theta_1^*+r_nu_1,\theta_3^*,\theta_1^*,\theta_2^*,\theta_3^*)
-\widetilde{\Psi}_{3,1}(\theta_1^*,\theta_3^*,\theta_1^*,\theta_2^*,\theta_3^*)
{\colorg\big\}}
\eea
on $\{u_1\in\bbR^{\sfp_1};\>|u_1|<1\}$ for any sequence of positive numbers $r_n\to0$,  
Sobolev's inequality gives 
\beas 
\sup_{u_1: |u_1|<1}|\Phi_n^{(\ref{201906091137})}(u_1)| 
&=&
o_p(1)
\eeas
{\colorg with the help of orthogonality.}
In particular, 
\bea\label{201906091150}
R^{(\ref{201906030100})}_n(\hat{\theta}_1^0,\hat{\theta}_2^0)
&=&
n^{-1/2}\widetilde{\Psi}_{3,1}(\theta_1^*,\theta_3^*,\theta_1^*,\theta_2^*,\theta_3^*)+o_p(1).
\nn\\&&
\eea
This implies 
\beas 
R^{(\ref{201906030100})}_n(\hat{\theta}_1^0,\hat{\theta}_2^0)
&=&
M_n^{(3)}+o_p(1). 
\eeas
}
Simpler is that 
$R^{(\ref{201906022139})}_n(\hat{\theta}_1^0,\hat{\theta}_2^0)={\colorr O_p(n^{1/2}h)}$. 
Similarly, 
\beas
R^{(\ref{201906030102})}_n(\hat{\theta}_1^0,\hat{\theta}_2^0)
&=&
n^{-1/2}h^{\colorg 1/2}\Psi_{3,3}(\theta_1^*,\theta_3^*,\theta_1^*,\theta_2^*,\theta_3^*)
+O_p(h^{1/2})
\\&=&
O_p({\colorr h^{1/2}}). 
\eeas
Thus, we obtained the result. 
\qed\halflineskip

\begin{en-text}
\beas
n^{-1/2}h\>\partial_3\bbH_n^{(2,3)}(\hat{\theta}_2^0,\theta_3^*)
&=&
O_p(h^{1/2})+O_p(n^{1/2}h\sqrt{h})
\\&=&
O_p(h^{1/2})
\eeas
where 
\bea
R^{(\ref{201906022139})}_n(\hat{\theta}_1,\hat{\theta}_2^0)
&=&
n^{-1/2}h^{1/2}\sum_{j=1}^n\hat{S}(Z_\tjm,\theta_3^*)^{-1}
\left[\cald_j(\hat{\theta}_2^0,\theta_3^*),\>
\left(\begin{array}{c}0\\
2^{-1}h\partial_3L_H(Z_\tjm,\hat{\theta}_1,\hat{\theta}_2^0,\theta_3^*)
\end{array}\right)
\right]
\nn\\&=&
n^{-1/2}h^{1/2}\sum_{j=1}^n\hat{S}(Z_\tjm,\theta_3^*)^{-1}
\left[\cald_j(\theta_2^*,\theta_3^*),\>
\left(\begin{array}{c}0\\
2^{-1}h\partial_3L_H(Z_\tjm,\hat{\theta}_1,\hat{\theta}_2^0,\theta_3^*)
\end{array}\right)
\right]
\nn\\&&
+O_p(n^{1/2}h^2)
\nn\\&=&
n^{-1/2}h^{1/2}\sum_{j=1}^n\hat{S}(Z_\tjm,\theta_3^*)^{-1}
\left[
\left[\begin{array}{c}
h^{-1/2}\int_\tjm^\tj B(Z_t,\theta_2^*)dw_t\y
h^{-3/2}K(Z_\tjm,\theta_1^*,\theta_3^*)\zeta_j
\end{array}
\right]\right.
\nn\\&&
\hspace{5cm}\otimes
\left.
\left[\begin{array}{c}0\\
2^{-1}h\partial_3L_H(Z_\tjm,\hat{\theta}_1,\hat{\theta}_2^0,\theta_3^*)
\end{array}\right]
\right]
+O_p(n^{1/2}h^2)
\nn\\&=&
O_p(n^{1/2}h\times\sqrt{h})\yeq o_p(\sqrt{h})
\eea
\end{en-text}

{\colorb In what follows, we quite often use the estimates in 
Lemma \ref{201906021729} without mentioning it explicitly. }

\begin{lemma}\label{201906031658}
Suppose that $[A1]$ with $(i_A,j_A,i_B,j_B,i_H,j_H)=(1,1,{\colorg 2},1,3,2)
$, $[A2]$ and $[A4]$ are satisfied. Then 
\beas 
\sup_{\theta_3\in B_n}\left|
n^{-1}h\>\partial_3^2\bbH_n^{(3)}(\theta_3)
+\Gamma_{33}
\right|
&\to^p&
0
\eeas
for any sequence of balls $B_n$ in $\bbR^{\sfp_3}$ shrinking to $\theta_3^*$. 
\end{lemma}
\proof 
\begin{en-text}
Since 
\beas &&
n^{-1}h\>\partial_3\bbH_n^{(3)}(\theta_3)
\\&=&
n^{-1}h^{1/2}\sum_{j=1}^n\hat{S}(Z_\tjm,\theta_3)^{-1}
\left[\cald_j(\hat{\theta}_1^0,\hat{\theta}_2^0,\theta_3),\>
\left[\begin{array}{c}0\\
\partial_3H(Z_\tjm,\theta_3)+2^{-1}h\partial_3L_H(Z_\tjm,\hat{\theta}_1^0,\hat{\theta}_2^0,\theta_3)
\end{array}\right]
\right]
\\&&
+\half n^{-1}h\sum_{j=1}^n
\big(\hat{S}^{-1}(\partial_3\hat{S})\hat{S}^{-1}\big)(Z_\tjm,\theta_3)
\big[\cald_j(\hat{\theta}_1^0,\hat{\theta}_2^0,\theta_3)^{\otimes2}-\hat{S}(Z_\tjm,\theta_3)\big],
\eeas
\end{en-text}
We have 
\beas 
n^{-1}h\>\partial_3^2\bbH_n^{(3)}(\theta_3)
&=&
n^{-1}\Psi_{33,1}(\hat{\theta}_1^0,\hat{\theta}_2^0,\theta_3)
\\&&
+n^{-1}h^{1/2}\Psi_{33,2}(\hat{\theta}_1^0,\hat{\theta}_2^0,\theta_3,\hat{\theta}_1^0,\hat{\theta}_2^0,\theta_3)
\\&&
+n^{-1}h\Psi_{33,3}(\hat{\theta}_1^0,\theta_3)
\\&&
+n^{-1}h^{1/2}\Psi_{33,4}(\hat{\theta}_1^0,\hat{\theta}_2^0,\theta_3,\hat{\theta}_1^0,\hat{\theta}_2^0,\theta_3)
\nn\\&&
+n^{-1}h\Psi_{33,5}(\hat{\theta}_1^0,\theta_3,\hat{\theta}_1^0,\hat{\theta}_2^0,\theta_3)
.
\eeas
For $\cald_j(\hat{\theta}_1^0,\hat{\theta}_2^0,\theta_3)$ in the above expression, 
we use Lemma \ref{201906031516} (b) to replace $\hat{\theta}_i^0$ by $\theta_i^*$ 
for $i=1,2$, and Lemma \ref{201906040435} to replace 
$\theta_3\in B_n$ by $\theta_3^*$ with an error uniform in $\theta_3\in B_n$. 
Next we use Lemma \ref{201906021438} (b). 
Then 
\beas 
n^{-1}h\>\partial_3^2\bbH_n^{(3)}(\theta_3)
&=&
-n^{-1}\sum_{j=1}^n\hat{S}(Z_\tjm,\theta_3)^{-1}
\left[
\left(\begin{array}{c}0\\
\partial_3H(Z_\tjm,\theta_3)
\end{array}\right)^{\otimes2}
\right]
+r_n^{(\ref{201906031646})}(\theta_3)
\eeas
where 
\bea\label{201906031646}
\sup_{\theta_3\in\overline{\Theta}_3}\big|r_n^{(\ref{201906031646})}(\theta_3)\big|
&=&
o_p(1). 
\eea
Now we obtain the result by using $[A2]$ and estimating the functions 
$\partial_3S$ and $\partial_3^2H$ uniformly in $(\theta_1,\theta_3)$. 
\qed\halflineskip

\begin{theorem}\label{201906031711}
Suppose that $[A1]$ with $(i_A,j_A,i_B,j_B,i_H,j_H)=(1,1,2,1,3,2)$, $[A2]$, 
$[A3]$ $(iii)$ and $[A4]$ are satisfied. Then 
\beas 
n^{1/2}h^{-1/2}\big(\hat{\theta}_3^0-\theta_3^*\big)
- \Gamma_{33}^{-1}M_n^{(3)}
&\to^p& 0
\eeas
as $n\to\infty$. In particular, 
\beas 
n^{1/2}h^{-1/2}\big(\hat{\theta}_3^0-\theta_3^*\big)
&\to^d&
N(0,\Gamma_{33}^{-1})
\eeas
as $n\to\infty$. 
\end{theorem}
\proof 
Use Lemmas \ref{201906031714} and \ref{201906031658}. 
\qed\halflineskip

\begin{remark}\label{201906181743}\rm 
It is possible to construct a QMLE $\hat{\vartheta}_3$ for $\theta_3$ based on 
the quasi-log likelihood function
\beas 
\calh_n^{(3)}(\theta_3)
&=&
-\half\sum_{j=1}^n\bigg\{
3V(Z_\tjm,\hat{\theta}_1^0,\theta_3)
\big[\big\{h^{-3/2}
\big(\Delta_jY-hG_n(\hat{\theta}_1^0,\hat{\theta}_2^0,\theta_3\big)\big\}^{\otimes2}
\big]
\\&&
\qquad\qquad+\log\big(3^{-1}V(Z_\tjm,\hat{\theta}_1^0,\theta_3)\big)\bigg\}.
\eeas
Then, under a certain set of conditions, we have 
\beas 
n^{1/2}h^{-1/2}\big(\hat{\vartheta}_3-\theta_3^*\big)
&\to^d&
N(0,4\Gamma_{33}^{-1}). 
\eeas
Therefore $\hat{\theta}_3^0$ is superior to $\hat{\vartheta}_3$. 
\end{remark}

\section{Adaptive one-step estimator 
{\colorr 
for $(\theta_1,\theta_2,\theta_3)$ }}\label{202001141623}
In this section, we will consider a one-step estimator for 
${\colorr\theta=(\theta_1,\theta_2,\theta_3)}$ 
given an initial estimators $(\hat{\theta}_1^0,\hat{\theta}_2^0,\hat{\theta}_3^0)$ for 
$(\theta_1,\theta_2,\theta_3)$ 
based on $(Z_\tj)_{j=0,1,...,n}$. 
{\colorr 
We will assume the following rate of convergence for each initial estimator 
\bd
\im[[A4$^\sharp$\!\!]] 
{\bf (i)} 
$\ds \hat{\theta}_1^0-\theta_1^* \yeq O_p(n^{-1/2})$ as $n\to\infty$
\bd
\im[(ii)] 
$\ds \hat{\theta}_2^0-\theta_2^*\yeq O_p(n^{-1/2}h^{-1/2})$  as $n\to\infty$
\im[(iii)] 
$\ds \hat{\theta}_3^0-\theta_3^*\yeq O_p(n^{-1/2}h^{1/2})$  as $n\to\infty$. 
\ed\ed

The initial estimator $\hat{\theta}_3^0$ is not necessarily 
the one defined in Section \ref{201906041935}, though we already know that one satisfies $[A4^\sharp]$ (iii). 
That is, the initial estimator $\hat{\theta}_3^0$ used in this section is requested to attain 
the convergence rate $n^{-1/2}h^{1/2}$ only, not to necessarily achieve 
the asymptotic variance equal to $\Gamma_{33}^{-1}$ or less. 
Thus, the estimator $\hat{\vartheta}_3$ mentioned by Remark \ref{201906181743},  
as well as $\hat{\theta}_3^0$ in Section \ref{201906041935}, 
can serve as the initial estimator of $\theta_3$.}
{\colorg 
As Section \ref{201905291607} recalls a construction of the initial estimator $\hat{\theta}_1^0$, 
in estimation of non-degenerate diffusion processes, 
there is an estimator of $\theta_1$ satisfying Condition $[A4^\sharp]$ (i) 
based on only the first equation of (\ref{201905262213}).
It is know that its information cannot be greater than the matrix 
\beas
\half\int  \rm{Tr}\big\{\big(C^{-1}(\partial_1C)C^{-1}\partial_1C\big)(z,\theta_1^*)\big\}\nu(dz).
\eeas
It will be turned out that 
the amount of information is increased  
by the one-step estimator. 
We will recall a standard construction of $\hat{\theta}_2^0$ in Section \ref{201906041938}. 
}

{\colorr 
Let 
\beas 
M_n^{(1)} &=& 
\half n^{-1/2}\sum_{j=1}^n
\big(S^{-1}(\partial_1S)S^{-1}\big)(Z_\tjm,\theta_1^*,\theta_3^*)
\big[\widetilde{\cald}_j(\theta_1^*,\theta_2^*,\theta_3^*)^{\otimes2}-S(Z_\tjm,\theta_1^*,\theta_3^*)\big]. 
\eeas
Let 
\beas 
\Gamma_{11}
&=&
\half\int  \rm{Tr}\big\{S^{-1}(\partial_1S)S^{-1}\partial_1S(z,\theta_1^*, \theta_3^*)\big\}\nu(dz)
\\&=&
\half\int \bigg[
 \rm{Tr}\big\{\big(C^{-1}(\partial_1C)C^{-1}\partial_1C\big)(z,\theta_1^*)\big\}
 \\&&\qquad\quad
 +\rm{Tr}\big\{\big(V^{-1}H_x(\partial_1C)H_x^\star V^{-1}H_x(\partial_1C)H_x^\star\big)(z,\theta_1^*,\theta_3^*)\big\}
 \bigg]\nu(dz).
\eeas
{\colorg
If $H_x$ is an invertible (square) matrix, then $\Gamma_{11}$ coincides with 
\beas
\int  \rm{Tr}\big\{\big(C^{-1}(\partial_1C)C^{-1}\partial_1C\big)(z,\theta_1^*)\big\}\nu(dz).
\eeas
Otherwise, it is not always true. 
}

Let  
\bea\label{202001281931} 
\Gamma_{22}
&=&
\int S(z,\theta_1^*,\theta_3^*)^{-1}
\left[
\left(\begin{array}{c}\partial_2A(z,\theta_2^*)\\
2^{-1}\partial_2L_H(z,\theta_1^*,\theta_2^*,\theta_3^*)
\end{array}\right)^{\otimes2}
\right]\nu(dz)
\nn\\&=&
\int \partial_2A(z,\theta_2^*)^\star C(z,\theta_1^*)^{-1}\partial_2A(z,\theta_2^*)\nu(dz).
\eea
Let 
$
\Gamma^J(\theta^*)
=
\text{diag}\big[
\Gamma_{11},\Gamma_{22},\Gamma_{33}
\big]$,  
where $\Gamma_{33}$ is defined by (\ref{201906190133}).

\begin{en-text}
\half\int \bigg[
 \rm{Tr}\big\{\big(C^{-1}(\partial_1C)C^{-1}\partial_1C\big)(z,\theta_1^*)\big\}
 \\&&\qquad\quad
 +\rm{Tr}\big\{\big(V^{-1}H_x(\partial_1C)H_x^\star V^{-1}H_x(\partial_1C)H_x^\star\big)(z,\theta_1^*,\theta_3^*)\big\}
 \bigg]\nu(dz),
 \\&&\quad\quad
 \int \partial_2A(z,\theta_2^*)^\star C(z,\theta_1^*)^{-1}\partial_2A(z,\theta_2^*)\nu(dz), 
\\&&\quad\quad
\int12\partial_3H(z,\theta_3^*)^\star V(z,\theta_1^*,\theta_3^*)^{-1}\partial_3H(z,\theta_3^*)\nu(dz)
\bigg]
\eeas
\end{en-text}
%

\begin{en-text}
\beas 
\bbH_n^{(1)}(\theta_1) 
&=& 
\bbH_n(\theta_1,\hat{\theta}_2,\hat{\theta}_3),
\eeas
where 
\end{en-text}
We will use the following random fields: 
\bea\label{201906060545} 
\bbH_n^{(1)}(\theta_1)
 &=& 
 -\half\sum_{j=1}^n \bigg\{
S(Z_\tjm,\theta_1,\hat{\theta}_3^0)^{-1}\big[\cald_j(\theta_1,\hat{\theta}_2^0,\hat{\theta}_3^0)^{\otimes2}\big]
 +\log\det S(Z_\tjm,\theta_1,\hat{\theta}_3^0)\bigg\}.
\eea
and 
\bea\label{201906030042} 
\bbH^{(2,3)}_n(\theta_2,\theta_3)
 &=& 
 -\half\sum_{j=1}^n 
 \hat{S}(Z_\tjm,\hat{\theta}_3^0)^{-1}\big[\cald_j(\hat{\theta}_1^0,\theta_2,\theta_3)^{\otimes2}\big].
\eea
Recall $\hat{S}(z,\theta_3)=S(z,\hat{\theta}_1^0,\theta_3)$. 
To construct one-step estimators, we consider the functions 
\beas 
\bbE_n(\theta_1)
&=&
\theta_1-\big[\partial_1^2\bbH^{(1)}_n(\theta_1)]^{-1}
\partial_1\bbH^{(1)}_n(\theta_1)
\eeas
and
\beas
\bbF_n(\theta_2,\theta_3) 
&=& 
\left(\begin{array}{c}\theta_2\\ \theta_3\end{array}\right)
-\big[\partial_{(\theta_2,\theta_3)}^2\bbH_n^{(2,3)}\big(\theta_2,\theta_3\big)\big]^{-1}
\partial_{(\theta_2,\theta_3)}\bbH_n^{(2,3)}\big(\theta_2,\theta_3\big)
\eeas
when both matrices $\partial_1^2\bbH^{(1)}_n(\theta_1)$ and 
$\partial_{(\theta_2,\theta_3)}^2\bbH_n^{(2,3)}\big(\theta_2,\theta_3\big)$ are invertible. 
{\colorg 
Let 
\beas 
\calx_n^{(1)} &=& \big\{\omega\in\Omega;\>\partial_1^2\bbH_n^{(1)}(\hat{\theta}_1^0)
\text{ is invertible and }\bbE_n(\hat{\theta}_1^0)\in\Theta_1\big\}
\eeas
and
\beas 
\calx_n^{(2,3)} &=& \big\{\omega\in\Omega;\>
\partial_{(\theta_2,\theta_3)}^2\bbH_n^{(2,3)}\big(\hat{\theta}_2^0,\hat{\theta}_3^0\big)
\text{ is invertible and }\bbF_n(\hat{\theta}_2^0,\hat{\theta}_3^0)\in\Theta_2\times\Theta_3\big\}
\eeas
Let $\calx_n=\calx_n^{(1)}\cap\calx_n^{(2,3)}$.}
\begin{en-text}
Let $\calx_n=\big\{\omega\in\Omega;\>$ both 
$\partial_1^2\bbH_n^{(1)}(\hat{\theta}_1^0)$ and 
$\partial_{(\theta_2,\theta_3)}^2\bbH_n^{(2,3)}\big(\hat{\theta}_2^0,\hat{\theta}_3^0\big)$ are invertible, 
and $\big(\bbE_n(\hat{\theta}_1^0), \bbF_n(\hat{\theta}_2^0,\hat{\theta}_3^0)\big)\in\Theta\big\}$. 
\end{en-text}
The event $\calx_n$ is a statistic because it is determined by the data $(Z_\tj)_{j=0,...,n}$ only. 
For $(\theta_1,\theta_2,\theta_3)$, the one-step estimator $(\hat{\theta}_1,\hat{\theta}_2,\hat{\theta}_3)$ 
with the initial estimator $(\hat{\theta}_1^0,\hat{\theta}_2^0,\hat{\theta}_3^0)$ 
is defined by
\beas 
\left(\begin{array}{c}\hat{\theta}_1\\\hat{\theta}_2\\ \hat{\theta}_3\end{array}\right)
&=&
\left\{\begin{array}{cl}
\left(\begin{array}{c} \bbE_n(\hat{\theta}_1^0)\\
\bbF_n\big(\hat{\theta}_2^0,\hat{\theta}_3^0\big)
\end{array}\right)
&\text{on }\calx_n\y
\upsilon&\text{on }\calx_n^c
\end{array}\right.
\eeas
where $\upsilon$ is an arbitrary value in $\Theta$. 
\begin{en-text}
\beas 
\left(\begin{array}{c}\hat{\theta}_2\\ \hat{\theta}_3\end{array}\right)
&=&
\left(\begin{array}{c}\hat{\theta}_2^0\\ \hat{\theta}_3^0\end{array}\right)
-\big[\partial_{(\theta_2,\theta_3)}^2\bbH_n\big(\hat{\theta}_2^0,\hat{\theta}_3^0\big)\big]^{-1}
\partial_{(\theta_2,\theta_3)}\bbH_n\big(\hat{\theta}_2^0,\hat{\theta}_3^0\big)
\eeas
on the event $\calx_n$
\end{en-text}

Let $\hat{\gamma}=\big(\hat{\theta}_2,\hat{\theta}_3\big)^\star$, 
$\hat{\gamma}^0=\big(\hat{\theta}_2^0,\hat{\theta}_3^0\big)^\star$ and 
$\gamma^*=\big(\theta_2^*,\theta_3^*\big)^\star$. 
Let $U$ be an open ball in $\bbR^{\sfp_2+\sfp_3}$ centered at $\gamma^*$ such that $U\subset\Theta_2\times\Theta_3$. 
Let $\calx_n^{*{\colorg(2,3)}}=\calx_n^{{\colorg(2,3)}}\cap\{\hat{\gamma}^0\in U\}$.

\begin{lemma}\label{201906040514}
Suppose that $[A1]$ with $(i_A,j_A,i_B,j_B,i_H,j_H)=(1,{\colorb 2},2,1,3,1)$, $[A2]$ $(i)$, 
{\colorb $(iii)$, $(iv)$} 
and $[A4^\sharp]$ are satisfied. 
Then 
\beas 
n^{-1/2}h^{-1/2}\>\partial_2\bbH_n^{(2,3)}(\hat{\gamma}^0)
&=& 
O_p(1)
\eeas
as $n\to\infty$. 
\end{lemma}
\proof 
By using Lemma \ref{201906040435} and Lemma \ref{201906031516} (b) 
together with the convergence rate of the initial estimators, we have 
\beas 
n^{-1/2}h^{-1/2}\>\partial_2\bbH_n^{(2,3)}(\hat{\gamma}^0)
&=&
n^{-1/2}\Psi_2(\hat{\theta}_1^0,\hat{\theta}_2^0,\hat{\theta}_3^0,\hat{\theta}_1^0,\hat{\theta}_2^0,\hat{\theta}_3^0)
\\&=&
n^{-1/2}\Psi_2(\hat{\theta}_1^0,\hat{\theta}_2^0,\hat{\theta}_3^0,\theta_1^*,\theta_2^*,\theta_3^*)
+O_p(1)
\\&=&
n^{-1/2}\widetilde{\Psi}_2(\hat{\theta}_1^0,\hat{\theta}_2^0,\hat{\theta}_3^0,\theta_1^*,\theta_2^*,\theta_3^*)
+O_p(1)
\eeas
by Lemma \ref{201906031518} and Lemma \ref{201906031516} (a). 

The open ball of radius $r$ centered at $\theta$ is denoted by $U(\theta,r)$. 
Define the random field 
\bea\label{201906040542} 
\Phi_n^{(\ref{201906040542})}(\theta)
&=& 
n^{-1/2}\widetilde{\Psi}_2(\theta_1,\theta_2,\theta_3,\theta_1^*,\theta_2^*,\theta_3^*)
\eea
on $\theta=(\theta_1,\theta_2,\theta_3)\in U(\theta^*,r)$ for a small number $r$ such that 
$U(\theta^*,r)\subset\Theta$. 
With the Burkholder-Davis-Gundy inequality 
{\colorb and in particular twice differentiability of $A$ in $\theta_2$}, we obtain 
\beas 
\sup_n\sum_{i=0,1}\sup_{\theta\in B(\theta^*,r)} \big\||\partial_\theta^i\Phi_n^{(\ref{201906040542})}(\theta)|\big\|_p &<& \infty
\eeas
for every $p>1$. Therefore, 
Sobolev's inequality ensures 
\beas 
\sup_n \bigg\|\sup_{\theta\in U(\theta^*,r)}|\Phi_n^{(\ref{201906040542})}(\theta)|\bigg\|_p &<& \infty
\eeas
Consequently, 
\beas 
\Phi_n^{(\ref{201906040542})}\big(\hat{\theta}_1^0,\hat{\theta}_2^0,\hat{\theta}_3^0\big)
1_{\big\{(\hat{\theta}_1^0,\hat{\theta}_2^0,\hat{\theta}_3^0)\in U(\theta^*,r)\big\}}
&=& 
O_p(1).
\eeas
This completes the proof. 
\qed\halflineskip

\begin{lemma}\label{201906040519}
Suppose that $[A1]$ with $(i_A,j_A,i_B,j_B,i_H,j_H)=(1,1,2,1,3,{\colorb 2})$, 
$[A2]$ $(i)$, {\colorb $(iii)$, $(iv)$} 
and $[A4^\sharp]$ are satisfied. Then 
\beas 
n^{-1/2}h^{1/2}\>\partial_3\bbH_n^{(2,3)}(\hat{\gamma}^0)
&=& 
O_p(1)
\eeas
as $n\to\infty$. 
\end{lemma}
\proof 
The proof is similar to that of Lemma \ref{201906040514}. 
First, 
\beas 
n^{-1/2}h^{1/2}\>\partial_3\bbH_n^{(2,3)}(\hat{\gamma}^0)
&=&
n^{-1/2}\Psi_3(\hat{\theta}_1^0,\hat{\theta}_2^0,\hat{\theta}_3^0,\hat{\theta}_1^0,\hat{\theta}_2^0,\hat{\theta}_3^0)
\\&=&
n^{-1/2}\widetilde{\Psi}_3(\hat{\theta}_1^0,\hat{\theta}_2^0,\hat{\theta}_3^0,\theta_1^*,\theta_2^*,\theta_3^*)+O_p(1).
\eeas
Then we can show the lemma in the same fashion as Lemma \ref{201906040514} 
with a random field. 
\qed\halflineskip

Let 
\beas 
B_n=U\big(\theta_1^*,n^{-1/2}\log (nh)\big)\times U\big(\theta_2^*,(nh)^{-1/2}\log (nh)\big)
\times U\big(\theta_3^*,n^{-1/2}h^{1/2}\log (nh)\big),
\eeas
\beas 
B_n'=U\big(\theta_2^*,(nh)^{-1/2}\log (nh)\big)
\times U\big(\theta_3^*,n^{-1/2}h^{1/2}\log (nh)\big)
\eeas 
and 
{\colorb
\beas 
B_n''=U\big(\theta_1^*,n^{-1/2}\log (nh)\big)
\times U\big(\theta_3^*,n^{-1/2}h^{1/2}\log (nh)\big).
\eeas 

We will use the following random fields. 
{\colorb
\beas
\Phi_{22,1}(\theta_1,\theta_3,\theta_1',\theta_2',\theta_3')
&=&
-\sum_{j=1}^nS(Z_\tjm,\theta_1,\theta_3)^{-1}
\left[
\left(\begin{array}{c}\partial_2A(Z_\tjm,\theta_2')\\
2^{-1}\partial_2L_H(Z_\tjm,\theta_1',\theta_2',\theta_3')
\end{array}\right)^{\otimes2}
\right]
\eeas
}
{\colorb 
\beas&&
\Phi_{22,2}(\theta_1,\theta_3,\theta_1',\theta_2',\theta_3',
\theta_2'',\theta_3'')
\nn\\&=&
\sum_{j=1}^nS(Z_\tjm,\theta_1,\theta_3)^{-1}
\left[\cald_j(\theta_1',\theta_2',\theta_3'),\>
\left(\begin{array}{c}\partial_2^2A(Z_\tjm,\theta_2'')\\
2^{-1}H_x(Z_\tjm,\theta_3'')\big[\partial_2^2A(Z_\tjm,\theta_2'')\big]
\end{array}\right)
\right]
\eeas
}
{\colorb 
\beas&&
\widetilde{\Phi}_{22,2}(\theta_1,\theta_3,\theta_1',\theta_2',\theta_3',
\theta_2'',\theta_3'')
\nn\\&=&
\sum_{j=1}^nS(Z_\tjm,\theta_1,\theta_3)^{-1}
\left[\widetilde{\cald}_j(\theta_1',\theta_2',\theta_3'),\>
\left(\begin{array}{c}\partial_2^2A(Z_\tjm,\theta_2'')\\
2^{-1}H_x(Z_\tjm,\theta_3'')\big[\partial_2^2A(Z_\tjm,\theta_2'')\big]
\end{array}\right)
\right]
\eeas

\beas &&
\Phi_{23,1}(\theta_1,\theta_3,\theta_1',\theta_2',\theta_3',\theta_2'',\theta_3'')
\\&=&
-\sum_{j=1}^nS(Z_\tjm,\theta_1,\theta_3)^{-1}
\left[\left(\begin{array}{c}0\\
2^{-1}\partial_3L_H(Z_\tjm,\theta_1',\theta_2',\theta_3')
\end{array}\right)\>\right.
\\&&\hspace{5cm}\left.\otimes
\left(\begin{array}{c}\partial_2A(Z_\tjm,\theta_2'')\\
2^{-1}H_x(Z_\tjm,\theta_3'')\big[\partial_2A(Z_\tjm,\theta_2'')\big]
\end{array}\right)
\right]
\eeas
\beas&&
\Phi_{23,2}(\theta_1,\theta_3,\theta_1',\theta_2',\theta_3',\theta_2'',\theta_3'')
\\&=&
\sum_{j=1}^nS(Z_\tjm,\theta_1,\theta_3)^{-1}
\left[\cald_j(\theta_1',\theta_2',\theta_3'),\>
\left(\begin{array}{c}0\\
2^{-1}\partial_3H_x(Z_\tjm,\theta_3'')\big[\partial_2A(Z_\tjm,\theta_2'')\big]
\end{array}\right)
\right]
\eeas
\beas
\Phi_{33,1}(\theta_1,\theta_3,\theta_1',\theta_2',\theta_3')
&=&
-\sum_{j=1}^nS(Z_\tjm,\theta_1,\theta_3)^{-1}
\left[
\left(\begin{array}{c}0\\
\partial_3H(Z_\tjm,\theta_3')+2^{-1}h\partial_3L_H(Z_\tjm,\theta_1',\theta_2',\theta_3')
\end{array}\right)^{\otimes2}
\right]
\eeas
\beas
\Phi_{33,2}(\theta_1,\theta_3,\theta_1',\theta_2',\theta_3',\theta_1'',\theta_2'',\theta_3'')
&=&
\sum_{j=1}^nS(Z_\tjm,\theta_1,\theta_3)^{-1}
\bigg[\cald_j(\theta_1',\theta_2',\theta_3')
\\&&\qquad\otimes
\left(\begin{array}{c}0\\
\partial_3^2 H(Z_\tjm,\theta_3'')
+2^{-1}h\partial_3^2 L_H(Z_\tjm,\theta_1'',\theta_2'',\theta_3'')
\end{array}\right)\bigg]
\eeas
}

}
\begin{lemma}\label{201906040536}
Suppose that $[A1]$ with $(i_A,j_A,i_B,j_B,i_H,j_H)=(1,{\colorr3},2,{\colorr1},3,{\colorr1})$,  $[A2]$ and $[A4^\sharp]$ are satisfied. 
Then 
\beas 
\sup_{(\theta_2,\theta_3)\in B_n'}\big|
n^{-1}h^{-1}\>\partial_2^2\bbH_n^{(2,3)}(\theta_2,\theta_3)
+
\Gamma_{22}
\big|
&\to^p&
0
\eeas
as $n\to\infty$.
\end{lemma}
\proof 
\begin{en-text}
Let 
\bea\label{201906040605}
\Phi^{(\ref{201906040605})}_n(\theta_1,\theta_2,\theta_3)
&=&
-n^{-1}\sum_{j=1}^nS(Z_\tjm,\theta_1,\theta_3)^{-1}
\left[
\left(\begin{array}{c}\partial_2A(Z_\tjm,\theta_2)\\
2^{-1}\partial_2L_H(Z_\tjm,\theta_1,\theta_2,\theta_3)
\end{array}\right)^{\otimes2}
\right]
\eea
and
\bea\label{201906040634}&&
\Phi^{(\ref{201906040634})}_n(\theta_1,\theta_2,\theta_3,\theta_1',\theta_2',\theta_3')
\nn\\&=&
n^{-1}h^{-1/2}\sum_{j=1}^nS(Z_\tjm,\theta_1',\theta_3')^{-1}
\left[\cald_j(\theta_1,\theta_2,\theta_3),\>
\left(\begin{array}{c}\partial_2^2A(Z_\tjm,\theta_2')\\
2^{-1}\partial_2^2L_H(Z_\tjm,\theta_1',\theta_2',\theta_3')
\end{array}\right)
\right].
\eea
\end{en-text}
{\colorb 
We have 
\bea\label{201906040705} 
n^{-1}h^{-1}\>\partial_2^2\bbH_n^{(2,3)}(\theta_2,\theta_3)
&=&
n^{-1}\Phi_{22,1}(\hat{\theta}_1^0,\hat{\theta}_3^0,\hat{\theta}_1^0,\theta_2,\theta_3)
\nn\\&&
+n^{-1}h^{-1/2}\Phi_{22,2}(\hat{\theta}_1^0,\hat{\theta}_3^0,\hat{\theta}_1^0,\theta_2,\theta_3,\theta_2,\theta_3)
\eea
Apply Lemma \ref{201906040435} and Lemma \ref{201906031516} (b) to obtain 
\bea\label{201906040706}&&
\sup_{(\theta_1,\theta_3)\in B_n''}
\sup_{(\theta_1',\theta_2',\theta_3')\in B_n}
\sup_{(\theta_2'',\theta_3'')\in B_n'}
\big|
n^{-1}h^{-1/2}\Phi_{22,2}(\theta_1,\theta_3,\theta_1',\theta_2',\theta_3',
\theta_2'',\theta_3'')
\nn\\&&\qquad\qquad\qquad
-
n^{-1}h^{-1/2}\Phi_{22,2}(\theta_1,\theta_3,\theta_1^*,\theta_2^*,\theta_3^*,
\theta_2'',\theta_3'')
\big|
\nn\\&=&
o_p(1).
\nn\\&&
\eea
Here we used the assumption that the functions are bound by a polynomial in $z$ uniformly in 
the parameters, and the count
\beas 
n^{-1}h^{-1/2}\times n \times h^{-1/2}\times n^{-1/2}h^{1/2}\log(nh)
&=&
    \frac{\log(nh)}{\sqrt{nh}}
\eeas
to estimate the error when replacing $\theta_3'$ by $\theta_3^*$, 
{\colorg as well a similar count when replacing $(\theta_1',\theta_2')$ by $(\theta_1^*,\theta_2^*)$.}

We apply Lemmas \ref{201906031518} and \ref{201906031516} (a) to obtain 
\bea\label{201906040707}&&
\sup_{(\theta_1,\theta_3)\in B_n''}
\sup_{(\theta_2'',\theta_3'')\in B_n'}
\big|
n^{-1}h^{-1/2}\Phi_{22,2}(\theta_1,\theta_3,\theta_1^*,\theta_2^*,\theta_3^*,
\theta_2'',\theta_3'')
\nn\\&&\qquad\qquad\qquad
-
n^{-1}h^{-1/2}\widetilde{\Phi}_{22,2}(\theta_1,\theta_3,\theta_1^*,\theta_2^*,\theta_3^*,
\theta_2'',\theta_3'')
\big|
\nn\\&=&
{\colorg O_p\big((nh)^{-1/2}\log(nh)\big)}
\yeq o_p(1).
\eea
Since $\widetilde{\cald}_j(\theta_1^*,\theta_2^*,\theta_3^*)$ 
in $\widetilde{\Phi}_{22,2}$ are martingale differences with respect to a suitable filtration, 
we can conclude by the random field argument with the Sobolev space of index $(1,p)$, 
$p>1$,  
that 
\bea\label{202001131615}
\sup_{(\theta_1,\theta_3)\in B_n''}
\sup_{(\theta_2'',\theta_3'')\in B_n'}
\big|n^{-1}h^{-1/2}\widetilde{\Phi}_{22,2}(\theta_1,\theta_3,\theta_1^*,\theta_2^*,\theta_3^*,
\theta_2'',\theta_3'')\big|
&=&
O_p((nh)^{-1/2})\yeq o_p(1)
\nn\\&&
\eea

On the other hand, 
\bea\label{201906040709}
\sup_{(\theta_1,\theta_3)\in B_n''}
\sup_{(\theta_1',\theta_2',\theta_3')\in B_n}
\big|n^{-1}\Phi_{22,1}(\theta_1,\theta_3,\theta_1',\theta_2',\theta_3')
-n^{-1}\Phi_{22,1}(\theta_1^*,\theta_3^*,\theta_1^*,\theta_2^*,\theta_3^*)
\big|
&=&
o_p(1)
\eea
From (\ref{201906040705})-(\ref{201906040709}) 
and $[A4^\sharp]$ (i), (iii), we obtain 
\bea\label{201906040711} 
\sup_{(\theta_2,\theta_3)\in B_n'}
\big|n^{-1}h^{-1}\>\partial_2^2\bbH_n^{(2,3)}(\theta_2,\theta_3)
-
n^{-1}\Phi_{22,1}(\theta_1^*,\theta_3^*,\theta_1^*,\theta_2^*,\theta_3^*)\big|
&=&
o_p(1).
\eea
Now the assertion of the lemma is easy to obtain if one uses 
$[A1]$, $[A2]$ and Lemma \ref{201906041850}. 
\qed\halflineskip

Let 
\bea\label{201909181526} 
i(z,\theta)
&=&
\left(\begin{array}{cc}
\partial_2A(z,\theta_2)^\star&2^{-1}\partial_2L_H(z,\theta_1,\theta_2,\theta_3)^\star \\ 
O& \partial_3H(z,\theta_3)^\star
\end{array}\right)
 S(z,\theta_1,\theta_3)^{-1}
\nn\\&&\times
 \left(\begin{array}{cc}
\partial_2A(z,\theta_2)&O \\ 
2^{-1}\partial_2L_H(z,\theta_1,\theta_2,\theta_3)& \partial_3H(z,\theta_3)
\end{array}\right).
\nn\\&&
\eea
Then simple calculus with (\ref{201906040958}) and 
\beas 
\partial_2L_H(z,\theta_1,\theta_2,\theta_3)
&=&
H_x(z,\theta_3)\big[\partial_2A(z,\theta_2)\big]
\eeas
yield 
\bea\label{201909181528}
i(z,\theta)
\nn&=&
\left(\begin{array}{cc}
\partial_2A(z,\theta_2)^\star C(z,\theta_1)^{-1}\partial_2A(z,\theta_2)&O \\ 
O& 12\partial_3H(z,\theta_3)^\star V(z,\theta_1,\theta_3)^{-1}\partial_3H(z,\theta_3)
\end{array}\right). 
\\&&
\eea

}

\begin{lemma}\label{201906040731}
Suppose that $[A1]$ with $(i_A,j_A,i_B,j_B,i_H,j_H)=(1,1,{\colorg2},1,3,1)
$ and $[A2]$ are satisfied. 
Then 
\beas 
\sup_{(\theta_2,\theta_3)\in B_n'}\big|
n^{-1}\>\partial_3\partial_2\bbH_n^{(2,3)}(\theta_2,\theta_3)
\big|
&\to^p&
0
\eeas
as $n\to\infty$. 
\end{lemma}
\proof 
{\colorb From (\ref{201909181526}) and (\ref{201909181528}), 
we see}
\beas 
S(z,\theta_1,\theta_3)^{-1}
\left[\left(\begin{array}{c}0\\
\partial_3H(z,\theta_3)
\end{array}\right),\>
\left(\begin{array}{c}\partial_2A(z,\theta_2)\\
2^{-1}\partial_2L_H(z,\theta_1,\theta_2,\theta_3)
\end{array}\right)
\right]
&=&0.
\eeas
Then, by definition, 
{\colorb 
\beas
n^{-1}\>\partial_3\partial_2\bbH_n^{(2,3)}(\theta_2,\theta_3)
&=&
n^{-1}h\Phi_{23,1}(\hat{\theta}_1^0,\hat{\theta}_3^0,
\hat{\theta}_1^0,\theta_2,\theta_3,\theta_2,\theta_3)
\\&&
+
n^{-1}h^{1/2}\Phi_{23,2}(\hat{\theta}_1^0,\hat{\theta}_3^0,
\hat{\theta}_1^0,\theta_2,\theta_3,\theta_2,\theta_3). 
\eeas
}\noindent
Now it is not difficult to show the desired result. 
\qed\halflineskip

\begin{lemma}\label{201906040742}
Suppose that $[A1]$ with $(i_A,j_A,i_B,j_B,i_H,j_H)={\coloro(1,1,{\colorg2},1,3,2)}
$ 
and $[A2]$ are satisfied. 
Then 
\beas 
\sup_{(\theta_2,\theta_3)\in B_n'}\big|
n^{-1}h\>\partial_3^2\bbH_n^{(2,3)}(\theta_2,\theta_3)
+
\Gamma_{33}
\big|
&\to^p&
0
\eeas
as $n\to\infty$. 
\end{lemma}
\proof 
By definition, 
{\colorb 
\beas 
n^{-1}h\>\partial_3^2\bbH_n^{(2,3)}(\theta_2,\theta_3)
&=&
n^{-1}\Phi_{33,1}(\hat{\theta}_1^0,\hat{\theta}_3^0,\hat{\theta}_1^0,
\theta_2,\theta_3)
\\&&
+
n^{-1}h^{1/2}\Phi_{33,2}(\hat{\theta}_1^0,\hat{\theta}_3^0,
\hat{\theta}_1^0,\theta_2,\theta_3,\hat{\theta}_1^0,\theta_2,\theta_3). 
\eeas
{\coloro $\Phi_{33,1}$ involves the first derivative $\partial_3$, and $\Phi_{33,2}$ does 
the second derivative $\partial_3^2$.} 
First applying Lemma \ref{201906040435} and Lemma \ref{201906031516} (b), 
and next Lemma \ref{201906021438} (b), we have 
\beas &&
\sup_{(\theta_2,\theta_3)\in B_n'}\big|
n^{-1}h^{1/2}\Phi_{33,2}(\hat{\theta}_1^0,\hat{\theta}_3^0,
\hat{\theta}_1^0,\theta_2,\theta_3,\hat{\theta}_1^0,\theta_2,\theta_3)
\big|
\\&\leq&
\sup_{(\theta_2,\theta_3)\in B_n'}\big|n^{-1}h^{1/2}\Phi_{33,2}(\hat{\theta}_1^0,\hat{\theta}_3^0,
\theta_1^*,\theta_2^*,\theta_3^*,\hat{\theta}_1^0,\theta_2,\theta_3)
\big|+O_p(n^{-1/2}h^{1/2}\log(nh))
\\&=&
O_p(h^{1/2}).
\eeas
Moreover, it is easy to show 
\beas 
\sup_{(\theta_2,\theta_3)\in B_n'}\big|
n^{-1}\Phi_{33,1}(\hat{\theta}_1^0,\hat{\theta}_3^0,\hat{\theta}_1^0,
\theta_2,\theta_3)
{\colorg+}\Gamma_{33}\big|
&\to^p&
0
\eeas
from $[A1]$, $[A2]$ with the aid of Lemma \ref{201906041850}. 
}
\qed\halflineskip

Let 
\beas 
a_n &=& \left( \begin{array}{cc}
 {\colorr n^{-1/2}}h^{{\colorr-1/2}}&0\\
 0&{\colorr n^{-1/2}}h^{{\colorr1/2}}
\end{array} \right).
\eeas

\begin{lemma}\label{201906040756}
Suppose that $[A1]$ with $(i_A,j_A,i_B,j_B,i_H,j_H)=
(1,{\colorb 3},2,{\coloro1},3,{\colorb2})
$ 
and $[A2]$ are satisfied. 
Then 
\bea\label{201906040804} 
\sup_{(\theta_2,\theta_3)\in B_n'}\big|
a_n
\partial_{(\theta_2,\theta_3)}^2\bbH^{(2,3)}_n(\theta_2,\theta_3)
a_n
+
\Gamma^{(2,3)}(\theta^*)
 \big|
\>\to^p\>0
\eea
where 
{\colorb
\beas
\Gamma^{(2,3)}(\theta^*)
&=&
\left(\begin{array}{cc}
\Gamma_{22}&O \\ 
O&\Gamma_{33}
\end{array}\right). 
\eeas
}
\begin{en-text}
\beas
\Gamma^{(2,3)}(\theta^*)
&=&
\left[\begin{array}{cc}
\int \partial_2A(z,\theta_2^*)^\star C(z,\theta_1^*)^{-1}\partial_2A(z,\theta_2^*)\nu(dz)&O \\ 
O& \int12\partial_3H(z,\theta_3^*)^\star V(z,\theta_1^*,\theta_3^*)^{-1}\partial_3H(z,\theta_3^*)\nu(dz)
\end{array}\right]. 
\eeas
\end{en-text}
\end{lemma}
\proof 
The convergence (\ref{201906040804}) follows from 
Lemmas \ref{201906040536}, \ref{201906040731} and \ref{201906040742}. 
\qed\halflineskip

\begin{lemma}\label{201906041013}
Suppose that $[A1]$ with $(i_A,j_A,i_B,j_B,i_H,j_H)=
{\coloro(1,3,2,1,3,2)}
$, $[A2]$ and $[A4^{\colorb\sharp}]$ are satisfied. 
Then 
$P[\calx_n^{*{\colorg(2,3)}}]\to1$ as $n\to\infty$. 
\end{lemma}
\proof 
By Lemmas \ref{201906040514} and \ref{201906040519}, 
\beas
a_n\partial_{(\theta_2,\theta_3)}\bbH_n^{(2,3)}(\hat{\gamma}^0)
&=& 
O_p(1)
\eeas
and by Lemma \ref{201906040756}, 
\beas
\big(a_n
\partial_{(\theta_2,\theta_3)}^2\bbH^{(2,3)}_n(\hat{\gamma}^0)
a_n\big)^{-1}
&=&
O_p(1).
\eeas
Therefore, 
\beas
\big(\partial_{(\theta_2,\theta_3)}^2\bbH^{(2,3)}_n(\hat{\gamma}^0)\big)^{-1}
\partial_{(\theta_2,\theta_3)}\bbH_n^{(2,3)}(\hat{\gamma}^0)
&=&
O_p((nh)^{-1/2})
\eeas
as $n\to\infty$. 
This means $P[\calx_n^{*{\colorg(2,3)}}]\to1$. 
\qed\halflineskip

Let 
\bea\label{202001281745}
M^{(2)}_n
&=&
n^{-1/2}\sum_{j=1}^nS(Z_\tjm,\theta_1^*,\theta_3^*)^{-1}
\left[
\left(\begin{array}{c}
h^{-1/2}B(Z_\tjm,\theta_1^*)\Delta_jw\y
h^{-3/2}\kappa(Z_\tjm,\theta_1^*,\theta_3^*)\zeta_j
\end{array}\right)
,\>
\left(\begin{array}{c}\partial_2A(Z_\tjm,\theta_2^*)\\
2^{-1}\partial_2L_H(Z_\tjm,\theta_1^*,\theta_2^*,\theta_3^*)
\end{array}\right)
\right]
\nn\\&=&
n^{-1/2}\sum_{j=1}^nC(Z_\tjm,\theta_1^*)^{-1}\big[h^{-1/2}B(Z_\tjm,\theta_1^*)\Delta_jw,\>
\partial_2A(Z_\tjm,\theta_2^*)\big].
\eea

\begin{lemma}\label{201906041017}
Suppose that $[A1]$ with $(i_A,j_A,i_B,j_B,i_H,j_H)=(1,1,2,1,3,1)$, $[A2]$ and $[A4]$ are satisfied. 
Then 
\beas 
n^{-1/2}h^{-1/2}\>\partial_2\bbH_n^{(2,3)}(\theta_2^*,\theta_3^*)
-
M^{(2)}_n
&\to^p&
0
\eeas
as $n\to\infty$. 
\end{lemma}
\proof 
By using Lemma \ref{201906031516} (b) 
together with the convergence rate of the estimators {\cred $\hat{\theta}_1^0$ and $\hat{\theta}_3^0$}, 
and next by Lemma \ref{201906031516} (a) and Lemma \ref{201906031518}, 
we have 
\bea\label{201906041140} &&
n^{-1/2}h^{-1/2}\>\partial_2\bbH_n^{(2,3)}(\theta_2^*,\theta_3^*)
\nn\\&=&
n^{-1/2}\sum_{j=1}^n\hat{S}(Z_\tjm,{\colorg\hat{\theta}_3^0})^{-1}
\left[\cald_j(\hat{\theta}_1^0,\theta_2^*,\theta_3^*),\>
\left(\begin{array}{c}\partial_2A(Z_\tjm,\theta_2^*)\\
2^{-1}\partial_2L_H(Z_\tjm,\hat{\theta}_1^0,\theta_2^*,\theta_3^*)
\end{array}\right)
\right]
\nn\\&=&
n^{-1/2}\sum_{j=1}^n\hat{S}(Z_\tjm,{\colorg\hat{\theta}_3^0})^{-1}
\left[\cald_j(\theta_1^*,\theta_2^*,\theta_3^*),\>
\left(\begin{array}{c}\partial_2A(Z_\tjm,\theta_2^*)\\
2^{-1}\partial_2L_H(Z_\tjm,\hat{\theta}_1^0,\theta_2^*,\theta_3^*)
\end{array}\right)
\right]
+O_p(h^{1/2}).
\nn\\&=&
n^{-1/2}\sum_{j=1}^n\hat{S}(Z_\tjm,\theta_3^*)^{-1}
\left[\widetilde{\cald}_j(\theta_1^*,\theta_2^*,\theta_3^*),\>
\left(\begin{array}{c}\partial_2A(Z_\tjm,\theta_2^*)\\
2^{-1}\partial_2L_H(Z_\tjm,\hat{\theta}_1^0,\theta_2^*,\theta_3^*)
\end{array}\right)
\right]
+O_p(\sqrt{n}h)
+O_p(h^{1/2}).
\nn\\&&
\eea
{\colorg Here we used the derivative $\partial_1H$.}

We consider the random field  
\bea\label{201906041116}
\Phi_n^{(\ref{201906041116})}(u_1)
&=&
n^{-1/2}
\widetilde{\Psi}_2(\theta_1(u_1),\theta_2^*,\theta_3^*,\theta_1^*,\theta_2^*,\theta_3^*)
\nn\\&&
\eea
on $\{u_1\in\bbR^{\sfp_1};\>|u_1|<1\}$, 
where 
$\theta_1(u_1)=\theta_1^*+n^{-1/2}(\log n) u_1$. 
Then $L^p$-estimate of 
\beas 
\partial_1^i\{\Phi_n^{(\ref{201906041116})}(u_1)
-\Phi_n^{(\ref{201906041116})}(0)\}\quad (i=0,1) 
\eeas
yields 
\beas 
\sup_{
{\colorb u_1\in U(0,1)}
}
\big|\Phi_n^{(\ref{201906041116})}(u_1)
-\Phi_n^{(\ref{201906041116})}(0)\big|
&\to^p&
0,
\eeas
in particular, 
\beas 
\Phi_n^{(\ref{201906041116})}(u_1^\dagger)
-\Phi_n^{(\ref{201906041116})}(0)
&\to^p&
0
\eeas
where 
$u_1^\dagger=n^{1/2}(\log n)^{-1}(\hat{\theta}_1-\theta_1^*)$. 
Obviously, $M^{(2)}_n-\Phi_n^{(\ref{201906041116})}(0)\to^p0$. 
Since the first term on the right-hand side of (\ref{201906041140}) is nothing but 
$\Phi_n^{(\ref{201906041116})}(u_1^\dagger)$ 
on an event the probability of which goes to $1$, we have already obtained the result. 
\qed\halflineskip

\begin{lemma}\label{201906041155}
Suppose that $[A1]$ with $(i_A,j_A,i_B,j_B,i_H,j_H)=(1,1,2,1,3,1)$, $[A2]$ and $[A4]$ are satisfied. 
Then 
\beas 
n^{-1/2}h^{1/2}\>\partial_3\bbH_n^{(2,3)}(\theta_2^*,\theta_3^*)
-
M^{(3)}_n
&\to^p&
0
\eeas
as $n\to\infty$. 
\end{lemma}
\proof 
\begin{en-text}
Recalling the expression of $n^{-1/2}h^{1/2}\>\partial_3\bbH_n^{(2,3)}(\theta_2^*,\theta_3^*)$ 
that is similar to the one in the proof of Lemma \ref{201906040519}. 
\end{en-text}
{\colorg We have
\beas 
n^{-1/2}h^{1/2}\>\partial_3\bbH_n^{(2,3)}(\theta_2^*,\theta_3^*)
&=&
n^{-1/2}\sum_{j=1}^nS(Z_\tjm,\hat{\theta}_1^0,\hat{\theta}_3^0)^{-1}
\bigg[\cald_j(\hat{\theta}_1^0,\theta_2^*,\theta_3^*)
\\&&\qquad\otimes
\left(\begin{array}{c}0\\
\partial_3 H(Z_\tjm,\theta_3^*)
+2^{-1}h\partial_3 L_H(Z_\tjm,\hat{\theta}_1^0,\theta_2^*,\theta_3^*)
\end{array}\right)\bigg].
\eeas
}
Then this lemma can be proved in the same way as Lemma \ref{201906041017}. 
\qed\halflineskip

Let 
\beas 
M_n^{(2,3)} &=& 
\left(\begin{array}{c}M^{(2)}_n\\ M^{(3)}_n\end{array}\right).
\eeas
Combining Lemmas \ref{201906041017} and \ref{201906041155}, we obtain the following lemma. 
\begin{lemma}\label{201906041207}
Suppose that $[A1]$ with $(i_A,j_A,i_B,j_B,i_H,j_H)=(1,1,2,1,3,1)$, $[A2]$ and $[A4^{\coloro\sharp}]$ are satisfied. 
Then 
\beas 
a_n\partial_{(\theta_2,\theta_3)}\bbH_n^{(2,3)}(\theta_2^*,\theta_3^*)
-M_n^{(2,3)}
&\to^p&
0
\eeas
and $M_n^{(2,3)}\to^dN(0,\Gamma^{(2,3)}(\theta^*))$ 
as $n\to\infty$. 
In particular, 
\beas
a_n\partial_{(\theta_2,\theta_3)}\bbH_n^{(2,3)}(\theta_2^*,\theta_3^*)
&\to^d& 
N(0,\Gamma^{(2,3)}(\theta^*))
\eeas
as $n\to\infty$. 
\end{lemma}

\begin{theorem}\label{201906091921}
Suppose that $[A1]$ with $(i_A,j_A,i_B,j_B,i_H,j_H)={\coloro(1,3,2,1,3,2)}
$, $[A2]$ and $[A4^{\coloro\sharp}]$ are satisfied. 
Then 
\bea\label{201911181122} 
a_n^{-1}(\hat{\gamma}-\gamma^*) - (\Gamma^{(2,3)}(\theta^*))^{-1}M_n^{(2,3)}
&\to^p&
0
\eea
as $n\to\infty$. In particular, 
\bea\label{201911181123}
a_n^{-1}(\hat{\gamma}-\gamma^*) &\to^d& 
N(0,(\Gamma^{(2,3)}(\theta^*))^{-1})
\eea
as $n\to\infty$.
\end{theorem}
\proof
Let 
\beas 
\calx_n^{**{\colorg{(2,3)}}}
&=&
\calx_n^{{\colorg{*(2,3)}}}\cap\big\{{\coloro(\hat{\theta}_1^0,\hat{\gamma}^0}
)\in B_n\big\}\cap 
\bigg\{\sup_{\gamma\in B_n'}
\big|a_n\partial_{(\theta_2,\theta_3)}^2\bbH_n^{(2,3)}(\gamma)a_n
{\cred+} \Gamma^{(2,3)}(\theta^*)\big|
<c
\bigg\}.
\eeas
{\colorb Here $c$ is a postive constant and we will make it sufficiently small.} 
Then $P[\calx_n^{**}]\to1$ thanks to Lemmas \ref{201906041013} and \ref{201906040756}. 
On the event $\calx_n^{**{\colorg(2,3)}}$, 
{\coloro we apply Taylor's formula to obtain 
\beas&&
a_n^{-1}(\hat{\gamma}-\gamma^*)
\\&=&
\big[a_n\partial_{(\theta_2,\theta_3)}^2\bbH^{(2,3)}_n(\hat{\gamma}^0)a_n\big]^{-1}
\bigg\{
-a_n\partial_{(\theta_2,\theta_3)}\bbH^{(2,3)}_n(\gamma^*)
\\&&
{\colorg +}
a_n\int_0^1\big[
\partial_{(\theta_2,\theta_3)}^2\bbH^{(2,3)}_n(\hat{\gamma}^0)
-\partial_{(\theta_2,\theta_3)}^2\bbH^{(2,3)}_n(\hat{\gamma}(u))\big]dua_n\>
a_n^{-1}\big(\hat{\gamma}^0-\gamma^*\big)\bigg\}
\eeas
where $\hat{\gamma}(u)=\gamma^*+u(\hat{\gamma}^0-\gamma^*)$.
}
\begin{en-text}
we have $\partial_{(\theta_2,\theta_3)}\bbH^{(2,3)}_n(\hat{\gamma})=0$. 
Applying Taylor's formula, we obtain 
\beas
a_n^{-1}(\hat{\gamma}-\gamma^*)
&=&
\bigg\{-\int_0^1a_n\partial_{(\theta_2,\theta_3)}^2\bbH^{(2,3)}_n(\gamma(u))a_ndu\bigg\}^{-1}
a_n\partial_{(\theta_2,\theta_3)}\bbH^{(2,3)}_n(\gamma^*)
\eeas
{\colorb for large $n$,} 
where $\gamma(u)=\gamma^*+u(\hat{\gamma}-\gamma^*)$. 
\beas 
a_n^{-1}(\hat{\gamma}-\gamma^*) - (\Gamma^{(2,3)}(\theta^*))^{-1}M_n^{(2,3)}
&\to^p&
0.
\eeas
\end{en-text}
Then Lemmas \ref{201906040756} and \ref{201906041207} give 
{\coloro (\ref{201911181122}). Then the martingale central limit theorem gives (\ref{201911181123}).} 
\qed\halflineskip

Let 
\beas 
b_n &=& \left( \begin{array}{ccc}
 n^{-1/2}&0&0\\
 0&{\colorr n^{-1/2}}h^{{\colorr-1/2}}&0\\
 0&0&{\colorr n^{-1/2}}h^{{\colorr1/2}}
\end{array} \right).
\eeas

\begin{en-text}
Let 
\beas 
\Gamma(\theta^*) &=& \text{diag}[\Gamma^{(1)},\Gamma^{(2,3)}(\theta^*)],
\eeas
that is, 
\beas 
\Gamma(\theta^*) &=& \text{diag}\bigg[\half\int\text{Tr}\big\{C^{-1}(\partial_1C)C^{-1}(\partial_1C)(z,\theta_1^*)\big\}\nu(dz),
\\&&
\int \partial_2A(z,\theta_2^*)^\star C(z,\theta_1^*)^{-1}\partial_2A(z,\theta_2^*)\nu(dz), 
\\&&
\int12\partial_3H(z,\theta_3^*)^\star V(z,\theta_1^*,\theta_3^*)^{-1}\partial_3H(z,\theta_3^*)\nu(dz)
\bigg]
\eeas
\begin{theorem}
Let $\hat{\theta}_1$ be the estimator for $\theta_1$ in Theorem \ref{201905291601}. 
Suppose that $[A1]$ with $(i_A,j_A,i_B,j_B,i_H,j_H)=(1,1,2,3,3,3)$, $[A2]$, $[A3]$ and $[A4]$ are satisfied. 
Then 
\beas 
b_n^{-1}(\hat{\theta}-\theta^*) &\to^d&
N(0,(\Gamma(\theta^*))^{-1})
\eeas
as $n\to\infty$.
\end{theorem}
\end{en-text}

{\coloro The following notation for random fields will be used. 
\beas 
\Psi_{1,1}(\theta_1,\theta_2,\theta_3,\theta_1',\theta_2',\theta_3')
&=&
\Psi_1(\theta_1,\theta_2,\theta_3,\theta_1',\theta_2',\theta_3')
\\&=&
\sum_{j=1}^nS(Z_\tjm,\theta_1,\theta_3)^{-1}
\left[\cald_j(\theta_1',\theta_2',\theta_3'),\>
\left(\begin{array}{c}0\\
2^{-1}\partial_1L_H(Z_\tjm,\theta_1,\theta_2,\theta_3)
\end{array}\right)
\right]
\eeas
\beas 
\Psi_{1,2}(\theta_1,\theta_3,\theta_1',\theta_2',\theta_3')
&=&
\half\sum_{j=1}^n
\big(S^{-1}(\partial_1S))S^{-1}\big)(Z_\tjm,\theta_1,\theta_3)
\big[\cald_j(\theta_1',\theta_2',\theta_3')^{\otimes2}-S(Z_\tjm,\theta_1,\theta_3)\big]
\eeas
\beas 
\Psi_{11,1}(\theta_1,\theta_3,\theta_1',\theta_2',\theta_3')
&=&
\sum_{j=1}^nS(Z_\tjm,\theta_1,\theta_3)^{-1}
\left[
\left(\begin{array}{c}0\\
2^{-1}\partial_1L_H(Z_\tjm,\theta_1',\theta_2',\theta_3')
\end{array}\right)^{\otimes2}
\right]
\eeas
\beas 
\Psi_{11,2}(\theta_1,\theta_3,\theta_1',\theta_2',\theta_3',\theta_1'',\theta_2'',\theta_3'')
&=&
\sum_{j=1}^nS(Z_\tjm,\theta_1,\theta_3)^{-1}
\bigg[\cald_j(\theta_1',\theta_2',\theta_3'),\>
\\&&\otimes
\left(\begin{array}{c}0\\
2^{-1}\partial_1^2L_H(Z_\tjm,\theta_1'',\theta_2'',\theta_3'')
\end{array}\right)
\bigg]
\eeas
\beas 
\Psi_{11,3}(\theta_1,\theta_3,\theta_1',\theta_2',\theta_3')
&=&
\sum_{j=1}^n
\partial_1\big\{S^{-1}(\partial_1S)S^{-1}(Z_\tjm,\theta_1,\theta_3)\big\}
\big[\cald_j(\theta_1',\theta_2',\theta_3')^{\otimes2}-S(Z_\tjm,\theta_1',\theta_3')\big]
\eeas
\beas
\Psi_{11,4}(\theta_1,\theta_3)
&=&
\sum_{j=1}^n(S^{-1}(\partial_1S)S^{-1})(Z_\tjm,\theta_1,\theta_3)
\big[\partial_1S(Z_\tjm,\theta_1, \theta_3)\big]
\eeas
\beas
\Psi_{11,5}(\theta_1,\theta_3,\theta_1',\theta_2',\theta_3',\theta_1'',\theta_2'',\theta_3'')
&=&
\sum_{j=1}^n
(S^{-1}(\partial_1S)S^{-1})(Z_\tjm,\theta_1,\theta_3)
\bigg[\cald_j(\theta_1',\theta_2',\theta_3'),\>
\\&&
\hspace{3cm}\otimes
\left(\begin{array}{c}0\\
2^{-1}\partial_1L_H(Z_\tjm,\theta_1'',\theta_2'',\theta_3'')
\end{array}\right)
\bigg].
\eeas
}

{\coloro
\begin{lemma}\label{201911250513}
Suppose that $[A1]$ with $(i_A,j_A,i_B,j_B,i_H,j_H)=\colorg{(1,1,2,3,3,1)}
$, $[A2]$ and $[A4^\sharp]$
are satisfied. 
Then, for any sequence of positive numbers $r_n$ tending to $0$,  
\bea
\sup_{\theta_1\in U(\theta_1^*,r_n)}\big|
n^{-1}\>\partial_1^2\bbH_n^{(1)}(\theta_1)
+
\Gamma_{11}
\big|
&\to^p&
0
\eea
as $n\to\infty$.  
\end{lemma}
\proof 
By definition, 
\beas 
n^{-1}\>\partial_1^2\bbH_n^{(1)}(\theta_1)
&=&
-n^{-1}h{\coloro \Psi_{11,1}(\theta_1,\hat{\theta}_3^0,\theta_1,\hat{\theta}_2^0,\hat{\theta}_3^0)}
\\&&
+n^{-1}h^{1/2}
{\coloro \Psi_{11,2}(\theta_1,\hat{\theta}_3^0,\theta_1,\hat{\theta}_2^0,\hat{\theta}_3^0,\theta_1,\hat{\theta}_2^0,\hat{\theta}_3^0)}
\\&&
-\half n^{-1}{\coloro \Psi_{11,3}(\theta_1,\hat{\theta}_3^0,\theta_1,\hat{\theta}_2^0,\hat{\theta}_3^0)}
\\&&
-\half n^{-1}{\coloro \Psi_{11,4}(\theta_1,\hat{\theta}_3^0)}
\quad(\text{this term will remain})
\\&&
- n^{-1}h^{1/2}
{\coloro \Psi_{11,5}(\theta_1,\hat{\theta}_3^0,\theta_1,\hat{\theta}_2^0,\hat{\theta}_3^0,\theta_1,\hat{\theta}_2^0,\hat{\theta}_3^0)}
\eeas
\begin{en-text}
We may assume that $n$ is sufficiently large and $r_n\geq(nh)^{-1/2}\geq h^{1/2}$. 
\end{en-text}

We will use Condition $[A4]$ {\colorg for $\hat{\theta}_2^0$ and $\hat{\theta}_3^0$, 
and the estimate $|\theta_1-\theta_1^*|<r_n$ for $\theta_1\in U(\theta_1^*,r_n)$.}
Then 
\beas &&
\sup_{\theta_1\in U(\theta_1^*,r_n)}\big|n^{-1}\>\partial_1^2\bbH_n^{(1)}(\theta_1)+\Gamma_{11}\big|
\\&\leq&
O_p(h)
\\&&
+n^{-1}h^{1/2}
\sup_{\theta_1\in U(\theta_1^*,r_n)}\big|\Psi_{11,2}(\theta_1,\hat{\theta}_3^0,\theta_1^*,\theta_2^*,\theta_3^*,\theta_1,\hat{\theta}_2^0,\hat{\theta}_3^0)\big|
+h^{1/2}O_p(n^{-1/2}+h^{1/2})
\\&&
\quad(\text{Lemmas }\ref{201906040435}\text{ and } \ref{201906031516}(b)
)
\\&&
+n^{-1}\sup_{\theta_1\in U(\theta_1^*,r_n)}\big|\Psi_{11,3}(\theta_1,\hat{\theta}_3^0,\theta_1^*,\theta_2^*,\theta_3^*)\big|
+O_p(h^{1/2}+n^{-1/2}h^{1/2})
\\&&
\quad(\text{Lemmas }\ref{201906040435},\ref{201906031516}(b)\text{ and } \ref{201906021438}(b))
\\&&
+\bigg(-\half n^{-1}{\coloro \Psi_{11,4}(\theta_1^*,\theta_3^*)}+\Gamma_{11}\bigg)+O_p(r_n)
\\&&
+n^{-1}h^{1/2}
\sup_{\theta_1\in U(\theta_1^*,r_n)}\big|\Psi_{11,5}(\theta_1,\hat{\theta}_3^0,\theta_1,\theta_2^*,\theta_3^*,\theta_1,\hat{\theta}_2^0,\hat{\theta}_3^0)\big|
+O_p(h^{1/2}+n^{-1/2})
\\&&
\quad(\text{Lemmas }\ref{201906040435}\text{ and } \ref{201906031516}(b))
\\&=&
O_p(h)
\\&&
+O_p(h^{1/2})\quad(\text{Lemma }\ref{201906021438}(b))
\\&&
+O_p(n^{-1/2})+O_p(r_n)\quad(\text{random field argument with orthogonality})
\\&&
+o_p(1)\quad(\text{Lemma }\ref{201906041850}(a))
\\&&
+O_p(h^{1/2})\quad(\text{Lemma }\ref{201906031516}(b))
\\&=&
o_p(1)
\eeas
We remark that the used lemmas and appearing functions here 
require the regularity indices $(i_A,j_A,i_B,j_B,i_H,j_H)$ 
for $[A1]$ as follows: 
$(1,0,1,0,3,0)$ for Lemma \ref{201906021438}(b); 
$({\colorg1,1,2},1,2,{\colorg0})$ for Lemma \ref{201906031516}(b); 
$(0,0,0,0,2,1)$ for Lemma \ref{201906040435}; 
$j_B=3$, $j_H=1$ for random field argument for $\Psi_{11,3}$. 
\begin{en-text}
\koko{\coloro If we apply the same machinery as in the proof of Lemma \ref{201906081400}, 
it is easy to obtain the result. 
$[$It is remarked that $\partial_1^2$ appears in $\Psi_{11,2}$ and $\Psi_{11,3}$. 
Uniform-in-$\theta_1$ estimate for $\Psi_{11,2}$ is simple since it has the factor $h^{1/2}$ in front of it.  
On the other hand, we use random field argument for $\Psi_{11,3}$ after making the martingale differences. 
We need $\partial_1^3$ at this stage. $]$
For the second assertion, the argument becomes local by Lemma \ref{201906070327}, then 
Lemma \ref{201906081400} and the convergence (\ref{201911122311}) gives it 
by Taylor's formula. 
}
\end{en-text}
\qed\halflineskip
}

{\colorg 

\begin{lemma}\label{202001181833}
Suppose that $[A1]$ with $(i_A,j_A,i_B,j_B,i_H,j_H)=\colorg{(1,1,2,1,2,1)}$, $[A2]$ and $[A4^\sharp]$
are satisfied. 
Then 
\bea\label{202001181432} 
n^{-1/2}\partial_1\bbH^{(1)}_n(\hat{\theta}_1^0)&=& O_p(1)
\eea
as $n\to\infty$.
\end{lemma}
\proof 
We have the expression 
\beas 
n^{-1/2}\partial_1\bbH^{(1)}_n(\hat{\theta}_1^0)
&=&
n^{-1/2}h^{1/2}\Psi_{1,1}(\hat{\theta}_1^0,\hat{\theta}_2^0,\hat{\theta}_3^0,\hat{\theta}_1^0,\hat{\theta}_2^0,\hat{\theta}_3^0)
+n^{-1/2}\Psi_{1,2}(\hat{\theta}_1^0,\hat{\theta}_3^0,\hat{\theta}_1^0,\hat{\theta}_2^0,\hat{\theta}_3^0). 
\eeas
We use $[A4^\sharp]$ together with Lemmas \ref{201906040435} and \ref{201906031516} (b) 
to show 
\beas 
n^{-1/2}h^{1/2}\Psi_{1,1}(\hat{\theta}_1^0,\hat{\theta}_2^0,\hat{\theta}_3^0,\hat{\theta}_1^0,\hat{\theta}_2^0,\hat{\theta}_3^0)
&=& 
{\cred 
n^{-1/2}h^{1/2}\Psi_{1,1}(\hat{\theta}_1^0,\hat{\theta}_2^0,\hat{\theta}_3^0,\theta_1^*,\theta_2^*,\theta_3^*)+o_p(1)
}
\\&{\cred=}&
{\cred o_p(1)\yeq O_p(1)
}
\eeas
and 
\beas
n^{-1/2}\Psi_{1,2}(\hat{\theta}_1^0,\hat{\theta}_3^0,\hat{\theta}_1^0,\hat{\theta}_2^0,\hat{\theta}_3^0)
&=&
n^{-1/2}\Psi_{1,2}(\hat{\theta}_1^0,\hat{\theta}_3^0,\theta_1^*,\theta_2^*,\theta_3^*)
+O_p(1)
\\&=&
O_p(1)
\eeas
as $n\to\infty$. 
{\cred Here random field argument was used. }
\qed\halflineskip

\begin{lemma}\label{202001181156}
Suppose that $[A1]$ with $(i_A,j_A,i_B,j_B,i_H,j_H)=\colorg{(1,1,2,1,3,1)}$, $[A2]$ and $[A4^\sharp]$
are satisfied. 
Then 
\bea\label{202001181417} 
n^{-1/2}\partial_1\bbH^{(1)}_n(\theta_1^*)-M_n^{(1)} &\to^p& 0
\eea
as $n\to\infty$. 
In particular, 
\bea\label{20200118141}
n^{-1/2}\partial_1\bbH^{(1)}_n(\theta_1^*) &\to^d& N\big(0,\Gamma_{11}\big)
\eea
as $n\to\infty$. 
\end{lemma}
\proof 
We have 
\beas 
{\sf E}_j(\theta_2,\theta_3) 
&:=&
\cald_j(\theta_1^*,\theta_2,\theta_3) -\cald_j(\theta_1^*,\theta_2^*,\theta_3^*) 
\\&=&
\left(\begin{array}{c}
h^{1/2}\big(A(Z_\tjm,\theta_2^*)-A(Z_\tjm,\theta_2)\big)
\\
\left\{\begin{array}{c}
h^{-1/2}\big(H(Z_\tjm,\theta_3^*)-H(Z_\tjm,\theta_3)\big)
\\
+2^{-1}h^{1/2}\big(L_H(Z_\tjm,\theta_1^*,\theta_2^*,\theta_3^*)
-L_H(Z_\tjm,\theta_1^*,\theta_2,\theta_3)\big)
\end{array}
\right\}
\end{array}\right)
\eeas
Define the random field $\Xi_n(u_2,u_3)$ on $(u_2,u_3)\in U(0,1)^2$ by 
\beas &&
\Xi_n(u_2,u_3)
\\&=&
n^{-1/2}\sum_{j=1}^n\big(S^{-1}(\partial_1S)S^{-1}\big)(Z_\tjm,\theta_1^*,\theta_3^*+r_n^{(3)}u_3)
\bigg[\cald_j(\theta_1^*,\theta_2^*,\theta_3^*)
\otimes{\sf E}_j(\theta_2^*+r_n^{(2)}u_2,\theta_3^*+r_n^{(3)}u_3)\bigg]
\eeas
whre $r_n^{(2)}=(nh)^{-1/2}\log(nh)$ and $r_n^{(3)}=n^{-1/2}h^{1/2}\log(nh)$. 
Then the Burkholder-Davis-Gundy inequality gives 
\beas 
\lim_{n\to\infty}\sup_{(u_2,u_3)\in U(0,1)^2}\sum_{i=0,1}
\big\|\partial_{(u_2,u_3)}^i\Xi_n(u_2,u_3)\big\|_p
&=&
0,
\eeas
which implies 
\beas 
\sup_{(u_2,u_3)\in U(0,1)^2}\big|\Xi_n(u_2,u_3)\big|  &\to^p& 0
\eeas
under $[A4^\sharp]$, 
and hence 
\bea\label{202001181318}
n^{-1/2}\sum_{j=1}^n\big(S^{-1}(\partial_1S)S^{-1}\big)(Z_\tjm,\theta_1^*,\hat{\theta}_3^0)
\bigg[\cald_j(\theta_1^*,\theta_2^*,\theta_3^*)
\otimes{\sf E}_j(\hat{\theta}_2^0,\hat{\theta}_3^0)\bigg]
&\to^p&
0
\eea
as $n\to\infty$. 
It is easier to see 
\bea\label{202001181322}
n^{-1/2}\sum_{j=1}^n\big(S^{-1}(\partial_1S)S^{-1}\big)(Z_\tjm,\theta_1^*,\hat{\theta}_3^0)
\big[{\sf E}_j(\hat{\theta}_2^0,\hat{\theta}_3^0)^{\otimes2}\big]
&\to^p&
0
\eea
as $n\to\infty$. 
From (\ref{202001181318}) and (\ref{202001181322}), 
\bea\label{202001181326}
n^{-1/2}\Psi_{1,2}(\theta_1^*,\hat{\theta}_3^0,\theta_1^*,\hat{\theta}_2^0,\hat{\theta}_3^0)
&=&
n^{-1/2}\Psi_{1,2}(\theta_1^*,\hat{\theta}_3^0,\theta_1^*,\theta_2^*,\theta_3^*)
+o_p(1)
\nn\\&=&
n^{-1/2}\Psi_{1,2}(\theta_1^*,\theta_3^*,\theta_1^*,\theta_2^*,\theta_3^*)
+o_p(1)
\eea
as $n\to\infty$, where the last equality is by $[A4^\sharp]$. 

On the other hand, by $[A4^\sharp]$ and Lemmas \ref{201906040435} and \ref{201906031516} (b), 
we obtain 
\bea\label{202001181400}
n^{-1/2}h^{1/2}\Psi_{1,1}(\theta_1^*,\hat{\theta}_2^0,\hat{\theta}_3^0,\theta_1^*,\hat{\theta}_2^0,\hat{\theta}_3^0)
&=&
n^{-1/2}h^{1/2}\Psi_{1,1}(\theta_1^*,\hat{\theta}_2^0,\hat{\theta}_3^0,\theta_1^*,\theta_2^*,\theta_3^*)
+o_p(1)
\eea
By random field argument applied to the first term on the right-hand side of (\ref{202001181400}),  
\bea\label{202001181405}
n^{-1/2}h^{1/2}\Psi_{1,1}(\theta_1^*,\hat{\theta}_2^0,\hat{\theta}_3^0,\theta_1^*,\hat{\theta}_2^0,\hat{\theta}_3^0)
&=&
o_p(1). 
\eea

Consequently, from (\ref{202001181326}) and (\ref{202001181405}), 
we obtain the convergence (\ref{202001181417}) since 
\beas 
n^{-1/2}\partial_1\bbH^{(1)}_n(\theta_1^*)
&=&
n^{-1/2}h^{1/2}\Psi_{1,1}(\theta_1^*,\hat{\theta}_2^0,\hat{\theta}_3^0,\theta_1^*,\hat{\theta}_2^0,\hat{\theta}_3^0)
+n^{-1/2}\Psi_{1,2}(\theta_1^*,\hat{\theta}_3^0,\theta_1^*,\hat{\theta}_2^0,\hat{\theta}_3^0)
\\&=&
n^{-1/2}\Psi_{1,2}(\theta_1^*,\theta_3^*,\theta_1^*,\theta_2^*,\theta_3^*)
+o_p(1)
\\&=&
M_n^{(1)}+o_p(1)
\eeas
by using Lemmas \ref{201906030207} and \ref{201906031516} (a). 
Convergence (\ref{20200118141}) follows from this fact and 
Lemma \ref{201906041850} with $[A2]$, 
\qed\halflineskip
}

{\cred 
Finally, we obtain a limit theorem for the joint adaptive one-step estimator.  
}
\begin{theorem}
Suppose that $[A1]$ with $(i_A,j_A,i_B,j_B,i_H,j_H)=(1,{\colorg3},2,3,3,{\colorg2})
$, $[A2]$, $[A3]$ and $[A4]$ are satisfied. 
Then 
\beas 
b_n^{-1}(\hat{\theta}-\theta^*) &\to^d&
N(0,(\Gamma^J(\theta^*))^{-1})
\eeas
as $n\to\infty$.
\end{theorem}
\proof
{\colorg 
Let 
\beas 
\calx_n^{***}
&=&
\calx_n^{(1)}\cap 
\calx_n^{**(2,3)}\cap\big\{(\hat{\theta}_1^0,\hat{\gamma}^0)\in B_n\big\}
\cap
\bigg\{\sup_{\theta_1\in B_n'''}
\big|n^{-1}\partial_1^2\bbH^{(1)}_n(\theta_1){\cred+}\Gamma_{11}\big|<c_1\bigg\}
\eeas
where 
$B_n'''={\cred U}(\theta_1^*,n^{-1/2}\log n)$, and 
$c_1$ is a sufficiently small number such that 
$|{\sf A}{\cred+}\Gamma_{11}|<c_1$ implies $\det{\sf A}\not=0$ any ${\sf p}_1\times{\sf p}_1$ matrix ${\sf A}$. 
We obtain $P[\calx_n^{***}]\to1$ from Lemmas \ref{202001181833} and \ref{201911250513}.

On the event $\calx_n^{***}$, 
we apply Taylor's formula to obtain 
\beas&&
n^{1/2}(\hat{\theta}_1-\theta_1^*)
\\&=&
\big[n^{-1}\partial_{\theta_1}^2\bbH^{(1)}_n(\hat{\theta}_1^0)\big]^{-1}
\bigg\{
-n^{-1/2}\partial_{\theta_1}\bbH^{(1)}_n(\theta_1^*)
\\&&
+
n^{-1}\int_0^1\big[
\partial_1^2\bbH^{(1)}_n(\hat{\theta}_1^0)
-\partial_1^2\bbH^{(1)}_n(\hat{\theta}_1(u))\big]du\>
n^{1/2}\big(\hat{\theta}_1^0-\theta_1^*\big)\bigg\}
\eeas
where $\hat{\theta}_1(u)=\theta_1^*+u(\hat{\theta}_1^0-\theta^*)$.
Then we obtain 
\bea\label{201906091920} 
n^{1/2}\big(\hat{\theta}_1-\theta_1^*\big)
-\Gamma_{11}
^{-1}M_n^{(1)}
&\to^p& 
0
\eea
as $n\to\infty$ 
{\cred from Lemmas \ref{201911250513} and \ref{202001181156}.}
Therefore the convergence of $b_n^{-1}(\hat{\theta}-\theta^*)$ follows from 
the martingale central limit theorem and the relations 
(\ref{201911181122}) and 
(\ref{201906091920}) 
\qed\halflineskip
}

{\coloro
\section{Non-adaptive estimator}
\label{202001141632}
In this section, we consider a non-adaptive joint quasi-maximum likelihood estimator. 
This method does not require initial estimators. 
From computational point of view, adaptive methods {\cred often} have merits but 
the non-adaptive method is still theoretically interesting. 
We will work with the quasi-log likelihood function $\bbH_n(\theta)$ given by 
{\colorb
\bea\label{201909110433} 
\bbH_n(\theta)
 &=& 
 -\half\sum_{j=1}^n \bigg\{
S(Z_\tjm,\theta_1,\theta_3)^{-1}\big[\cald_j(\theta_1,\theta_2,\theta_3)^{\otimes2}\big]
 +\log\det S(Z_\tjm,\theta_1,\theta_3)\bigg\}
\eea
for $\theta=(\theta_1,\theta_2,\theta_3)$. 
}
{\colorg 
Suppose that a function $\hat{\theta}^J{\cred =(\hat{\theta}_1^J,\hat{\theta}_2^J,\hat{\theta}_3^J)}$ 
of the data maximizes $\bbH_n(\theta)$ in $\overline{\Theta}$.
}
\begin{en-text}
\bea\label{201906060545} 
\bbH_n(\theta)
 &=& 
 -\half\sum_{j=1}^n \bigg\{
S(Z_\tjm,\theta_1,\theta_3)^{-1}\big[\cald_j(\theta_1,\theta_2,\theta_3)^{\otimes2}\big]
 +\log\det S(Z_\tjm,\theta_1,\theta_3)\bigg\}.
\eea
\end{en-text}
\begin{en-text}
Let 
\beas 
G_n(z,\theta_1,\theta_2,\theta_3) 
&=&
H(z,\theta_3)+\frac{h}{2}L_H(z,\theta_1,\theta_2,\theta_3).
\eeas
\end{en-text}
Let 
\beas 
\bbD_n(\theta_1,\theta_2,\theta_3,\theta_1',\theta_2',\theta_3')
&=&
\bbH_n(\theta_1,\theta_2,\theta_3)-\bbH_n(\theta_1',\theta_2',\theta_3').
\eeas
Let 
\beas &&
\bbD_n^{[1]}(\theta_1,\theta_2,\theta_3,\theta_1',\theta_2',\theta_3')
\\&=&
-\half \sum_{j=1}^nS(Z_\tjm,\theta_1,\theta_3)^{-1}\big[\big(\cald_j(\theta_1,\theta_2,\theta_3)-\cald_j(\theta_1',\theta_2',\theta_3')\big)^{\otimes2}\big]
\\&=&
-\half\sum_{j=1}^nS(Z_\tjm,\theta_1,\theta_3)^{-1}
\left[
\left(\begin{array}{c}
h^{1/2}\big(A(Z_\tjm,\theta_2)-A(Z_\tjm,\theta_2')\big)
\\
\left\{\begin{array}{c}
h^{-1/2}\big(H(Z_\tjm,\theta_3)-H(Z_\tjm,\theta_3')\big)
\\
+2^{-1}h^{1/2}\big(L_H(Z_\tjm,\theta_1,\theta_2,\theta_3)-L_H(Z_\tjm,\theta_1',\theta_2',\theta_3')\big)
\end{array}\right\}
\end{array}\right)^{\otimes2}\right],
\eeas
\beas
\bbD_n^{[2]}(\theta_1,\theta_3,\theta_1',\theta_2',\theta_3')
&=&
h^{-1/2}\sum_{j=1}^nS(Z_\tjm,\theta_1,\theta_3)^{-1}
\left[
\cald_j(\theta_1',\theta_2',\theta_3'),\>
\left(\begin{array}{c}
0\\
H(Z_\tjm,\theta_3)-H(Z_\tjm,\theta_3')
\end{array}\right)\right],
\eeas
\beas&&
\bbD_n^{[3]}(\theta_1,\theta_2,\theta_3,\theta_1',\theta_2',\theta_3')
\\&=&
h^{1/2}\sum_{j=1}^nS(Z_\tjm,\theta_1,\theta_3)^{-1}
\left[
\cald_j(\theta_1',\theta_2',\theta_3'),\>
\left(\begin{array}{c}
A(Z_\tjm,\theta_2)-A(Z_\tjm,\theta_2')\\
2^{-1}\big(L_H(Z_\tjm,\theta_1,\theta_2,\theta_3)-L_H(Z_\tjm,\theta_1',\theta_2',\theta_3')\big)
\end{array}\right)\right],
\eeas
and 
\beas
\bbD_n^{[4]}(\theta_1,\theta_3,\theta_1',\theta_2',\theta_3')
&=&
-\half\sum_{j=1}^n\bigg\{\big(S(Z_\tjm,\theta_1,\theta_3)^{-1}-S(Z_\tjm,\theta_1',\theta_3')^{-1}\big)
\big[\cald_j(\theta_1',\theta_2',\theta_3')^{\otimes2}\big]
\\&&
+\log\frac{\det S(Z_\tjm,\theta_1,\theta_3)}{\det S(Z_\tjm,\theta_1',\theta_3')}\bigg\}.
\eeas
Then 
\bea\label{201906070103}
\bbD_n(\theta_1,\theta_2,\theta_3,\theta_1',\theta_2',\theta_3')
&=&
\bbD_n^{[1]}(\theta_1,\theta_2,\theta_3,\theta_1',\theta_2',\theta_3')
+\bbD_n^{[2]}(\theta_1,\theta_3,\theta_1',\theta_2',\theta_3')
\nn\\&&
+\bbD_n^{[3]}(\theta_1,\theta_2,\theta_3,\theta_1',\theta_2',\theta_3')
+\bbD_n^{[4]}(\theta_1,\theta_3,\theta_1',\theta_2',\theta_3').
\eea

\begin{en-text}
We have
\beas&&
n^{-1}\big\{\bbH_n(\theta_1,\theta_2,\theta_3)-\bbH_n(\theta_1^*,\theta_2^*,\theta_3^*)\big\}
\\&=&
-\half n^{-1}\sum_{j=1}^nS(Z_\tjm,\theta_1,\theta_3)^{-1}\big[\big(\cald_j(\theta_1,\theta_2,\theta_3)-\cald_j(\theta_1^*,\theta_2^*,\theta_3^*)\big)^{\otimes2}\big]
\\&&
+n^{-1}h^{-1/2}\sum_{j=1}^nS(Z_\tjm,\theta_1,\theta_3)^{-1}
\left[
\cald_j(\theta_1^*,\theta_2^*,\theta_3^*),\>
\left(\begin{array}{c}
0\\
H(Z_\tjm,\theta_3)-H(Z_\tjm,\theta_3^*)
\end{array}\right)\right]
\\&&
+n^{-1}h^{1/2}\sum_{j=1}^nS(Z_\tjm,\theta_1,\theta_3)^{-1}
\left[
\cald_j(\theta_1^*,\theta_2^*,\theta_3^*),\>
\left(\begin{array}{c}
0\\
2^{-1}\big(L(Z_\tjm,\theta_1,\theta_2,\theta_3)-L(Z_\tjm,\theta_1,\theta_2,\theta_3^*)\big)
\end{array}\right)\right]
\\&&
-\half n^{-1}\sum_{j=1}^n\big\{S(Z_\tjm,\theta_1,\theta_3)^{-1}-S(Z_\tjm,\theta_1^*,\theta_3^*)^{-1}\big\}
\big[\cald_j(\theta_1^*,\theta_2^*,\theta_3^*)^{\otimes2}\big]
\\&&
-\half n^{-1}\sum_{j=1}^n\log\frac{\log\det S(Z_\tjm,\theta_1,\theta_3)}{\log\det S(Z_\tjm,\theta_1^*,\theta_3^*)}.
\eeas
\end{en-text}

\begin{lemma}\label{201906070326}
Suppose that $[A1]$ with $(i_A,j_A,i_B,j_B,i_H,j_H)=(0,0,1,1,2,1)$ and $[A2]$ are satisfied. 
Then 
\bd\im[(a)] As $n\to\infty$, 
\bea\label{201906070150} 
\sup_{\theta\in \overline{\Theta}}\big|
n^{-1}h\big\{\bbH_n(\theta_1,\theta_2,\theta_3)-\bbH_n(\theta_1,\theta_2,\theta_3^*)\big\}
-\bbY^{(3)}(\theta_3)
\big|
\>\to^p\>0
\eea

\im[(b)] If $[A3]$ $(iii)$ is satisfied, then $\hat{\theta}^J_3\to^p\theta_3^*$ as $n\to\infty$. 
\ed
\end{lemma}
\proof 
We have 
\beas
n^{-1}h\big\{\bbH_n(\theta_1,\theta_2,\theta_3)-\bbH_n(\theta_1,\theta_2,\theta_3^*)\big\}
&=&
n^{-1}h\bbD_n(\theta_1,\theta_2,\theta_3,\theta_1,\theta_2,\theta_3^*)
\\&=&
n^{-1}h\bbD_n^{[1]}(\theta_1,\theta_2,\theta_3,\theta_1,\theta_2,\theta_3^*)
+n^{-1}h\bbD_n^{[2]}(\theta_1,\theta_3,\theta_1,\theta_2,\theta_3^*)
\nn\\&&
+n^{-1}h\bbD_n^{[3]}(\theta_1,\theta_2,\theta_3,\theta_1,\theta_2,\theta_3^*)
+n^{-1}h\bbD_n^{[4]}(\theta_1,\theta_3,\theta_1,\theta_2,\theta_3^*).
\eeas

By definition, 
\beas&&
n^{-1}h\bbD_n^{[1]}(\theta_1,\theta_2,\theta_3,\theta_1,\theta_2,\theta_3^*)
\\&=&
-\half n^{-1}\sum_{j=1}^nS(Z_\tjm,\theta_1,\theta_3)^{-1}
\left[
\left(\begin{array}{c}
0
\\
\left\{\begin{array}{c}
\big(H(Z_\tjm,\theta_3)-H(Z_\tjm,\theta_3^*)\big)
\\
+2^{-1}h\big(L_H(Z_\tjm,\theta_1,\theta_2,\theta_3)-L_H(Z_\tjm,\theta_1,\theta_2,\theta_3^*)\big)
\end{array}\right\}
\end{array}\right)^{\otimes2}\right].
\eeas
We apply Sobolev's inequality to uniformly estimate the ``$S$''-part and the ``$H$''-part; 
these estimates involve $\partial_1^i\partial_z^jB$ and $\partial_3^i\partial_z^jH$ for $i\in\{0,1\}$ and $j\in\{0,1\}$. 
For the ``$L_H$''-part, we use the assumption that the function is bound by a polynomial in $Z_\tjm$ uniformly in $\theta$. 
{\cred More precisely,} 
we obtain 
\beas 
\sup_{\theta\in\Theta}\big|n^{-1}h\bbD_n^{[1]}(\theta_1,\theta_2,\theta_3,\theta_1,\theta_2,\theta_3^*)
-\Phi^{(\ref{201906070229})}_n(\theta_1,\theta_3)
\big|
&=&
o_p(1)
\eeas
where 
\bea\label{201906070229}
\Phi^{(\ref{201906070229})}_n(\theta_1,\theta_3)
&=&
-\half n^{-1}\sum_{j=1}^nS(Z_\tjm,\theta_1,\theta_3)^{-1}
\left[
\left(\begin{array}{c}
0
\\
H(Z_\tjm,\theta_3)-H(Z_\tjm,\theta_3^*)
\end{array}\right)^{\otimes2}\right].
\eea
With the help of Lemma \ref{201906041850} (a), Taylor's formula and $[A2]$ give 
\beas 
\sup_{(\theta_1,\theta_3)\in\overline{\Theta}_1\times\overline{\Theta}_3}
\bigg|\Phi^{(\ref{201906070229})}_n(\theta_1,\theta_3)
+\int 6V(z,\theta_1,\theta_3)^{-1}\big[\big(H(z,\theta_3)-H(z,\theta_3^*)\big)^{\otimes2}\big]\nu(dz)\bigg|
&=&
o_p(1).
\eeas
{\cred The uniform-in-$(\theta_1,\theta_3)$ convergence follows from the point-wise convergence 
with the aid of the derivatives with respect to $(\theta_1,\theta_3)$. 
Remark that $\partial_xV$ therefore $\partial_xB$ is used, and $L_H$ has $H_{xx}$ in its expression. }
}

It is easy to see 
\beas 
\sup_{\theta\in\Theta}\big|n^{-1}h\bbD_n^{[2]}(\theta_1,\theta_3,\theta_1,\theta_2,\theta_3^*)
\big|
&=&
O_p(h^{1/2}),
\eeas
\beas 
\sup_{\theta\in\Theta}\big|n^{-1}h\bbD_n^{[3]}(\theta_1,\theta_2,\theta_3,\theta_1,\theta_2,\theta_3^*)
\big|
&=&
O_p(h^{3/2})
\eeas
and  
\beas 
\sup_{\theta\in\Theta}\big|n^{-1}h\bbD_n^{[4]}(\theta_1,\theta_3,\theta_1,\theta_2,\theta_3^*)
\big|
&=&
O_p(h). 
\eeas
This completes the proof of (a). 
The assertion (b) is a consequence of (a). 
In fact, for $\ep>0$, 
\beas 
\{|\hat{\theta}_3^J-\theta_3^{{\cred*}}|>\ep\}
&\subset& 
\{\bbY^{(3)}(\hat{\theta}_3^J)<-\chi_3\ep^2\}
\\&\subset& 
\bigg\{
\sup_{\theta\in\overline{\Theta}}\big|
n^{-1}h\big\{\bbH_n(\theta_1,\theta_2,\theta_3)-\bbH_n(\theta_1,\theta_2,\theta_3^*)\big\}
-\bbY^{(3)}(\theta_3)
\big|>\chi_3\ep^2/2\bigg\}
\eeas
since 
$\bbH_n(\hat{\theta}_1^J,\hat{\theta}_2^J,\hat{\theta}_3^J)-\bbH_n(\hat{\theta}_1^J,\hat{\theta}_2^J,\theta_3^*)
\geq0$. 
\qed\halflineskip

{\colorr We will derive a rate of convergence of $\hat{\theta}_3^J$ by 
the random field $\bbH_n(\hat{\theta}_1^J,\hat{\theta}_2^J,\theta_3)$. 

\begin{lemma}\label{201906071009}
Suppose that $[A1]$ with $(i_A,j_A,i_B,j_B,i_H,j_H)=(1,1,2,1,3,1)$ and {\coloro$[A2]$} 
are satisfied. 
Then 
\beas 
n^{-1/2}h^{1/2}\>\partial_3\bbH_n(\hat{\theta}_1^J,\hat{\theta}_2^J,\theta_3^*)
&=& 
{\colorr O_p(n^{1/2}h^{1/2})}
\eeas
as $n\to\infty$. 
\end{lemma}
\proof 
We first use Lemmas \ref{201906031516} (b) and \ref{201906021438} (b), 
next take out the principal part of $\cald_j(\theta_1^*,\theta_2^*,\theta_3^*)$, 
and apply argument with a random field and the Burkholder-Davis-Gundy inequality. 
By this procedure, 
{\colorb
\beas &&
n^{-1/2}h^{1/2}\>\partial_3\bbH_n(\hat{\theta}_1^J,\hat{\theta}_2^J,\theta_3^*)
\\&=&
n^{-1/2}\Psi_3(\hat{\theta}_1^J,\hat{\theta}_2^J,\theta_3^*,\hat{\theta}_1^J,\hat{\theta}_2^J,\theta_3^*)
+n^{-1/2}h^{1/2}\Psi_{3,3}(\hat{\theta}_1^J,\theta_3^*,\hat{\theta}_1^J,\hat{\theta}_2^J,\theta_3^*)
\\&=&
n^{-1/2}\Psi_3(\hat{\theta}_1^J,\hat{\theta}_2^J,\theta_3^*,\theta_1^*,\theta_2^*,\theta_3^*)+O_p(n^{1/2}h^{1/2})
\\&=&
n^{-1/2}\widetilde{\Psi}_3(\theta_1^*,\theta_2^*,\theta_3^*,\theta_1^*,\theta_2^*,\theta_3^*)+O_p(n^{1/2}h^{1/2})
\\&=&
O_p(n^{1/2}h^{1/2}).
\eeas
}
\qed\halflineskip


\begin{lemma}\label{201906080411}
Suppose that $[A1]$ with $(i_A,j_A,i_B,j_B,i_H,j_H)=(1,1,2,{\colorg1}
,3,{\colorr2})$ and $[A2]$ are satisfied. 
Then, for any sequence of positive numbers $r_n$ tending to $0$,  
\bea\label{202001190344} 
\sup_{\theta_3\in U(\theta_3^{{\coloro*}},r_n)}\big|
n^{-1}h\>\partial_3^2\bbH_n(\hat{\theta}_1^J,\hat{\theta}_2^J,\theta_3)
+
\Gamma_{33}(\hat{\theta}_1^J,\theta_3^*)
\big|
&\to^p&
0
\eea
as $n\to\infty$, where 
\beas 
\Gamma_{33}(\theta_1,\theta_3^*)
&=&
\int 12V(z,\theta_1,\theta_3^*)^{-1}
\left[\partial_3H(z,\theta_3^*)^{\otimes2}\right]\nu(dz).
\eeas
{\colorg 
If $[A3](iii')$ is satisfied, then $\Gamma_{33}(\theta_1,\theta_3^*)$ is non-degenerate uniformly in $\theta_1$ and 
}
$\hat{\theta}_3^J-\theta_3^*=O_p(h)$. 
\end{lemma}
\proof 
By definition, 
{\colorb
\beas
n^{-1}h\>\partial_3^2\bbH_n(\hat{\theta}_1^J,\hat{\theta}_2^J,\theta_3)
&=&
n^{-1}\Psi_{33,1}(\hat{\theta}_1^J,\hat{\theta}_2^J,\theta_3)
+n^{-1}h^{{\cred 1/2}}\Psi_{33,2}(\hat{\theta}_1^J,\hat{\theta}_2^J,\theta_3,\hat{\theta}_1^J,\hat{\theta}_2^J,\theta_3)
\\&&
+n^{-1}h\Psi_{33,3}(\hat{\theta}_1^J,\theta_3)
+n^{-1}h^{1/2}\Psi_{33,4}(\hat{\theta}_1^J,\hat{\theta}_2^J,\theta_3,\hat{\theta}_1^J,\hat{\theta}_2^J,\theta_3)
\\&&
+n^{-1}h\Psi_{33,5}(\hat{\theta}_1^J,\theta_3,\hat{\theta}_1^J,\hat{\theta}_2^J,\theta_3).
\eeas
By Lemmass \ref{201906031516} (b) and \ref{201906040435}, we have 
\beas
n^{-1}h^{1/2}\Psi_{33,4}(\hat{\theta}_1^J,\hat{\theta}_2^J,\theta_3,\hat{\theta}_1^J,\hat{\theta}_2^J,\theta_3)
&=&
n^{-1}h^{1/2}\Psi_{33,4}(\hat{\theta}_1^J,\hat{\theta}_2^J,\theta_3,\theta_1^*,\theta_2^*,\theta_3^*)+O_p(r_n)+O_p(h)
\\&=&
O(h^{1/2}+r_n)
\eeas
and this error is uniform in $\theta_3\in B(\theta_3^*,r_n)$. 
Here Lemma \ref{201906021729} was applied to estimate 
the factor 
$(S^{-1}(\partial_3S)S^{-1})(Z_\tjm,\hat{\theta}_1^J,\theta_3)$. 
{\cred Estimation of the term involving $\Psi_{33,2}$ is similar. }
Estimation of other terms is simpler. 
The term $n^{-1}\Psi_{33,1}(\hat{\theta}_1^J,\hat{\theta}_2^J,\theta_3)$ 
is approximated by $\Gamma_{33}(\hat{\theta}_1^J,\theta_3^*)$ uniformly in $B(\theta_3^*,r_n)$.
}
{\coloro
Remark that $\partial_3^2$ appears in $\Psi_{33,2}$ and $\Psi_{33,5}$. 
We do not need further differentiation with respect to $\theta_3$ to estimate them, 
because they are accompanied with the factor $h$ and the uniform-in-$\theta_3$ estimate for each term 
is carried out by simple $L^p$ estimate without random field argument. 
}

{\colorg 
Condition $[A3]$ (iii$'$) implies $[A3]$ (iii). 
We obtain the rate of convergence of $\hat{\theta}_3^J$ from 
the consistency given in Lemma \ref{201906070326} (b), 
Lemma \ref{201906071009} and (\ref{202001190344}), 
if applying the Taylor formula and $\partial_3\bbH_n(\hat{\theta}_1^J,\hat{\theta}_2^J,\hat{\theta}_3^J)=0$ 
on an event with probability tending to $1$. 
}
\qed\halflineskip

\begin{en-text}
Let 
\beas 
\bbY^{(J,1)}(\theta_1)
&=& 
-\half\int \bigg\{\text{Tr}\big(C(z,\theta_1)^{-1}C(z,\theta_1^*)+V(z,\theta_1,\theta_3^*)^{-1}V(z,\theta_1^*,\theta_3^*)
-I_{\sfd_Z}\big)
\\&&
+\log \frac{\det C(z,\theta_1)\det V(\theta_1,\theta_3^*)}{\det C(z,\theta_1^*)\det V(\theta_1^*,\theta_3^*)}
\bigg\}
\nu(dz)
\eeas
\end{en-text}

\begin{lemma}\label{201906070327}
Suppose that $[A1]$ with $(i_A,j_A,i_B,j_B,i_H,j_H)=
(1,1,2,{\colorg 1},3,2)
$, $[A2]$ {\coloro and $[A3](iii')$} are satisfied. 
Then 
\bd\im[(a)] As $n\to\infty$, 
\bea
\sup_{(\theta_1,\theta_2)\in \overline{\Theta}_1\times\overline{\Theta}_2}\big|
n^{-1}\big\{\bbH_n(\theta_1,\theta_2,\hat{\theta}_3^J)-\bbH_n(\theta_1^*,\theta_2,\hat{\theta}_3^J)\big\}
-\bbY^{(J,1)}(\theta_1)
\big|
\>\to^p\>0
\eea

\im[(b)] If $[A3]$ $(i')$ is satisfied, then $\hat{\theta}^J_1\to^p\theta_1^*$ as $n\to\infty$. 
\ed
\end{lemma}
\proof 
We have 
\beas&&
n^{-1}\big\{\bbH_n(\theta_1,\theta_2,\hat{\theta}_3^J)-\bbH_n(\theta_1^*,\theta_2,\hat{\theta}_3^J)\big\}
\\&=&
n^{-1}\bbD_n(\theta_1,\theta_2,\hat{\theta}_3^J,\theta_1^*,\theta_2,\hat{\theta}_3^J)
\\&=&
n^{-1}\bbD_n^{[1]}(\theta_1,\theta_2,\hat{\theta}_3^J,\theta_1^*,\theta_2,\hat{\theta}_3^J)
+n^{-1}\bbD_n^{[2]}(\theta_1,\hat{\theta}_3^J,\theta_1^*,\theta_2,\hat{\theta}_3^J)
\nn\\&&
+n^{-1}\bbD_n^{[3]}(\theta_1,\theta_2,\hat{\theta}_3^J,\theta_1^*,\theta_2,\hat{\theta}_3^J)
+n^{-1}\bbD_n^{[4]}(\theta_1,\hat{\theta}_3^J,\theta_1^*,\theta_2,\hat{\theta}_3^J).
\eeas

We have 
\beas &&
n^{-1}\bbD_n^{[1]}(\theta_1,\theta_2,\hat{\theta}_3^J,\theta_1^*,\theta_2,\hat{\theta}_3^J)
\\&=&
-\half n^{-1}h\sum_{j=1}^nS(Z_\tjm,\theta_1,\hat{\theta}_3^J)^{-1}
\left[
\left(\begin{array}{c}
0
\\
2^{-1}\big(L_H(Z_\tjm,\theta_1,\theta_2,\hat{\theta}_3^J)-L_H(Z_\tjm,\theta_1^*,\theta_2,\hat{\theta}_3^J)\big)
\end{array}\right)^{\otimes2}\right].
\eeas
Therefore, 
\beas 
\sup_{(\theta_1,\theta_2)\in\overline{\Theta}_1\times\overline{\Theta}_2}\big|n^{-1}\bbD_n^{[1]}(\theta_1,\theta_2,\hat{\theta}_3^J,\theta_1^*,\theta_2,\hat{\theta}_3^J)
\big|
&=&
O_p(h).
\eeas
By definition, $\bbD_n^{[2]}(\theta_1,\hat{\theta}_3^J,\theta_1^*,\theta_2,\hat{\theta}_3^J)=0$. 
Moreover, by using the preliminary estimate $\hat{\theta}_3^J-\theta_3^*=O_p(h)$ 
{\coloro provided by Lemma \ref{201906080411},} 
{\colorb and the expression}
\beas&&
n^{-1}\bbD_n^{[3]}(\theta_1,\theta_2,\hat{\theta}_3^J,\theta_1^*,\theta_2,\hat{\theta}_3^J)
\\&=&
n^{-1}h^{1/2}\sum_{j=1}^nS(Z_\tjm,\theta_1,\theta_3^J)^{-1}
\left[
\cald_j(\theta_1^*,\theta_2,\hat{\theta}_3^J),\>
\left(\begin{array}{c}
0\\
2^{-1}\big(L(Z_\tjm,\theta_1,\theta_2,\hat{\theta}_3^J)-L(Z_\tjm,\theta_1^*,\theta_2,\hat{\theta}_3^J)\big)
\end{array}\right)\right],
\eeas
{\colorb we obtain }
\beas 
\sup_{(\theta_1,\theta_2)\in\overline{\Theta}_1\times\overline{\Theta}_2}
\big|{\cred n^{-1}}
\bbD_n^{[3]}(\theta_1,\theta_2,\hat{\theta}_3^J,\theta_1^*,\theta_2,\hat{\theta}_3^J)
\big|
&=&
O_p(h)
\eeas
{\coloro by using Lemmas \ref{201906040435}, \ref{201906031516} (b) and \ref{201906021438} (b). 
}

Now
\beas
\bbD_n^{[4]}(\theta_1,\hat{\theta}_3^J,\theta_1^*,\theta_2,\hat{\theta}_3^J)
&=&
-\half\sum_{j=1}^n\bigg\{\big(S(Z_\tjm,\theta_1,\hat{\theta}_3^J)^{-1}-S(Z_\tjm,\theta_1^*,\hat{\theta}_3^J)^{-1}\big)
\big[\cald_j(\theta_1^*,\theta_2,\hat{\theta}_3^J)^{\otimes2}\big]
\\&&
+\log\frac{\det S(Z_\tjm,\theta_1,\hat{\theta}_3^J)}{\det S(Z_\tjm,\theta_1^*,\hat{\theta}_3^J)}\bigg\}.
\eeas
Once again by using $\hat{\theta}_3^J-\theta_3^*=O_p(h)$ 
{\coloro provided by Lemma \ref{201906080411}}, 
we obtain the result {\coloro with the help of Taylor's formula and Lemma \ref{201906041850}.} 
\qed\halflineskip

{\colorr We shall deduce a tentative rough estimate 
{\colorg$o_p({\colorb n^{-1/2}}h^{-1/2})$} 
for 
the error of 
$\hat{\theta}_1^J$. }

\begin{lemma}\label{201906081400}
Suppose that $[A1]$ with $(i_A,j_A,i_B,j_B,i_H,j_H)={\coloro(1,1,2,{\colorg1},
3,2)}$, {\coloro$[A2]$} 
{\coloro and $[A3](iii')$} are satisfied. 
Then 
\beas 
n^{-1/2}h^{1/2}\>\partial_1\bbH_n(\theta_1^*,\hat{\theta}_2^J,\hat{\theta}_3^J)
&=& 
{\colorr o_p(1)}
\eeas
as $n\to\infty$.
\end{lemma}
\proof 
We use the tentative estimate of $\hat{\theta}_3^J-\theta_3^*$ given by Lemma \ref{201906080411}. 
Then 
\bea\label{201911120854} &&
n^{-1/2}h^{1/2}\>\partial_1\bbH_n(\theta_1^*,\hat{\theta}_2^J,\hat{\theta}_3^J)
\nn\\&=&
n^{-1/2}h
{\coloro \Psi_{1,1}(\theta_1^*,\hat{\theta}_2^J,\hat{\theta}_3^J,\theta_1^*,\hat{\theta}_2^J,\hat{\theta}_3^J)}
+n^{-1/2}h^{1/2}
{\coloro\Psi_{1,2}({\colorb \theta_1^*},\hat{\theta}_3^J,\theta_1^*,\hat{\theta}_2^J,\hat{\theta}_3^J)}
\nn\\&=&
O_p(n^{1/2}h){\coloro+O_p(h^{1/2})}
\\&=&
o_p(1)
\nn
\eea
since $nh^2\to0$. 
{\coloro 
In the equality (\ref{201911120854}), we used the following estimates for the second term: 
\beas &&
n^{-1/2}h^{1/2}\Psi_{1,2}({\colorb \theta_1^*},\hat{\theta}_3^J,\theta_1^*,\hat{\theta}_2^J,\hat{\theta}_3^J)
\\&=&
n^{-1/2}h^{1/2}\Psi_{1,2}({\colorb \theta_1^*},\hat{\theta}_3^J,\theta_1^*,\theta_2^*,\theta_3^*)
+O_p(n^{1/2}h)\quad(\because\text{ Lemmas \ref{201906080411}}
\text{, \ref{201906040435}, \ref{201906031516}(b) and \ref{201906021438}(b)})
\\&=&
n^{-1/2}h^{1/2}\Psi_{1,2}({\colorb \theta_1^*},\theta_3^*,\theta_1^*,\theta_2^*,\theta_3^*)
+O_p(n^{1/2}h)\quad(\because\text{ Lemmas \ref{201906080411}}
\text{ and random field argument})
\\&=&
O_p(n^{1/2}h)+O_p(h^{1/2})
\quad(\because\text{ orthogonality}).
\eeas
A similar estimate applies to the first term on (\ref{201911120854}). 
}
\qed\halflineskip

Recall 
\beas 
\Gamma_{11}
&=&
\half\int  \rm{Tr}\big\{S^{-1}(\partial_1S)S^{-1}\partial_1S(z,\theta_1^*, \theta_3^*)\big\}\nu(dz)
\\&=&
\half\int \bigg[
 \rm{Tr}\big\{\big(C^{-1}(\partial_1C)C^{-1}\partial_1C\big)(z,\theta_1^*)\big\}
 \\&&\qquad\quad
 +\rm{Tr}\big\{\big(V^{-1}H_x(\partial_1C)H_x^\star V^{-1}H_x(\partial_1C)H_x^\star\big)(z,\theta_1^*,\theta_3^*)\big\}
 \bigg]\nu(dz).
\eeas

\begin{lemma}\label{201906081402}
Suppose that $[A1]$ with $(i_A,j_A,i_B,j_B,i_H,j_H)=
(1,1,2,{\coloro3},3,{\colorr2})
$, $[A2]$,  
{\coloro $[A3](i')$ and $[A3](iii')$} 
are satisfied. 
Then, for any sequence of positive numbers $r_n$ tending to $0$,  
\bea\label{201911122311}
\sup_{\theta_1\in U(\theta_1^*,r_n)}\big|
n^{-1}\>\partial_1^2\bbH_n(\theta_1,\hat{\theta}_2^J,{\coloro\hat{\theta}_3^J})
+
\Gamma_{11}
\big|
&\to^p&
0
\eea
as $n\to\infty$.  
{\coloro In particular, }
$\hat{\theta}_1^J-\theta_1^*={\colorg o}_p(n^{-1/2}h^{-1/2})$. 
\end{lemma}
\proof 
By definition, 
\beas 
n^{-1}\>\partial_1^2\bbH_n(\theta_1,\hat{\theta}_2^J,\hat{\theta}_3^J)
&=&
-n^{-1}h{\coloro \Psi_{11,1}(\theta_1,\hat{\theta}_3^J,\theta_1,\hat{\theta}_2^J,\hat{\theta}_3^J)}
\\&&
+n^{-1}h^{1/2}
{\coloro \Psi_{11,2}(\theta_1,\hat{\theta}_3^J,\theta_1,\hat{\theta}_2^J,\hat{\theta}_3^J,\theta_1,\hat{\theta}_2^J,\hat{\theta}_3^J)}
\\&&
-\half n^{-1}{\coloro \Psi_{11,3}(\theta_1,\hat{\theta}_3^J,\theta_1,\hat{\theta}_2^J,\hat{\theta}_3^J)}
\\&&
-\half n^{-1}{\coloro \Psi_{11,4}(\theta_1,\hat{\theta}_3^J)}
\quad(\text{this term will remain})
\\&&
- n^{-1}h^{1/2}
{\coloro \Psi_{11,{\cred5}}(\theta_1,\hat{\theta}_3^J,\theta_1,\hat{\theta}_2^J,\hat{\theta}_3^J,\theta_1,\hat{\theta}_2^J,\hat{\theta}_3^J)}.
\eeas
{\coloro If we apply the same machinery as in the proof of Lemma \ref{201906081400}, 
it is easy to obtain the result. 
$[$It is remarked that $\partial_1^2$ appears in $\Psi_{11,2}$ and $\Psi_{11,3}$. 
Uniform-in-$\theta_1$ estimate for $\Psi_{11,2}$ is simple since it has the factor $h^{1/2}$ in front of it.  
On the other hand, we use random field argument for $\Psi_{11,3}$ after making the martingale differences. 
We need $\partial_1^3$ at this stage. $]$
For the second assertion, the argument becomes local by Lemma \ref{201906070327}, then 
Lemma \ref{201906081400} and the convergence (\ref{201911122311}) gives it 
by Taylor's formula. 
}
\qed\halflineskip

{\colorr
\begin{lemma}\label{201906091323}
Suppose that $[A1]$ with $(i_A,j_A,i_B,j_B,i_H,j_H)={\coloro (1,1,2,{\colorg3},3,2)}
$, {\coloro$[A2]$,} 
{\coloro $[A3](i')$ and $[A3](iii')$} are satisfied. 
Then 
\beas 
n^{-1/2}h^{1/2}\>\partial_3\bbH_n(\hat{\theta}_1^J,\hat{\theta}_2^J,\theta_3^*)
&=& 
{\colorr O_p(1)}
\eeas
as $n\to\infty$. 
In particular, $\hat{\theta}_3^J-\theta_3^*=O_p(n^{-1/2}h^{1/2})$. 
\end{lemma}
\proof 
First using an algebraic identity similar to (\ref{201906091302}), next using Lemma \ref{201906081402} 
{\coloro and once again using Lemma \ref{201906081402} with Lemma \ref{201906031516}(b)}, 
we have 
{\coloro 
\beas
n^{-1/2}\Psi_{3,1}(\hat{\theta}_1^J,\theta_3^*,\hat{\theta}_1^J,\hat{\theta}_2^J,\theta_3^*)
&=&
n^{-1/2}\Psi_{3,1}(\hat{\theta}_1^J,\theta_3^*,\hat{\theta}_1^J,\theta_2^*,\theta_3^*)
\\&=&
n^{-1/2}\Psi_{3,1}(\hat{\theta}_1^J,\theta_3^*,\theta_1^*,\theta_2^*,\theta_3^*)+O_p(1)
\\&=&
n^{-1/2}\Psi_{3,1}(\theta_1^*,\theta_3^*,\theta_1^*,\theta_2^*,\theta_3^*)+O_p(1).
\eeas
Then, from the representation of $\cald_j(\theta_1^*,\theta_2^*,\theta_3^*)$ given by 
Lemmas \ref{201906031518} and \ref{201906031516} (a) 
with the aid of the orthogonality of the martingale parts, we obtain 
\beas 
n^{-1/2}\Psi_{3,1}(\hat{\theta}_1^J,\theta_3^*,\hat{\theta}_1^J,\hat{\theta}_2^J,\theta_3^*)
&=&
O_p(1). 
\eeas
Lemmas \ref{201906031516}(b) and \ref{201906021438} easily ensures 
\beas 
n^{-1/2}
\Psi_{3,2}(\hat{\theta}_1^J,\hat{\theta}_2^J,\theta_3^*,\hat{\theta}_1^J,\hat{\theta}_2^J,\theta_3^*)
&=&
O_p(1)
\eeas
Lemmas \ref{201906031516}(b), \ref{201906021438} and \ref{201906081402} give
\beas 
n^{-1/2}h^{1/2}\Psi_{3,3}(\hat{\theta}_1^J,\theta_3^*,\hat{\theta}_1^J,\hat{\theta}_2^J,\theta_3^*)
&=&
n^{-1/2}h^{1/2}\Psi_{3,3}(\theta_1^*,\theta_3^*,\theta_1^*,\theta_2^*,\theta_3^*)
+O_p(n^{1/2}h)+O_p(1),
\eeas
and the representation of $\cald_j(\theta_1^*,\theta_2^*,\theta_3^*)$ in 
Lemmas \ref{201906031518} and \ref{201906031516} (a) and the orthogonality between 
the martingale differences, we see 
\beas 
n^{-1/2}h^{1/2}\Psi_{3,3}(\hat{\theta}_1^J,\theta_3^*,\hat{\theta}_1^J,\hat{\theta}_2^J,\theta_3^*)
&=&
O_p(1).
\eeas
Consequently, 
\beas 
n^{-1/2}h^{1/2}\>\partial_3\bbH_n(\hat{\theta}_1^J,\hat{\theta}_2^J,\theta_3^*)
&=&
n^{-1/2}
\Psi_{3,1}(\hat{\theta}_1^J,\theta_3^*,\hat{\theta}_1^J,\hat{\theta}_2^J,\theta_3^*)
+n^{-1/2}
\Psi_{3,2}(\hat{\theta}_1^J,\hat{\theta}_2^J,\theta_3^*,\hat{\theta}_1^J,\hat{\theta}_2^J,\theta_3^*)
\\&&
+{\coloro n^{-1/2}h^{1/2}
\Psi_{3,3}(\hat{\theta}_1^J,\theta_3^*,\hat{\theta}_1^J,\hat{\theta}_2^J,\theta_3^*)
}
\\&=&
O_p(1).
\eeas
}
\begin{en-text}
\beas &&
n^{-1/2}h^{1/2}\>\partial_3\bbH_n(\hat{\theta}_1^J,\hat{\theta}_2^J,\theta_3^*)
\\&=&
n^{-1/2}
\Psi_{3,1}(\hat{\theta}_1^J,\theta_3^*,\hat{\theta}_1^J,\hat{\theta}_2^J,\theta_3^*)
+n^{-1/2}
\Psi_{3,2}(\hat{\theta}_1^J,\hat{\theta}_2^J,\theta_3^*,\hat{\theta}_1^J,\hat{\theta}_2^J,\theta_3^*)
\\&&
+{\coloro n^{-1/2}h^{1/2}
\Psi_{3,3}(\hat{\theta}_1^J,\theta_3^*,\hat{\theta}_1^J,\hat{\theta}_2^J,\theta_3^*)
}
\\&=&
n^{-1/2}
\Psi_{3,1}(\hat{\theta}_1^J,\theta_3^*,\hat{\theta}_1^J,\theta_2^*,\theta_3^*)%
+n^{-1/2}\Psi_{3,2}(\hat{\theta}_1^J,\hat{\theta}_2^J,\theta_3^*,\hat{\theta}_1^J,\hat{\theta}_2^J,\theta_3^*)
\\&&
+\half n^{-1/2}h^{1/2}\sum_{j=1}^n
\big(S^{-1}(\partial_3S))S^{-1}\big)(Z_\tjm,\hat{\theta}_1^J,\theta_3^*)
\big[\cald_j(\hat{\theta}_1^J,\hat{\theta}_2^J,\theta_3^*)^{\otimes2}-\hat{S}(Z_\tjm,\theta_3^*)\big]
\\&=&
n^{-1/2}
\Psi_{3,1}(\theta_1^*,\theta_3^*,\theta_1^*,\theta_2^*,\theta_3^*)+O_p(1)
\\&&
+n^{-1/2}h\sum_{j=1}^nS(Z_\tjm,\theta_1^*,\theta_3^*)^{-1}
\left[\cald_j(\theta_1^*,\theta_2^*,\theta_3^*),\>
\left(\begin{array}{c}0\\
2^{-1}\partial_3L_H(Z_\tjm,\hat{\theta}_1^J,\hat{\theta}_2^J,\theta_3^*)
\end{array}\right)
\right]+o_p(h^{1/2})
\\&&
+\half n^{-1/2}h^{1/2}\sum_{j=1}^n
\big(S^{-1}(\partial_3S))S^{-1}\big)(Z_\tjm,\theta_1^*,\theta_3^*)
\big[\cald_j(\theta_1^*,\theta_2^*,\theta_3^*)^{\otimes2}-S(Z_\tjm,\theta_1^*,\theta_3^*)\big]
+o_p(1)
\\&=&
O_p(1).
\eeas
\end{en-text}
For the last assertion, we may apply Lemma \ref{201906080411}. 
\qed\halflineskip
}
\begin{en-text}
Next, 
\beas&&
n^{-1}\big\{\bbH_n(\theta_1,\theta_2,\hat{\theta}_3)-\bbH_n(\theta_1^*,\theta_2,\hat{\theta}_3)\big\}
\eeas
\beas&&
n^{-1}\big\{\bbH_n(\theta_1,\theta_2,\hat{\theta}_3)-\bbH_n(\theta_1^*,\theta_2,\hat{\theta}_3)\big\}
\\&=&
-\half n^{-1}\sum_{j=1}^nS(Z_\tjm,\theta_1,\hat{\theta}_3)^{-1}\big[\big(\cald_j(\theta_1,\theta_2,\hat{\theta}_3)-\cald_j(\theta_1^*,\theta_2,\hat{\theta}_3)\big)^{\otimes2}\big]
\\&&
+n^{-1}h^{1/2}\sum_{j=1}^nS(Z_\tjm,\theta_1,\hat{\theta}_3)^{-1}
\left[
\cald_j(\theta_1^*,\theta_2,\hat{\theta}_3),\>
\left(\begin{array}{c}
0\\
2^{-1}\big(L(Z_\tjm,\theta_1,\theta_2,\hat{\theta}_3)-L(Z_\tjm,\theta_1^*,\theta_2,\hat{\theta}_3)\big)
\end{array}\right)\right]
\\&&
-\half n^{-1}\sum_{j=1}^n\big\{S(Z_\tjm,\theta_1,\hat{\theta}_3)^{-1}-S(Z_\tjm,\theta_1^*,\hat{\theta}_3)^{-1}\big\}
\big[\cald_j(\theta_1^*,\theta_2,\hat{\theta}_3)^{\otimes2}\big]
\\&&
-\half n^{-1}\sum_{j=1}^n\log\frac{\log\det S(Z_\tjm,\theta_1,\hat{\theta}_3)}{\log\det S(Z_\tjm,\theta_1^*,\hat{\theta}_3)}
\\&\to^p&
-\half \int \nu(dz).
\eeas
\end{en-text}

{\colorr 

Recall
\beas 
M_n^{(1)} &=& 
\half n^{-1/2}\sum_{j=1}^n
\big(S^{-1}(\partial_1S)S^{-1}\big)(Z_\tjm,\theta_1^*,\theta_3^*)
\big[\widetilde{\cald}_j(\theta_1^*,\theta_2^*,\theta_3^*)^{\otimes2}-S(Z_\tjm,\theta_1^*,\theta_3^*)\big]. 
\eeas

\begin{lemma}\label{201906091405}
Suppose that $[A1]$ with $(i_A,j_A,i_B,j_B,i_H,j_H)=
{\coloro(1,1,2,{\colorg3},3,2)}
$, {\coloro$[A2],$} 
{\coloro $[A3](i')$ and $[A3](iii')$} are satisfied. 
Then 
\beas 
n^{-1/2}\>\partial_1\bbH_n(\theta_1^*,\hat{\theta}_2^J,\hat{\theta}_3^J)
&=& 
{\colorr O_p(1)}
\eeas
as $n\to\infty$. 
Moreover, 
\beas 
n^{1/2}\big(\hat{\theta}_1^J-\theta_1^*\big)
-\Gamma_{11}
^{-1}M_n^{(1)}
&\to^p& 
0
\eeas
as $n\to\infty$. 
In particular, $\hat{\theta}_1^J-\theta_1^*=O_p(n^{-1/2})$ as $n\to\infty$. 
\end{lemma}
\proof 
We are in the same situation as Lemma \ref{201906081400} but 
we can use the convergence rate $\hat{\theta}_3^J-\theta_3^*=O_p(n^{-1/2}h^{1/2})$ 
elaborated by Lemma \ref{201906091323}.  
Then 
{\coloro 
\beas &&
n^{-1/2}\>\partial_1\bbH_n(\theta_1^*,\hat{\theta}_2^J,\hat{\theta}_3^J)
\\&=&
n^{-1/2}h^{1/2}
\Psi_{1,1}(\theta_1^*,\hat{\theta}_2^J,\hat{\theta}_3^J,\theta_1^*,\hat{\theta}_2^J,\hat{\theta}_3^J)
+n^{-1/2}
\Psi_{1,2}(\theta_1^*,\hat{\theta}_3^J,\theta_1^*,\hat{\theta}_2^J,\hat{\theta}_3^J)
\\&=&
n^{-1/2}h^{1/2}\Psi_{1,1}(\theta_1^*,\theta_2^*,\theta_3^*,\theta_1^*,\theta_2^*,\theta_3^*)
+O_p(n^{1/2}h)+O_p(h^{1/2})
\quad(\text{Lemmas \ref{201906031516}(b), \ref{201906040435} and \ref{201906091323}})
\\&&
+n^{-1/2}\Psi_{1,2}(\theta_1^*,\theta_3^*,\theta_1^*,\hat{\theta}_2^J,\theta_3^*)+O_p(1)
\quad(\text{Lemmas \ref{201906040435}, \ref{201906021438}(b) and \ref{201906091323}
}
)
\eeas
{\colorg 
For the last term, we can use the decomposition 
\beas 
\cald_j(\theta_1^*,\hat{\theta}_2^J,\hat{\theta}_3^J)^{\otimes2}-\cald_j(\theta_1^*,\hat{\theta}_2^J,\theta_3^*)^{\otimes2}
&=&
\big\{\cald_j(\theta_1^*,\hat{\theta}_2^J,\hat{\theta}_3^J)-\cald_j(\theta_1^*,\hat{\theta}_2^J,\theta_3^*)\big\}^{\otimes2}
\\&&
+2\bigg\{\cald_j(\theta_1^*,\theta_2^*,\theta_3^*)
+\big(\cald_j(\theta_1^*,\hat{\theta}_2^J,\theta_3^*)-\cald_j(\theta_1^*,\theta_2^*,\theta_3^*)\big)\bigg\}
\\&&
\qquad{\cred\otimes_{sym}}
\big\{\cald_j(\theta_1^*,\hat{\theta}_2^J,\hat{\theta}_3^J)-\cald_j(\theta_1^*,\hat{\theta}_2^J,\theta_3^*)\big\},
\eeas
{\cred 
where $\otimes_{sym}$ means the symmetrized tensor product.}
}

We have 
\beas &&
\cald_j(\theta_1^*,\hat{\theta}_2^J,\theta_3^*)^{\otimes2}
-\cald_j(\theta_1^*,\theta_2^*,\theta_3^*)^{\otimes2}
\\&=&
2\big(\cald_j(\theta_1^*,\hat{\theta}_2^J,\theta_3^*)-\cald_j(\theta_1^*,\theta_2^*,\theta_3^*)\big)\otimes_{sym}\cald_j(\theta_1^*,\theta_2^*,\theta_3^*)
+\big(\cald_j(\theta_1^*,\hat{\theta}_2^J,\theta_3^*)-\cald_j(\theta_1^*,\theta_2^*,\theta_3^*)\big)^{\otimes2}
\\&=&
2h^{1/2}\left(\begin{array}{c}A(Z_\tjm,\theta_2^*)-A(Z_\tjm,\hat{\theta}_2^J)
\\ 2^{-1}{\cred\big(}
L_H(Z_\tjm,\theta_1^*,\theta_2^*,\theta_3^*)-L_H(Z_\tjm,\theta_1^*,\hat{\theta}_2^J,\theta_3^*)
{\cred\big)}
\end{array}\right)\otimes_{sym}\cald_j(\theta_1^*,\theta_2^*,\theta_3^*)+
O_{L^\inftym}(h)
\eeas
To estimate $n^{-1/2}\Psi_{1,2}(\theta_1^*,\theta_3^*,\theta_1^*,\hat{\theta}_2^J,\theta_3^*)$, 
we introduce the random field
\bea\label{201911171819} 
\Xi_n^{(\ref{201911171819})}(\theta_2) 
&=& 
n^{-1/2}h^{1/2}
\sum_{j=1}^n
\big(S^{-1}(\partial_1S))S^{-1}\big)(Z_\tjm,\theta_1^*,\theta_3^*)
\nn\\&&
\cdot
\bigg[
\left(\begin{array}{c}A(Z_\tjm,\theta_2^*)-A(Z_\tjm,\theta_2)
\nn\\ 2^{-1}{\cred\big(}
L_H(Z_\tjm,\theta_1^*,\theta_2^*,\theta_3^*)-L_H(Z_\tjm,\theta_1^*,\theta_2,\theta_3^*)
{\cred\big)}
\end{array}\right)\otimes_{sym}\cald_j(\theta_1^*,\theta_2^*,\theta_3^*)\bigg].
\nn\\&&
\eea
With the aid of the representation of $\cald_j(\theta_1^*,\theta_2^*,\theta_3^*)$ and 
the orthogonality between martingale differences, a random field argument concludes 
\beas 
\sup_{\theta_2\in\Theta_2}|\Xi_n^{(\ref{201911171819})}(\theta_2) | &=& O_p(h^{1/2}),
\eeas
in particular, 
\beas 
n^{-1/2}\Psi_{1,2}(\theta_1^*,\theta_3^*,\theta_1^*,\hat{\theta}_2^J,\theta_3^*)
&=&
n^{-1/2}\Psi_{1,2}(\theta_1^*,\theta_3^*,\theta_1^*,\theta_2^*,\theta_3^*)
+o_p(1).
\eeas

The orthogonality further applied gives 
\beas 
n^{-1/2}h^{1/2}\Psi_{1,1}(\theta_1^*,\theta_2^*,\theta_3^*,\theta_1^*,\theta_2^*,\theta_3^*)
&=&
O_p(h^{1/2}). 
\eeas
Consequently, 
\bea\label{201911171552}
n^{-1/2}\>\partial_1\bbH_n(\theta_1^*,\hat{\theta}_2^J,\hat{\theta}_3^J)
&=& 
n^{-1/2}\Psi_{1,2}(\theta_1^*,\theta_3^*,\theta_1^*,\theta_2^*,\theta_3^*)
+o_p(1)
\nn\\&=&
M_n^{(1)}+o_p(1)
\\&=&
\nn O_p(1)
\eea
as $n\to\infty$. 

Since $\hat{\theta}_1^J\to^p\theta_1^*$ by e.g. Lemma \ref{201906070327}, 
we can show the first order efficiency of $\hat{\theta}_1^J$ by using Taylor's formula 
combined with (\ref{201911171552}) and Lemma \ref{201906081402}. 
}
\begin{en-text}
\beas&&
\\&=&
n^{-1/2}h^{1/2}\sum_{j=1}^nS(Z_\tjm,\theta_1^*,\hat{\theta}_3^J)^{-1}
\left[\cald_j(\theta_1^*,\hat{\theta}_2^J,\hat{\theta}_3^J),\>
\left(\begin{array}{c}0\\
2^{-1}\partial_1L_H(Z_\tjm,\theta_1^*,\hat{\theta}_2^J,\hat{\theta}_3^J)
\end{array}\right)
\right]
\\&&
+\half n^{-1/2}\sum_{j=1}^n
\big(S^{-1}(\partial_1S))S^{-1}\big)(Z_\tjm,\hat{\theta}_1^J,\hat{\theta}_3^J)
\big[\cald_j(\theta_1^*,\hat{\theta}_2^J,\hat{\theta}_3^J)^{\otimes2}-S(Z_\tjm,\theta_1^*,\hat{\theta}_3^J)\big]
\\&=&
n^{-1/2}h^{1/2}\sum_{j=1}^nS(Z_\tjm,\theta_1^*,\theta_3^*)^{-1}
\left[\cald_j(\theta_1^*,\hat{\theta}_2^J,\theta_3^*),\>
\left(\begin{array}{c}0\\
2^{-1}\partial_1L_H(Z_\tjm,\theta_1^*,\hat{\theta}_2^J,\theta_3^*)
\end{array}\right)
\right]+O_p(h^{1/2})
\\&&
+\half n^{-1/2}\sum_{j=1}^n
\big(S^{-1}(\partial_1S))S^{-1}\big)(Z_\tjm,\hat{\theta}_1^J,\theta_3^*)
\big[\cald_j(\theta_1^*,\hat{\theta}_2^J,\theta_3^*)^{\otimes2}-S(Z_\tjm,\theta_1^*,\theta_3^*)\big]
+O_p(1)
\\&=&
n^{-1/2}h^{1/2}\sum_{j=1}^nS(Z_\tjm,\theta_1^*,\theta_3^*)^{-1}
\left[\cald_j(\theta_1^*,\theta_2^*,\theta_3^*),\>
\left(\begin{array}{c}0\\
2^{-1}\partial_1L_H(Z_\tjm,\theta_1^*,\hat{\theta}_2^J,\theta_3^*)
\end{array}\right)
\right]
\\&&+O_p(h^{1/2})+O_p(n^{1/2}h)
\\&&
+\half n^{-1/2}\sum_{j=1}^n
\big(S^{-1}(\partial_1S))S^{-1}\big)(Z_\tjm,\hat{\theta}_1^J,\theta_3^*)
\big[\cald_j(\theta_1^*,\theta_2^*,\theta_3^*)^{\otimes2}-S(Z_\tjm,\theta_1^*,\theta_3^*)\big]
\\&&
+O_p(1)+O_p(n^{1/2}h)
\\&=&
\half n^{-1/2}\sum_{j=1}^n
\big(S^{-1}(\partial_1S))S^{-1}\big)(Z_\tjm,\hat{\theta}_1^J,\theta_3^*)
\big[\widetilde{\cald}_j(\theta_1^*,\theta_2^*,\theta_3^*)^{\otimes2}-S(Z_\tjm,\theta_1^*,\theta_3^*)\big]
\\&&
+O_p(1)
\eeas
since $nh^2\to0$. 
Fandom field argument works to replace $\hat{\theta}_1^J$ by $\theta_1^*$ 
for the first term of the right-hand side. 
Thus, we obtain the first order efficiency. 
Then, totally, we know the last expression if of $O_p(1)$. 
Now Lemma \ref{201906081402} gives $\hat{\theta}_1^J-\theta_1^*=O_p(n^{-1/2})$. 
\end{en-text}
\qed\halflineskip

\begin{lemma}\label{201906091542}
Suppose that $[A1]$ with $(i_A,j_A,i_B,j_B,i_H,j_H)=(1,1,2,{\cred3},3,2)
$, {\coloro$[A2]$}, 
{\coloro $[A3](i')$ and $[A3](iii')$} 
are satisfied. 
Then 
\bea\label{201911171901} 
n^{-1/2}h^{1/2}\>\partial_3\bbH_n(\hat{\theta}_1^J,\hat{\theta}_2^J,\theta_3^*)
-M_n^{(3)}
&\to^p& 
0
\eea
as $n\to\infty$. 
In particular, 
\bea\label{201911171902} 
n^{1/2}h^{-1/2}\big(\hat{\theta}_3^J-\theta_3^*\big)
&\to^d&
N(0,\Gamma_{33}^{-1})
\eea
as $n\to\infty$. 
\end{lemma}
\proof 
We elaborate the estimate in the proof of Lemma \ref{201906091323}. 
Taking advantage of the convergence rate of $\hat{\theta}_1^J$ given by Lemma \ref{201906091405}, 
we see 
{\coloro
\beas&&
n^{-1/2}h^{1/2}\>\partial_3\bbH_n(\hat{\theta}_1^J,\hat{\theta}_2^J,\theta_3^*)
\\&=&
n^{-1/2}\Psi_{3,1}(\hat{\theta}_1^J,\theta_3^*,\hat{\theta}_1^J,{\colorb \theta_2^*},\theta_3^*)
+n^{-1/2}\Psi_{3,2}(\hat{\theta}_1^J,\hat{\theta}_2^J,\theta_3^*,\hat{\theta}_1^J,\hat{\theta}_2^J,\theta_3^*)
\\&&
+n^{-1/2}h^{1/2}\Psi_{3,3}(\hat{\theta}_1^J,\theta_3^*,\hat{\theta}_1^J,\hat{\theta}_2^J,\theta_3^*)
\\&=&
n^{-1/2}\Psi_{3,1}(\hat{\theta}_1^J,\theta_3^*,\theta_1^*,{\colorb \theta_2^*},\theta_3^*)+O_p(h^{1/2})
\quad(\text{Lemmas \ref{201906091405} and \ref{201906031516}(b)})
\\&&
+n^{-1/2}\Psi_{3,2}(\hat{\theta}_1^J,\hat{\theta}_2^J,\theta_3^*,\theta_1^*,\theta_2^*,\theta_3^*)
+o_p(h^{1/2})
\quad(\text{Lemma \ref{201906031516}(b)})
\\&&
+n^{-1/2}h^{1/2}\Psi_{3,3}(\hat{\theta}_1^J,\theta_3^*,\theta_1^*,\theta_2^*,\theta_3^*)+o_p(1)
\quad(\text{Lemmas \ref{201906031516}(b) and \ref{201906021438}(b)}).
\eeas
\begin{en-text}
Define $\Xi_n^{(\ref{201911171821})}(\theta_1)$ by 
\bea\label{201911171821}
\Xi_n^{(\ref{201911171821})}(\theta_1)
&=&
n^{-1/2}h^{1/2}\Psi_{3,3}(\theta_1,\theta_3^*,\theta_1^*,\theta_2^*,\theta_3^*)
\eea
Then we see 
$\sup_{\theta_1\in\Theta_1}|\Xi_n^{(\ref{201911171821})}(\theta_1)|=o_p(1)$ 
by a random field argument with Sobolev's inequality, the representation of 
$\cald_j(\theta_1^*,\theta_2^*,\theta_3^*)$ and the orthogonality. 
In particular, 
\end{en-text}
By Lemma \ref{201906091405}, the representation of 
$\cald_j(\theta_1^*,\theta_2^*,\theta_3^*)$ and the orthogonality, we obtain 
\beas 
n^{-1/2}h^{1/2}\Psi_{3,3}(\hat{\theta}_1^J,\theta_3^*,\theta_1^*,\theta_2^*,\theta_3^*)
&=&
n^{-1/2}h^{1/2}\Psi_{3,3}(\theta_1^*,\theta_3^*,\theta_1^*,\theta_2^*,\theta_3^*)
+O_p(h^{1/2})
\\&=&
O_p(h^{1/2}). 
\eeas
\begin{en-text}
For the random field 
\bea\label{201911171831}
\Xi_n^{(\ref{201911171831})}(\theta_1,\theta_2)
&=&
n^{-1/2}\Psi_{3,2}(\theta_1,\theta_2,\theta_3^*,\theta_1^*,\theta_2^*,\theta_3^*), 
\eea
we obtain 
\beas 
\sup_{(\theta_1,\theta_2)\in\Theta_1\times\Theta_2}\big|\Xi_n^{(\ref{201911171831})}(\theta_1,\theta_2)\big|
&=&
O_p(h^{1/2}).
\eeas
\end{en-text}
{\cred 
We have
\beas 
\sup_{(\theta_1,\theta_2)\in\Theta_1\times\Theta_2}\big|
n^{-1/2}\Psi_{3,2}(\theta_1,\theta_2,\theta_3^*,\theta_1^*,\theta_2^*,\theta_3^*)\big|
&=&
O_p(h).
\eeas
}
Next, we consider 
\bea\label{201911171830}
\Xi_n^{(\ref{201911171830})}(u_1)
&=&
n^{-1/2}\bigg\{\Psi_{3,1}(\theta_1^*+r_nu_1,\theta_3^*,\theta_1^*,{\colorb \theta_2^*},\theta_3^*)
-
\Psi_{3,1}(\theta_1^*,\theta_3^*,\theta_1^*,{\colorb \theta_2^*},\theta_3^*)\bigg\}
\eea
for any sequence $r_n$ of positive numbers such that $r_n\to0$. 
Then a random field argument with Sobolev's inequality ensures the convergence 
\beas 
\sup_{u_1\in B(0.1)}|\Xi_n^{(\ref{201911171830})}(u_1)|
&=& o_p(1). 
\eeas
Therefore, 
\beas 
n^{-1/2}\Psi_{3,1}(\hat{\theta}_1^J,\theta_3^*,\theta_1^*,{\colorb \theta_2^*},\theta_3^*)
&=&
n^{-1/2}\Psi_{3,1}(\theta_1^*,\theta_3^*,\theta_1^*,{\colorb \theta_2^*},\theta_3^*)+o_p(1)
\\&=&
M_n^{(3)}+o_p(1).
\eeas

From the above estimates, we already have (\ref{201911171901}). 
Moreover, Lemmas \ref{201906080411} and the martingale central limit theorem givens (\ref{201911171902}). 
}
\begin{en-text}
\beas &&
n^{-1/2}h^{1/2}\>\partial_3\bbH_n(\hat{\theta}_1^J,\hat{\theta}_2^J,\theta_3^*)
\\&=&
n^{-1/2}\sum_{j=1}^nS(Z_\tjm,\hat{\theta}_1^J,\theta_3^*)^{-1}
\left[\cald_j(\hat{\theta}_1^J,{\colorb \theta_2^*},\theta_3^*),\>
\left(\begin{array}{c}0\\
\partial_3H(Z_\tjm,\theta_3^*)
\end{array}\right)
\right]
\\&&
+n^{-1/2}h\sum_{j=1}^nS(Z_\tjm,\hat{\theta}_1^J,\theta_3^*)^{-1}
\left[\cald_j(\hat{\theta}_1^J,\hat{\theta}_2^J,\theta_3^*),\>
\left(\begin{array}{c}0\\
2^{-1}\partial_3L_H(Z_\tjm,\hat{\theta}_1^J,\hat{\theta}_2^J,\theta_3^*)
\end{array}\right)
\right]
\\&&
+\half n^{-1/2}h^{1/2}\sum_{j=1}^n
\big(S^{-1}(\partial_3S))S^{-1}\big)(Z_\tjm,\hat{\theta}_1^J,\theta_3^*)
\big[\cald_j(\hat{\theta}_1^J,\hat{\theta}_2^J,\theta_3^*)^{\otimes2}-\hat{S}(Z_\tjm,\theta_3^*)\big]
\\&=&
n^{-1/2}\sum_{j=1}^nS(Z_\tjm,\theta_1^*,\theta_3^*)^{-1}
\left[\cald_j(\theta_1^*,{\colorr \theta_2^*},\theta_3^*),\>
\left(\begin{array}{c}0\\
\partial_3H(Z_\tjm,\theta_3^*)
\end{array}\right)
\right]+{\colorb O_p(h^{1/2})}
\\&&
+n^{-1/2}h\sum_{j=1}^nS(Z_\tjm,\theta_1^*,\theta_3^*)^{-1}
\left[\cald_j(\theta_1^*,\theta_2^*,\theta_3^*),\>
\left(\begin{array}{c}0\\
2^{-1}\partial_3L_H(Z_\tjm,\hat{\theta}_1^J,\hat{\theta}_2^J,\theta_3^*)
\end{array}\right)
\right]+o_p(h^{1/2})
\\&&
+\half n^{-1/2}h^{1/2}\sum_{j=1}^n
\big(S^{-1}(\partial_3S))S^{-1}\big)(Z_\tjm,\theta_1^*,\theta_3^*)
\big[\cald_j(\theta_1^*,\theta_2^*,\theta_3^*)^{\otimes2}-S(Z_\tjm,\theta_1^*,\theta_3^*)\big]
+o_p(1)
\\&=&
M_n^{(3)}+o_p(1). 
\eeas
\end{en-text}
\qed\halflineskip
}

\begin{lemma}\label{201906080713}
Suppose that $[A1]$ with $(i_A,j_A,i_B,j_B,i_H,j_H)=(1,1,2,{\cred3},3,2)
$, $[A2]$, {\coloro $[A3](i')$ and $[A3](iii')$} are satisfied. 
Then 
\bd\im[(a)] As $n\to\infty$, 
\bea
\sup_{\theta_2\in \overline{\Theta}_2}\big|
n^{-1}h^{-1}\big\{\bbH_n(\hat{\theta}_1^J,\theta_2,\hat{\theta}_3^J)-\bbH_n(\hat{\theta}_1^J,\theta_2^*,\hat{\theta}_3^J)\big\}
-\bbY^{(2)}(\theta_2)
\big|
\>\to^p\>0
\eea

\im[(b)] If $[A3]$ $(ii)$ is satisfied, then $\hat{\theta}^J_2\to^p\theta_2^*$ as $n\to\infty$. 
\ed
\end{lemma}
\proof 
We have 
\beas &&
n^{-1}h^{-1}\big\{\bbH_n(\hat{\theta}_1^J,\theta_2,\hat{\theta}_3^J)-\bbH_n(\hat{\theta}_1^J,\theta_2^*,\hat{\theta}_3^J)\big\}
\\&=&
n^{-1}h^{-1}\bbD_n(\hat{\theta}_1^J,\theta_2,\hat{\theta}_3^J,\hat{\theta}_1^J,\theta_2^*,\hat{\theta}_3^J)
\\&=&
n^{-1}h^{-1}\bbD_n^{[1]}(\hat{\theta}_1^J,\theta_2,\hat{\theta}_3^J,\hat{\theta}_1^J,\theta_2^*,\hat{\theta}_3^J)
+n^{-1}h^{-1}\bbD_n^{[2]}(\hat{\theta}_1^J,\hat{\theta}_3^J,\hat{\theta}_1^J,\theta_2^*,\hat{\theta}_3^J)
\nn\\&&
+n^{-1}h^{-1}\bbD_n^{[3]}(\hat{\theta}_1^J,\theta_2,\hat{\theta}_3^J,\hat{\theta}_1^J,\theta_2^*,\hat{\theta}_3^J)
+n^{-1}h^{-1}\bbD_n^{[4]}(\hat{\theta}_1^J,\hat{\theta}_3^J,\hat{\theta}_1^J,\theta_2^*,\hat{\theta}_3^J).
\eeas

We have 
\beas &&
n^{-1}h^{-1}\bbD_n^{[1]}(\hat{\theta}_1^J,\theta_2,\hat{\theta}_3^J,\hat{\theta}_1^J,\theta_2^*,\hat{\theta}_3^J)
\\&=&
-n^{-1}h^{-1}\half \sum_{j=1}^nS(Z_\tjm,\hat{\theta}_1^J,\hat{\theta}_3^J)^{-1}\big[\big(\cald_j(\hat{\theta}_1^J,\theta_2,\hat{\theta}_3^J)-\cald_j(\hat{\theta}_1^J,\theta_2^*,\hat{\theta}_3^J)\big)^{\otimes2}\big]
\\&=&
-\half n^{-1}\sum_{j=1}^nS(Z_\tjm,\hat{\theta}_1^J,\hat{\theta}_3^J)^{-1}
\left[
\left(\begin{array}{c}
\big(A(Z_\tjm,\theta_2)-A(Z_\tjm,\theta_2^*)\big)
\\
2^{-1}\big(L_H(Z_\tjm,\hat{\theta}_1^J,\theta_2,\hat{\theta}_3^J)-L_H(Z_\tjm,\hat{\theta}_1^J,\theta_2^*,\hat{\theta}_3^J)\big)
\end{array}\right)^{\otimes2}\right]
\\&=&
-\half n^{-1}\sum_{j=1}^nS(Z_\tjm,\hat{\theta}_1^J,\hat{\theta}_3^J)^{-1}
\left[
\left(\begin{array}{c}
\big(A(Z_\tjm,\theta_2)-A(Z_\tjm,\theta_2^*)\big)
\\
2^{-1}H_x(Z_\tjm,\hat{\theta}_3^J)\big[A(Z_\tjm,\theta_2)-A(Z_\tjm,\theta_2^*)\big]
\end{array}\right)^{\otimes2}\right]
\\&=&
-\half n^{-1}\sum_{j=1}^n
C(Z_\tjm,\hat{\theta}_1^J)^{{\cred -1}}
\big[\big(A(Z_\tjm,\theta_2)-A(Z_\tjm,\theta_2^*)\big)^{\otimes2}\big],
\eeas
$\bbD_n^{[2]}(\hat{\theta}_1^J,\hat{\theta}_3^J,\hat{\theta}_1^J,\theta_2^*,\hat{\theta}_3^J)
=0$, 
$\bbD_n^{[4]}(\hat{\theta}_1^J,\hat{\theta}_3^J,\hat{\theta}_1^J,\theta_2^*,\hat{\theta}_3^J)
=0$
and 
{\colorr 
\beas &&
n^{-1}h^{-1}\bbD_n^{[3]}(\hat{\theta}_1^J,\theta_2,\hat{\theta}_3^J,\hat{\theta}_1^J,\theta_2^*,\hat{\theta}_3^J)
\yeq
n^{-1}h^{-1}\bbD_n^{[3]}(\hat{\theta}_1^J,\theta_2,\hat{\theta}_3^J,\theta_1^*,\theta_2^*,\theta_3^*)
+o_p(1)
\eeas
by Lemmas \ref{201906091323}, {\cred \ref{201906091405}, \ref{201906040435} and 
\ref{201906031516} (b),}
{\coloro where the order $o_p(1)$ is uniform in $\theta_2\in\Theta_2$.} 
The last expression is 
\beas &&
n^{-1}h^{-1}\bbD_n^{[3]}(\theta_1^*,\theta_2,\theta_3^*,\theta_1^*,\theta_2^*,\theta_3^*)
+o_p(1)
\eeas
by using the exact convergence rate of $\hat{\theta}_1^J$ and $\hat{\theta}_3^J$, 
{\coloro where $o_p(1)$ is uniform in $\theta_2\in\Theta_2$.} 
Random field argument shows that the last one converges in probability {\coloro to zero} 
uniformly in $\theta_2$. 
{\coloro This shows (a). The property (b) is now easy to deduce from (a).} 
\qed\halflineskip
}

\begin{en-text}
{\colorr 
\begin{fbox}
{$\up$ $\hat{\theta}_3^J$のレートが足りない．}
\end{fbox}
}
\end{en-text}

{\colorr 
We will derive a convergence rate of $\hat{\theta}_2^J$. 

\begin{lemma}\label{201906081515}
Suppose that $[A1]$ with $(i_A,j_A,i_B,j_B,i_H,j_H)=(1,1,2,{\cred1},3,2)
$, $[A2]$, {\coloro $[A3](i')$ and $[A3](iii')$} are satisfied. 
Then 
\beas 
n^{-1/2}h^{-1/2}\partial_2\bbH_n(\hat{\theta}_1^J,\theta_2^*,\hat{\theta}_3^J)
-M_n^{(2)}
&\to^p&
0
\eeas
as $n\to\infty$. 
\end{lemma}
\proof
{\coloro By simple algebra {\cred and Lemma \ref{201906041900}}, }
\beas &&
n^{-1/2}h^{-1/2}\partial_2\bbH_n(\hat{\theta}_1^J,\theta_2^*,\hat{\theta}_3^J)
\\&=&
n^{-1/2}\sum_{j=1}^nS(Z_\tjm,\hat{\theta}_1^J,\hat{\theta}_3^J)^{-1}
\left[\cald_j(\hat{\theta}_1^J,\theta_2^*,\hat{\theta}_3^J),\>
\left(\begin{array}{c} \partial_2A(Z_\tjm,\theta_2^*)\\
2^{-1}H_x(Z_\tjm,\hat{\theta}_3^J)[\partial_2A(Z_\tjm,\theta_2^*)]
\end{array}\right)
\right]
\\&=&
n^{-1/2}\sum_{j=1}^nC(Z_\tjm,\hat{\theta}_1^J)^{-1}\big[\partial_2A(Z_\tjm,\theta_2^*),\>
h^{-1/2}(\Delta_jX-hA(Z_\tjm,\theta_2^*))\big]
\\&=&
n^{-1/2}\sum_{j=1}^nC(Z_\tjm,\hat{\theta}_1^J)^{-1}\big[\partial_2A(Z_\tjm,\theta_2^*),\>
h^{-1/2}B(Z_\tjm,\theta_1^*)\Delta_jw\big]+o_p(1)
\\&=&
M_n^{(2)}+o_p(1).
\eeas
{\coloro Here the last equation can be verified by a $\theta_1$-random field argument 
using the consistency of $\hat{\theta}_1^J$ obtained in Lemma \ref{201906070327}. 
{\cred Remark that $M^{(2)}_n$ is defined by (\ref{202001281745}) on p.\pageref{202001281745}.}
}
\qed\halflineskip

\begin{lemma}\label{201906091659}
Suppose that $[A1]$ with $(i_A,j_A,i_B,j_B,i_H,j_H)=
(1,{\cred 3},2,{\coloro3},3,{\colorr2})
$, $[A2]$, 
{\coloro $[A3](i')$, $[A3](ii)$ and $[A3](iii')$} are satisfied. 
Then 
\bea\label{201911180914} 
\sup_{\theta_2\in U(\theta_2^*,r_n)}\big|
n^{-1}h^{-1}\>\partial_2^2\bbH_n(\hat{\theta}_1^J,\theta_2,\hat{\theta}_3^J)
+
\Gamma_{22}
\big|
&\to^p&
0
\eea
as $n\to\infty$, where $r_n$ is any sequence of positive numbers such that $r_n\to0$ and 
\beas 
\Gamma_{22}
&=&
\int S(z,\theta_1^*,\theta_3^*)^{-1}
\left[
\left(\begin{array}{c}\partial_2A(z,\theta_2^*)\\
2^{-1}\partial_2L_H(z,\theta_1^*,\theta_2^*,\theta_3^*)
\end{array}\right)^{\otimes2}
\right]\nu(dz)
\\&=&
\int C(z,\theta_1^*)^{-1}\big[(\partial_2A(z,\theta_2^*)^{\otimes2}\big]\nu(dz).
\eeas
Moreover, 
\bea\label{201911180915} 
n^{1/2}h^{1/2}\big(\hat{\theta}_2^J-\theta_2^*\big)
-\Gamma_{22}^{-1}M_n^{(2)}
&\to^p&
0
\eea
as $n\to\infty$. In particular, 
\beas 
n^{1/2}h^{1/2}\big(\hat{\theta}_2^J-\theta_2^*\big)
&\to^d&
N(0,\Gamma_{22}^{-1})
\eeas
as $n\to\infty$. 
\end{lemma}
\proof We see
\beas &&
n^{-1}h^{-1}\partial_2^2\bbH_n(\hat{\theta}_1^J,\theta_2,\hat{\theta}_3^J)
\\&=&
-n^{-1}\sum_{j=1}^nS(Z_\tjm,\hat{\theta}_1^J,\hat{\theta}_3^J)^{-1}
\left[
\left(\begin{array}{c} \partial_2A(Z_\tjm,\theta_2)\\
2^{-1}H_x(Z_\tjm,\hat{\theta}_3^J)[\partial_2A(Z_\tjm,\theta_2)]
\end{array}\right)^{\otimes2}
\right]
\\&&
+n^{-1}h^{-1/2}\sum_{j=1}^nS(Z_\tjm,\hat{\theta}_1^J,\hat{\theta}_3^J)^{-1}
\left[\cald_j(\hat{\theta}_1^J,\theta_2,\hat{\theta}_3^J),\>
\left(\begin{array}{c} \partial_2^2A(Z_\tjm,\theta_2)\\
2^{-1}H_x(Z_\tjm,\hat{\theta}_3^J)[\partial_2^2A(Z_\tjm,\theta_2)]
\end{array}\right)
\right]
\\&=&
-n^{-1}\sum_{j=1}^nS(Z_\tjm,\hat{\theta}_1^J,\hat{\theta}_3^J)^{-1}
\left[
\left(\begin{array}{c} \partial_2A(Z_\tjm,\theta_2)\\
2^{-1}H_x(Z_\tjm,\hat{\theta}_3^J)[\partial_2A(Z_\tjm,\theta_2)]
\end{array}\right)^{\otimes2}
\right]
\\&&
+n^{-1}h^{-1/2}\sum_{j=1}^nS(Z_\tjm,\theta_1^*,\theta_3^*)^{-1}
\left[\cald_j(\theta_1^*,\theta_2,\theta_3^*),\>
\left(\begin{array}{c} \partial_2^2A(Z_\tjm,\theta_2)\\
2^{-1}H_x(Z_\tjm,\theta_3^*)[\partial_2^2A(Z_\tjm,\theta_2)]
\end{array}\right)
\right]
\\&&+o_p(n^{-1/2}h^{-1/2})
\quad(\text{Lemmas \ref{201906091323} and \ref{201906091405}})
\\&=&
-n^{-1}\sum_{j=1}^nC(Z_\tjm,\hat{\theta}_1)^{-1}\big[\big(\partial_2A(Z_\tjm,\theta_2)\big)^{\otimes2}\big]
\\&&
+n^{-1}h^{-1/2}\sum_{j=1}^nC(Z_\tjm,\theta_1^*)^{-1}\big[\partial_2^2A(Z_\tjm,\theta_2),\>
h^{-1/2}(\Delta_jX-hA(Z_\tjm,\theta_2^*))\big]+o_p(1)
\\&=&
-n^{-1}\sum_{j=1}^nC(Z_\tjm,\theta_1^*)^{-1}\big[\big(\partial_2A(Z_\tjm,\theta_2)\big)^{\otimes2}\big]
\quad(\text{Lemma \ref{201906091405}})
\\&&
+n^{-1}h^{-1/2}\sum_{j=1}^nC(Z_\tjm,\theta_1^*)^{-1}\big[\partial_2^{{\cred 2}}A(Z_\tjm,\theta_2),\>
h^{-1/2}(\Delta_jX-hA(Z_\tjm,\theta_2^*))\big]
+o_p(1)
\\&=&
-n^{-1}\sum_{j=1}^nC(Z_\tjm,\theta_1^*)^{-1}\big[\big(\partial_2A(Z_\tjm,\theta_2)\big)^{\otimes2}\big]
+o_p(1)
\eeas
{\coloro The order $o_p(1)$ is uniform in $\theta_2\in\Theta_2$.}  
The last equation is verified by random field argument with the shrinking $B(\theta_2^*,r_n)$, 
{\cred where we need $\partial_2^3A$. }
{\coloro Since $\hat{\theta}_2^J\to^p\theta_2^*$ by Lemma \ref{201906080713}(b), 
applying Taylor's formula with $\partial_2^2A$, we obtain (\ref{201911180914}) 
with the help of Lemma \ref{201906041850} (a). 
Moreover, we obtain (\ref{201911180915}) by combining (\ref{201911180914}) with Lemma \ref{201906081515}. 
}
\qed\halflineskip

Let 
\beas 
\hat{\theta}^J &=& 
\left(\begin{array}{c}\hat{\theta}_1^J \\ \hat{\theta}_2^J \\ \hat{\theta}_3^J\end{array}\right).
\eeas
Recall 
\beas 
\Gamma^J(\theta^*)
&=&
\text{diag}\bigg[
\half\int \bigg[
 \rm{Tr}\big\{\big(C^{-1}(\partial_1C)C^{-1}\partial_1C\big)(z,\theta_1^*)\big\}
 \\&&\qquad\quad
 +\rm{Tr}\big\{\big(V^{-1}H_x(\partial_1C)H_x^\star V^{-1}H_x(\partial_1C)H_x^\star\big)(z,\theta_1^*,\theta_3^*)\big\}
 \bigg]\nu(dz),
 \\&&\quad\quad
 \int \partial_2A(z,\theta_2^*)^\star C(z,\theta_1^*)^{-1}\partial_2A(z,\theta_2^*)\nu(dz), 
\\&&\quad\quad
\int12\partial_3H(z,\theta_3^*)^\star V(z,\theta_1^*,\theta_3^*)^{-1}\partial_3H(z,\theta_3^*)\nu(dz)
\bigg]
\eeas

\begin{theorem} Suppose that $[A1]$ with $(i_A,j_A,i_B,j_B,i_H,j_H)=
(1,{\cred 3},2,3,3,2)
$, 
$[A2]$, {\coloro $[A3](i')$, $[A3](ii)$ and $[A3](iii')$} 
are satisfied. 
Then 
\beas 
b_n^{-1} \big(\hat{\theta}^J-\theta^*\big)
&\to^d& 
N\big(0,\big(\Gamma^J(\theta^*)\big)^{-1}\big)
\eeas
as $n\to\infty$. 
\end{theorem}
\proof 
By Lemmas \ref{201906091405}, \ref{201906091659} and \ref{201906091542}, we obtain the result. 
By simple linear calculus, we can see 
that $M_n^{(2)}$ and $M_n^{(3)}$ are asymptotically orthogonal. 
Since $M_n^{(1)}$ is written by the second Wiener chaos, it is asymptotically orthogonal 
to $M_n^{(2)}$ and $M_n^{(3)}$. 
\qed\halflineskip
}


\section{Estimation of $\theta_1$}\label{201905291607}
The purpose of this section is to recall a standard construction of 
estimator for $\theta_1$ and to clarify what conditions we mentioned 
validate its asymptotic properties. 
Let 
\beas 
\bbH^{(1)}_n(\theta_1) 
&=& 
-\half \sum_{j=1}^n \bigg\{C(Z_\tjm,\theta_1)^{-1}\big[h^{-1}(\Delta_jX)^{\otimes2}\big]
+\log\det C(Z_\tjm,\theta_1)\bigg\}
\eeas
where $\Delta_jX=X_\tj-X_\tjm$. 
{\cred It should be remarked that the present $\bbH^{(1)}_n(\theta_1) $ is different from 
the one given in (\ref{201906060545}) on p.\pageref{201906060545}. 
}
Under $[A1]$ and $[A2]$ (iii), $\bbH^{(1)}_n$ is a continuous function on $\overline{\Theta}_1$ a.s. 

Given the data $(Z_\tj)_{j=0,1,...,n}$, 
let us consider the quasi-maximum likelihood estimator (QMLE) $\hat{\theta}_1^0=\hat{\theta}_{1,n}^0$  
for $\theta_1$, that is, 
$\hat{\theta}_1^0$ is any measurable function of $(Z_\tj)_{j=0,1,...,n}$ satisfying 
\beas 
\bbH^{(1)}_n(\hat{\theta}_1^0) &=& \max_{\theta_1\in\overline{\Theta}_1}\bbH^{(1)}_n(\theta_1)\quad a.s.
\eeas

Routinely, $n^{1/2}$-consistency and asymptotic normality of $\hat{\theta}_1^0$ can be established. 
We will give a brief for self-containedness and for the later use. 
Let 
\bea\label{201905281059}
\Gamma^{(1)}[u_1^{\otimes2}]
&=&
\half\int_{\bbR^{\sfd_Z}}\text{Tr}\big\{C^{-1}(\partial_1C)[u_1]C^{-1}(\partial_1C)[u_1](z,\theta_1^*)\big\}\nu(dz)
\eea
for $u_1\in\bbR^{\sfp_1}$. 
We will see the existence and positivity of $\Gamma^{(1)}$ in the following theorem. 
\begin{theorem}\label{201905291601}
{\bf (a)} Suppose that $[A1]$ with $(i_A,j_A,i_B,j_B,i_H,j_H)=(0,0,1,1,0,0)$, $[A2]$ $(i), (ii), (iii)$, and 
$[A3]$ $(i)$ are satisfied. 
Then $\hat{\theta}_1^0\to^p\theta_1^*$ as $n\to\infty$. 
\bd\im[(b)] Suppose that $[A1]$ with $(i_A,j_A,i_B,j_B,i_H,j_H)=(1,0,2,3,0,0)$, $[A2]$ $(i), (ii), (iii)$, and 
$[A3]$ $(i)$ are satisfied. Then $\Gamma^{(1)}$ exists and is positive-definite, and 
\beas
\sqrt{n}\big(\hat{\theta}_1^0-\theta_1^*\big)
-(\Gamma^{(1)})^{-1}{\colorb \hat{M}}^{(1)}_n&\to^p&0
\eeas
as $n\to\infty$, where 
\beas 
{\colorb \hat{M}}^{(1)}_n 
&=&
\half n^{-1/2}\sum_{j=1}^n \big(C^{-1}(\partial_1C)C^{-1}\big)(Z_\tjm,\theta_1^*)
\big[\big(h^{-1/2}B(Z_\tjm,\theta_1^*)\Delta_jw\big)^{\otimes2}
-C(Z_\tjm,\theta_1^*)\big].
\eeas
Moreover, 
$
M^{(1)}_n
\to^d 
N_{\sfp_1}(0,\Gamma^{(1)})
$
as $n\to\infty$. In particular, 
\beas 
\sqrt{n}\big(\hat{\theta}_1^0-\theta_1^*\big)
&\to^d&
N_{\sfp_1}\big(0,(\Gamma^{(1)})^{-1}\big)
\eeas
as $n\to\infty$.
\ed
\end{theorem}
\proof (a): 
Let $\bbY^{(1)}_n(\theta_1)=n^{-1}\big(\bbH_n^{(1)}(\theta_1)-\bbH_n^{(1)}(\theta_1^*)\big)$. 
Suppose that $[A1]$ with $(i_A,j_A,i_B,j_B,i_H,j_H)=(0,0,1,1,0,0)$ and $[A2]$ (i), (ii), (iii). 
Use (\ref{201906021458}) and Lemma \ref{201906041850}, 
then 
\beas 
\sum_{i=0}^1\sup_{\theta_1\in\Theta_1}\bigg\|\big|\partial_1^i\big(\bbY^{(1)}_n(\theta_1)-\bbY^{(1)}(\theta_1)\big)
\big|\bigg\|_p
&\to& 0\quad(n\to\infty)
\eeas
for every $p>1$. By Sobolev's inequality, we obtain 
\beas
\bigg\|\sup_{\theta_1\in\Theta_1}|\bbY^{(1)}_n(\theta_1)-\bbY^{(1)}(\theta_1)|\bigg\|_p&\to&0
\eeas
for every $p>1$. Therefore, the identifiability condition $[A3]$ (i) ensures $\hat{\theta}_1^0\to^p\theta_1^*$ 
as $n\to\infty$. 
\halflineskip

(b): Under $[A1]$ with $(i_A,j_A,i_B,j_B,i_H,j_H)=(0,0,0,2,0,0)$, we have
\beas 
\partial_1\bbH^{(1)}_n(\theta_1) 
&=& 
\half \sum_{j=1}^n \big(C^{-1}(\partial_1C)C^{-1}\big)(Z_\tjm,\theta_1)\big[h^{-1}(\Delta_jX)^{\otimes2}
-C(Z_\tjm,\theta_1)\big]
\eeas
and 
\beas 
\partial_1^2\bbH^{(1)}_n(\theta_1) 
&=& 
\half \sum_{j=1}^n \partial_1\big(C^{-1}(\partial_1C)C^{-1}\big)(Z_\tjm,\theta_1)\big[h^{-1}(\Delta_jX)^{\otimes2}
-C(Z_\tjm,\theta_1)\big]
\\&&
-\half \sum_{j=1}^n \big(C^{-1}(\partial_1C)C^{-1}\big)(Z_\tjm,\theta_1)\big[\partial_1C(Z_\tjm,\theta_1)\big].
\eeas\

Now, by using orthogonality and the estimate (\ref{201905280613}), 
if $[A1]$ is satisfied for $(i_A,j_A,i_B,j_B,i_H,j_H)=(1,0,2,1,0,0)$, then 
\beas &&
n^{-1/2}\partial_1\bbH^{(1)}_n(\theta_1^*) 
\\&=& 
\half n^{-1/2} \sum_{j=1}^n \big(C^{-1}(\partial_1C)C^{-1}\big)(Z_\tjm,\theta_1^*)\big[h^{-1}(\Delta_jX)^{\otimes2}
-C(Z_\tjm,\theta_1^*)\big]
\\&=&
M^{(\ref{201905281048})}_n
+O_{L^\inftym}(n^{-1/2})+O_{L^\inftym}(n^{1/2}h)
\eeas
where 
\bea\label{201905281048}
M^{(\ref{201905281048})}_n
&=&
\half n^{-1/2}\sum_{j=1}^n \big(C^{-1}(\partial_1C)C^{-1}\big)(Z_\tjm,\theta_1^*)
\bigg[\bigg(h^{-1/2}\int_\tjm^\tj B(Z_t,\theta_1^*)dw_t\bigg)^{\otimes2}
-C(Z_\tjm,\theta_1^*)\bigg].
\nn\\&&
\eea
At the same time It\^o's formula gives
\bea\label{201905281045}
h^{-1/2}\int_\tjm^\tj B(Z_t,\theta_1^*)dw_t
&=&
h^{-1/2}B(Z_\tjm,\theta_1^*)\Delta_jw
+h^{-1/2}\int_\tjm^\tj \big(B(Z_t,\theta_1^*)-B(Z_\tjm,\theta_1^*)\big)dw_t
\nn\\&=&
h^{-1/2}B(Z_\tjm,\theta_1^*)\Delta_jw
+h^{-1/2}\int_\tjm^\tj \int_\tjm^tB_x(Z_s,\theta_1^*)[B(Z_s,\theta_1^*)dw_s]dw_t
\nn\\&&
+h^{-1/2}\int_\tjm^\tj \int_\tjm^tL_B(Z_s,\theta_1^*,\theta_2^*,\theta_3^*)dsdw_t
\eea
for $L_B(z,\theta_1,\theta_2,\theta_3)$ given by (\ref{201996160909}). 
The products of the first two terms on the right-hand side of (\ref{201905281045}) form martingale differences, 
and hence 
\beas 
M^{(\ref{201905281048})}_n
&=&
M^{(1)}_n
+O_{L^\inftym}(h^{1/2})+O_{L^\inftym}(n^{1/2}h).
\eeas
Under $[A2]$ (i), (ii), (iii), the martingale central limit theorem gives 
\beas
M^{(1)}_n
&\to^d& 
N_{\sfp_1}(0,\Gamma^{(1)})
\eeas
as $n\to\infty$. 
Consequently, 
\bea\label{201905281219}
n^{-1/2}\partial_1\bbH^{(1)}_n(\theta_1^*) &\to^d& N(0,\Gamma^{(1)})
\eea
if $[A1]$ with $(i_A,j_A,i_B,j_B,i_H,j_H)=(1,0,2,1,0,0)$ and $[A2]$ (i), (ii), (iii) are fulfilled.

Next, suppose that  $[A1]$ with $(i_A,j_A,i_B,j_B,i_H,j_H)=(0,0,2,3,0,0)$ and $[A2]$ (i), (ii), (iii) are fulfilled. 
It is rather simple to prove 
\bea\label{201905281202} 
\sum_{i=0}^1\sup_{u_1:\theta_1^*+\rho_nu_1\in\Theta_1, |u_1|<1}
\bigg\|\partial_{u_1}^i\bigg( n^{-1}\partial_1^2\bbH^{(1)}_n(\theta_1^*+\rho_nu_1) 
+\Gamma^{(1)}\bigg)\bigg\|_p
&\to&0
\eea
for every $p>1$ and any sequence $(\rho_n)_{n\in\bbN}$ of positive numbers tending to $0$ as $n\to\infty$. 
We apply Sobolev's embedding inequality 
to each component of the matrix valued random field 
$\partial_1^2\bbH^{(1)}_n(\theta_1^*+\rho_nu_1)$ on $\{u_1\in\bbR^{\sfp_1};\>|u_1|<1\}$ for large $n$. 
Then (\ref{201905281202}) gives 
\bea\label{201905281210} 
\bigg\|\sup_{u_1:\theta_1^*+\rho_nu_1\in\Theta_1, |u_1|<1}
\big|n^{-1}\partial_1^2\bbH^{(1)}_n(\theta_1^*+\rho_nu_1) 
+\Gamma^{(1)}\big|\bigg\|_p
&\to&0
\eea
for every $p>1$. 

Suppose that $[A1]$ with $(i_A,j_A,i_B,j_B,i_H,j_H)=(1,0,2,3,0,0)$, $[A2]$ $(i), (ii), (iii)$, and 
$[A3]$ $(i)$ are satisfied. Then differentiating $\bbY^{(1)}$ twice, we see, from $[A3]$ (i), 
that $\Gamma^{(1)}$ is positive-definite. 
By (a), $\hat{\theta}_1^0\to^p\theta_1^*$. 
With this fact, we obtain (b) from (\ref{201905281219}) and (\ref{201905281210}). 
\qed\halflineskip

\begin{remark}\rm
It is possible to show that the quasi-Bayesian estimator (QBE) also enjoys 
the same asymptotic properties as the QMLE in Theorem \ref{201905291601}, 
if we follows the argument in Yoshida \cite{yoshida2011polynomial}. 
This means we can use both estimators together with the estimator for $\theta_2$ 
e.g. given in Section \ref{201906041938}, 
to construct a one-step estimator for $\theta_3$ 
based on the scheme presented in Section \ref{201906041935}, 
and consequently we can construct a one-step estimator for 
$\theta=(\theta_1,\theta_2,\theta_3)$ by the method in Section \ref{202001141623}. 
\end{remark}

\section{Estimation of $\theta_2$}\label{201906041938}
This section will recall a standard construction of estimator for $\theta_2$. 
As usual, the scheme is adaptive. 
Suppose that an estimator $\hat{\theta}_1^0$ based on the data $(Z_\tj)_{j=0,1,...,n}$ 
satisfies Condition $[A4^\sharp]$ (i), i.e., 
\beas 
\hat{\theta}_1^0 -\theta_1^* &=& O_p(n^{-1/2})
\eeas
as $n\to\infty$. 
Obviously we can apply the estimator constructed in Section \ref{201905291607}, but 
any estimator satisfying this condition can be used. 

Define the random field $\bbH^{(2)}_n$ on $\overline{\Theta}_2$ by 
\bea\label{202002171935} 
\bbH^{(2)}_n(\theta_2) 
&=&
-\half\sum_{j=1}^n C(Z_\tjm,\hat{\theta}_1^0)^{{\cred -1}}
\big[h^{-1}\big(\Delta_jX-hA(Z_\tjm,\theta_2)\big)^{\otimes2}\big].
\eea
We will denote by $\hat{\theta}_2^0=\hat{\theta}_{2,n}^0$ any sequence of 
quasi-maximum likelihood estimator for $\bbH^{(2)}_n$, that is, 
\beas 
\bbH^{(2)}_n(\hat{\theta}_2^0) &=& \sup_{\theta_2\in\overline{\Theta}_2}\bbH^{(2)}_n(\theta_2).
\eeas

Let $\bbY^{(2)}_n(\theta_2)=T^{-1}\big(\bbH^{(2)}_n(\theta_2)-\bbH^{(2)}_n(\theta_2^*)\big)$, 
{\cred where $T=nh$.}

\begin{lemma}\label{202002171923}
Suppose that Conditions $[A1]$ with $(i_A,j_A,i_B,j_B,i_H,j_H)=(1,1,2,1,0,0)$, $[A2]$ and 
$[A4^\sharp]$ $(i)$. 
Then 
\beas
\sup_{\theta_2\in\overline{\Theta}_2}
\big|\bbH^{(2)}_n(\theta_2) -\bbY^{(2)}_n(\theta_2)\big|
&\to^p&
0
\eeas
as $n\to\infty$. 
If additionally $[A3]$ $(ii)$ is satisfied, then $\hat{\theta}_2^0\to^p\theta_2^*$ 
as $n\to\infty$. 
\end{lemma}
\proof 
\beas 
\bbY^{(2)}_n(\theta_2)
&=&
\begin{en-text}
-\frac{1}{2T}\sum_{j=1}^n h C(Z_\tjm,\hat{\theta}_1^0)^{{\cred -1}}
\big[\big(A(Z_\tjm,\theta_2)-A(Z_\tjm,\theta_2^*)\big)^{\otimes2}\big]
\\&&
+
T^{-1}\sum_{j=1}^n  C(Z_\tjm,\hat{\theta}_1^0)^{{\cred -1}}
\big[\Delta_jX-hA(Z_\tjm,\theta_2^*),A(Z_\tjm,\theta_2)-A(Z_\tjm,\theta_2^*)\big]
\end{en-text}
\Phi^{(\ref{202002171801})}_n(\theta_2)
+\Phi^{(\ref{202002171802})}_n(\hat{\theta}_1^0,\theta_2)
\eeas
where 
\bea\label{202002171801}
\Phi^{(\ref{202002171801})}_n(\theta_2)
&=&
-\frac{1}{2T}\sum_{j=1}^n h C(Z_\tjm,\hat{\theta}_1^0)^{{\cred -1}}
\big[\big(A(Z_\tjm,\theta_2)-A(Z_\tjm,\theta_2^*)\big)^{\otimes2}\big]
\eea
and 
\bea\label{202002171802}
\Phi^{(\ref{202002171802})}_n(\theta_1,\theta_2)
&=&
T^{-1}\sum_{j=1}^n  C(Z_\tjm,\theta_1)^{{\cred -1}}
\big[\Delta_jX-hA(Z_\tjm,\theta_2^*),A(Z_\tjm,\theta_2)-A(Z_\tjm,\theta_2^*)\big].
\nn\\&&
\eea

For $\bbY^{(2)}$ given by (\ref{202002171821}) on p.\pageref{202002171821}, 
\begin{en-text}
\bea\label{202002171824}
\Phi^{(\ref{202002171801})}_n(\theta_2)
&\to^p&
\bbY^{(2)}(\theta_2)
\eea
for every $\theta_2\in\overline{\Theta}_2$, 
and 
\end{en-text}
\bea\label{202002171825}
\sum_{i=0}^1\sup_{\theta_2\in\Theta_2}
\big\|\partial_2^i\big(\Phi^{(\ref{202002171801})}_n(\theta_2)-\bbY^{(2)}(\theta_2)
\big)\big\|_p
&\to& 0
\eea
for every $p>1$. 
Here 
Conditions 
$[A1]$ (i) with $(i_A,j_A,i_B,j_B,i_H,j_H)=(1,1,1,1,0,0)$, $[A2]$ (i)-(iii) and 
$[A4^\sharp]$ (i) were used. 
Then 
(\ref{202002171825}) implies
\bea\label{202002171852}
\sup_{\theta_2\in\overline{\Theta}_2}\big|
\Phi^{(\ref{202002171801})}_n(\theta_2)-\bbY^{(2)}(\theta_2)\big|
&\to^p&
0
\eea
as $n\to\infty$. 

We have 
\bea\label{202002171900}
\sum_{i=0}^1\sup_{(\theta_1,\theta_2)\in\Theta_1\times\Theta_2}
\big\|\partial_{(\theta_1,\theta_2)}^i\Phi^{(\ref{202002171802})}_n(\theta_1,\theta_2)\big\|_p
&\to^p&
0
\eea
for every $p>1$ 
from Lemma \ref{201906030207} applied to $\Delta_jX-hA(Z_\tjm,\theta_2^*)$ 
with the aid of orthogonality. 
The conditions we used include 
$[A1]$ with $(i_A,j_A,i_B,j_B,i_H,j_H)=(1,1,2,1,0,0)$. 
The embedding inequality makes 
\bea\label{202002171901}
\sup_{(\theta_1,\theta_2)\in\overline{\Theta}_1\times\overline{\Theta}_2}
\big|\Phi^{(\ref{202002171802})}_n(\theta_1,\theta_2)\big|
&\to^p&
0
\eea
from (\ref{202002171900}). 

The proof completes by the estimates (\ref{202002171852}) and (\ref{202002171901}). 
\qed\halflineskip

The matrix $\Gamma_{22}$ is defined by (\ref{202001281931}) on p.\pageref{202001281931}. 

\begin{lemma}\label{202002172006}
Under Conditions $[A1]$ with $(i_A,j_A,i_B,j_B,i_H,j_H)=(1,3,2,1,0,0)$, $[A2]$ and $[A4^\sharp]$, 
Then 
\bea\label{202002171957}
\sup_{\theta_2\in U(\theta_2^*,r_n)}\big|
T^{-1}\partial_2^2\bbH^{(2)}_n(\theta_2) -\Gamma_{22}\big|
&\to^p&
0
\eea
for any sequence of positive numbers $r_n$ satisfying $r_n\to0$ as $n\to\infty$. 
\end{lemma}
\proof
From (\ref{202002171935}), 
\bea\label{202002171936} 
\partial_2\bbH^{(2)}_n(\theta_2) 
&=&
\sum_{j=1}^n C(Z_\tjm,\hat{\theta}_1^0)^{{\cred -1}}
\big[\Delta_jX-hA(Z_\tjm,\theta_2),\partial_2A(Z_\tjm,\theta_2)\big]
\eea
and 
\bea\label{202002171937} 
\partial_2^2\bbH^{(2)}_n(\theta_2) 
&=&
-\sum_{j=1}^n C(Z_\tjm,\hat{\theta}_1^0)^{{\cred -1}}
\big[\big(\partial_2A(Z_\tjm,\theta_2)\big)^{\otimes2}\big]h
\nn\\&&
+\sum_{j=1}^n C(Z_\tjm,\hat{\theta}_1^0)^{{\cred -1}}
\big[\Delta_jX-hA(Z_\tjm,\theta_2),\partial_2^2A(Z_\tjm,\theta_2)\big]
\eea

Let 
\bea\label{202002171941} 
\hat{M}^{(2)}_n
&=&
T^{-1/2}\sum_{j=1}^n C(Z_\tjm,\theta_1^*)^{-1}
\big[B(Z_\tjm,\theta_1^*)\Delta_jw,\partial_2A(Z_\tjm,\theta_2^*)\big]
\eea

By random field argument for (\ref{202002171936}) 
with $\partial_1C(Z_\tjm,\theta_1)$ and $[A4^\sharp]$, 
we obtain 
\bea\label{202002171950} 
T^{-1/2}\partial_2\bbH^{(2)}_n(\theta_2^*) 
-\hat{M}^{(2)}_n
&\to^p&
0
\eea
as $n\to\infty$, if (\ref{201906030207}) in Lemma \ref{201906031518} 
applied to $\Delta_jX-hA(Z_\tjm,\theta_2^*)$.

Under Conditions $[A1]$ with $(i_A,j_A,i_B,j_B,i_H,j_H)=(1,3,2,1,0,0)$, $[A2]$ and $[A4^\sharp]$, 
we obtain the convergence (\ref{202002171957}) 
for any sequence of positive numbers $r_n$ satisfying $r_n\to0$ as $n\to\infty$. 
Here we applied random field argument to the second term on the right-hand side of 
(\ref{202002171937}). 
\qed\halflineskip

The matrix $\Gamma_{22}$ is given by (\ref{202001281931}) on p.\pageref{202001281931}. 
Form Lemmas \ref{202002171923}, \ref{202002171923} and \ref{202002172006}, 
the following theorem follows. 

\begin{theorem}\label{202002172010}
{\bf (a)} 
Suppose that Conditions $[A1]$ with $(i_A,j_A,i_B,j_B,i_H,j_H)=(1,1,2,1,0,0)$, $[A2]$, 
$[A3]$ $(ii)$ and 
$[A4^\sharp]$ $(i)$. 
Then 
$\hat{\theta}_2^0\to^p\theta_2^*$ as $n\to\infty$. 
\bd
\im[(b)]
Suppose that Conditions $[A1]$ with $(i_A,j_A,i_B,j_B,i_H,j_H)=(1,3,2,1,0,0)$, $[A2]$, 
$[A3]$ $(ii)$ and 
$[A4^\sharp]$ $(i)$. 
Then 
\beas
{\cred(nh)}^{1/2}\big(\hat{\theta}_2^0-\theta_2^*\big)
-\Gamma_{22}^{-1}\hat{M}^{(2)}_n
&\to^p&
0
\eeas
as $n\to\infty$. 
In particular, 
\beas 
{\cred(nh)}^{1/2}\big(\hat{\theta}_2^0-\theta_2^*\big)
&\to^d&
N(0,\Gamma_{22}^{-1})
\eeas
as $n\to\infty$. 
\ed
\end{theorem}
}

{\cred 
\begin{remark}\rm 
The estimator $\hat{\theta}_1^0$ in Section\ref{201905291607} is 
asymptotically orthogonal to $\hat{\theta}_2^0$ constructed in this section. 
Therefore, for that $\hat{\theta}_1^0$, we obtain the joint convergence 
\beas 
\big(n^{1/2}(\hat{\theta}_1^0-\theta_1^*),(nh)^{1/2}(\hat{\theta}_2^0-\theta_2^*)\big)
&\to^d& 
N\big(0,\text{diag}[(\Gamma^{(1)})^{-1},\Gamma_{22}^{-1}]\big)
\eeas
as $n\to\infty$, 
as is well known. 
\end{remark}
}

\bibliographystyle{spmpsci}      
\bibliography{bibtex-20180426-20180615bis}   

\color{black}
\newpage
\section{Symbols and Conditions} 
\subsection{Section \ref{202002201325}}
\beas
\left\{\begin{array}{ccl}
dX_t &=& A(Z_t,\theta_2)dt+B(Z_t,\theta_1)dw_t
\y
dY_t&=&H(Z_t,\theta_3)dt
\end{array}\right.
\eeas

\subsection{Section \ref{202001141611}}

\bd
\im[[A1\!\!]]
{\bf (i)} 
$A\in C^{i_A,j_A}_p(\bbR^{\sfd_Z}\times\Theta_2;\bbR^{\sfd_X})$ and 
$B\in C^{i_B,j_B}_p(\bbR^{\sfd_Z}\times\Theta_1;\bbR^{\sfd_X}\otimes\bbR^\sfr)$. 

\bd\im[(ii)] $H\in C^{i_H,j_H}_p(\bbR^{\sfd_Z}\times\Theta_3;\bbR^{\sfd_Y})$. 
\ed
\ed

\bd
\im[[A2\!\!]] {\bf (i)} 
$\sup_{t\in\bbR_+}\|Z_t\|_p<\infty$ 
for every $p>1$. 
\bd\im[(ii)] There exists a probability measure $\nu$ on $\bbR^{\sfd_Z}$ such that 
\beas 
\frac{1}{T}\int_0^T f(Z_t)\>dt &\to^p& \int f(z)\nu(dz)\qquad(T\to\infty)
\eeas
for any bounded measurable function $f:\bbR^{\sfd_Z}\to\bbR$. 

\im[(iii)] {\coloro The function $\theta_1\mapsto C(Z_t,\theta_1)^{-1}$ is continuous 
{\colorg on $\overline{\Theta}_1$}
a.s., and} 
\beas
\sup_{\theta_1\in\Theta_1}\sup_{t\in\bbR_+}\|\det C(Z_t,\theta_1)^{-1}\|_p&<&\infty
\eeas 
for every $p>1$. 

\im[(iv)] For the $\bbR^{{\colorr\sfd_Y}}\otimes\bbR^{{\colorr\sfd_Y}}$ valued function 
$V(z,\theta_1,\theta_3)=H_x(z,\theta_3)C(z,\theta_1)H_x(z,\theta_3)^\star$, 
{\coloro the function $(\theta_1,\theta_3)\mapsto V(Z_t,\theta_1,\theta_3)^{-1}$ is continuous 
{\colorg on $\overline{\Theta}_1\times\overline{\Theta}_3$} a.s., 
and }
\beas 
\sup_{(\theta_1,\theta_3)\in\Theta_1\times\Theta_3}\sup_{t\in\bbR_+}\|\det V(Z_t,\theta_1,\theta_3)^{-1}\|_p&<&\infty
\eeas 
for every $p>1$. 
\ed
\ed

\beas 
\bbY^{(1)}(\theta_1) 
&=& 
-\half\int \bigg\{
\text{Tr}\big(C(z,\theta_1)^{-1}C(z,\theta_1^*)\big)-\sfd_X
+\log\frac{\det C(z,\theta_1)}{\det C(z,\theta_1^*)}\bigg\}\nu(dz).
\eeas
Since 
$|\log x|\leq x+x^{-1}$ for $x>0$, $\bbY^{(1)}(\theta_1)$ is a continuous function on $\overline{\Theta}_1$ 
well defined 
under {\colorg $[A1]$ and }$[A2]$. 
Let 
\beas 
\bbY^{(J,1)}(\theta_1)
&=& 
-\half\int \bigg\{\text{Tr}\big(C(z,\theta_1)^{-1}C(z,\theta_1^*)\big)
+\text{Tr}\big(V(z,\theta_1,\theta_3^*)^{-1}V(z,\theta_1^*,\theta_3^*)\big)
-\sfd_Z
\\&&
+\log \frac{\det C(z,\theta_1)\det V(\theta_1,\theta_3^*)}{\det C(z,\theta_1^*)\det V(\theta_1^*,\theta_3^*)}
\bigg\}
\nu(dz)
\eeas
Let 
\beas 
\bbY^{(2)}(\theta_2)
&=&
-\half \int 
C(z,\theta_1^*)\big[\big(A(z,\theta_2)-A(z,\theta_2^*)\big)^{\otimes2}\big]\nu(dz)
\eeas
Let 
\beas 
\bbY^{(3)}(\theta_3) &=&
-\int 6V(z,\theta_1^*,\theta_3)^{-1}
\big[\big({\colorr H(z,\theta_3)-H(z,\theta_3^*)}\big)^{\otimes2}\big]\nu(dz). 
\eeas
The random field $\bbY^{(3)}$ is well defined under {\colorg $[A1]$ and }$[A2]$. 
Let 
\beas 
\bbY^{(J,3)}(\theta_1,\theta_3) &=&
-\int 6V(z,\theta_1,\theta_3)^{-1}
\big[\big(H(z,\theta_3)-H(z,\theta_3^*)\big)^{\otimes2}\big]\nu(dz). 
\eeas

We will assume all or some of the following identifiability conditions
\bd
\im[[A3\!\!]] 
{\bf (i)} There exists a positive constant $\chi_1$ such that 
\beas 
\bbY^{(1)}(\theta_1) &\leq& -\chi_1|\theta_1-\theta_1^*|^2\qquad(\theta_1\in\Theta_1). 
\eeas
\bd 
\im[(i$'$)] There exists a positive constant $\chi_1'$ such that 
\beas 
\bbY^{(J,1)}(\theta_1)&\leq& -\chi_1'|\theta_1-\theta_1^*|^2\qquad(\theta_1\in\Theta_1). 
\eeas

\im[(ii)] There exists a positive constant $\chi_2$ such that 
\beas 
\bbY^{(2)}(\theta_1) &\leq& -\chi_2|\theta_2-\theta_2^*|^2\qquad(\theta_2\in\Theta_2). 
\eeas

\im[(iii)]
There exists a positive constant $\chi_3$ such that 
\beas 
\bbY^{(3)}(\theta_3) &\leq& -\chi_3|\theta_3-\theta_3^*|^2\qquad(\theta_3\in\Theta_3). 
\eeas

\im[(iii$'$)]
There exists a positive constant $\chi_3$ such that 
\beas 
\bbY^{(J,3)}(\theta_1,\theta_3) &\leq& -\chi_3|\theta_3-\theta_3^*|^2
\qquad(\theta_1\in\Theta_1,\>\theta_3\in\Theta_3). 
\eeas
\ed
\ed

\subsection{Section \ref{202001141618}}
\beas
L_H(z,\theta_1,\theta_2,\theta_3)
&=&
H_x(z,\theta_3)[A(z,\theta_2)]+\half H_{xx}(z,\theta_3)[C(z,\theta_1)]+H_y(z,\theta_3)[H(z,\theta_3)].
\eeas
\beas 
G_n(z,\theta_1,\theta_2,\theta_3) 
&=& 
H(z,\theta_3)+\frac{h}{2} L_H\big(z,\theta_1,\theta_2,\theta_3\big).
\eeas
\beas 
\zeta_j &=& \sqrt{{\colorr{3}}}\int_\tjm^\tj\int_\tjm^t dw_sdt
\eeas
\beas
\cald_j(\theta_1,\theta_2,\theta_3) &=&
\left(\begin{array}{c}
h^{-1/2}\big(\Delta_jX-{\colorr h}A(Z_\tjm,\theta_2)\big)\y
h^{-3/2}\big(\Delta_jY-{\colorr h}G_n(Z_\tjm,\theta_1,\theta_2,\theta_3)\big)
\end{array}\right).
\eeas

\subsection{Section \ref{201906041935}}
\beas
S(z,\theta_1,\theta_3)&=&
\left(\begin{array}{ccc}
C(z,\theta_1)&& 2^{-1}C(z,\theta_1)H_x(z,\theta_3)^\star
\\
2^{-1}H_x(z,\theta_3)C(z,\theta_1)&& 3^{-1}H_x(z,\theta_3)C(z,\theta_1)H_x(z,\theta_3)^\star
\end{array}\right)
\eeas

\beas &&
S(z,\theta_1,\theta_3)^{-1}
\nn\\&=&
\left(\begin{array}{ccc}
C(z,\theta_1)^{-1}+3H_x(z,\theta_3)^\star V(z,\theta_1,\theta_3)^{-1}H_x(z,\theta_3)
&&-6H_x(z,\theta_3)^\star V(z,\theta_1,\theta_3)^{-1}\y
-6V(z,\theta_1,\theta_3)^{-1}H_x(z,\theta_3)&&12V(z,\theta_1,\theta_3)^{-1}
\end{array}\right)
\nn\\&&
\eeas

\beas
\hat{S}(z,\theta_3)&=&S(z,\hat{\theta}_1^0,\theta_3)
\eeas
\beas
\bbH^{(3)}_n(\theta_3)
 &=& 
 -\half\sum_{j=1}^n \bigg\{
 \hat{S}(Z_\tjm,\theta_3)^{-1}\big[\cald_j(\hat{\theta}_1^0,\hat{\theta}_2^0,\theta_3)^{\otimes2}\big]
 +\log\det \hat{S}(Z_\tjm,\theta_3)\bigg\}
\eeas

\beas
\Psi_1(\theta_1,\theta_2,\theta_3,\theta_1',\theta_2',\theta_3')
&=&
\sum_{j=1}^nS(Z_\tjm,\theta_1,\theta_3)^{-1}
\left[\cald_j(\theta_1',\theta_2',\theta_3'),\>
\left(\begin{array}{c}0\\
2^{-1}\partial_1 L_H(Z_\tjm,\theta_1,\theta_2,\theta_3)
\end{array}\right)\right]
\\&=&
\sum_{j=1}^nS(Z_\tjm,\theta_1,\theta_3)^{-1}
\left[\cald_j(\theta_1',\theta_2',\theta_3'),\>
\left(\begin{array}{c}0\\
2^{-1}H_{xx}(z,\theta_3)[\partial_1C(Z_\tjm,\theta_1)]
\end{array}\right)\right]
\eeas
\beas
\Psi_2(\theta_1,\theta_2,\theta_3,\theta_1',\theta_2',\theta_3')
&=&
\sum_{j=1}^nS(Z_\tjm,\theta_1,\theta_3)^{-1}
\left[\cald_j(\theta_1',\theta_2',\theta_3'),\>
\left(\begin{array}{c}\partial_2A(Z_\tjm,\theta_1,\theta_2)\\
2^{-1}\partial_2 L_H(Z_\tjm,\theta_1,\theta_2,\theta_3)
\end{array}\right)\right]
\\&=&
\sum_{j=1}^nS(Z_\tjm,\theta_1,\theta_3)^{-1}
\left[\cald_j(\theta_1',\theta_2',\theta_3'),\>
\left(\begin{array}{c}\partial_2A(Z_\tjm,\theta_1,\theta_2)\\
2^{-1}H_x(z,\theta_3)[\partial_2A(Z_\tjm,\theta_2)]
\end{array}\right)\right]
\eeas
\beas
\widetilde{\Psi}_2(\theta_1,\theta_2,\theta_3,\theta_1',\theta_2',\theta_3')
&=&
\sum_{j=1}^nS(Z_\tjm,\theta_1,\theta_3)^{-1}
\left[\widetilde{\cald}_j(\theta_1',\theta_2',\theta_3'),\>
\left(\begin{array}{c}\partial_2A(Z_\tjm,\theta_2)\\
2^{-1}\partial_2 L_H(Z_\tjm,\theta_1,\theta_2,\theta_3)
\end{array}\right)\right]
\\&=&
\sum_{j=1}^nS(Z_\tjm,\theta_1,\theta_3)^{-1}
\left[\widetilde{\cald}_j(\theta_1',\theta_2',\theta_3'),\>
\left(\begin{array}{c}\partial_2A(Z_\tjm,\theta_2)\\
2^{-1}H_x(z,\theta_3)[\partial_2A(Z_\tjm,\theta_2)]
\end{array}\right)\right],
\eeas
where 
\begin{en-text}
\beas 
\widetilde{\cald}_j(\theta_1^*,\theta_2^*,\theta_3^*)
&=&
\left(\begin{array}{c}
h^{-1/2}\big(\xi^{(\ref{201906030126})}_j+\xi^{(\ref{201906030150})}_j\big)\y
h^{-3/2}\big(\xi_j^{(\ref{201906030334})}+\xi_j^{(\ref{201906030345})}\big)
\end{array}\right).
\eeas
\end{en-text}
\beas 
\widetilde{\cald}_j(\theta_1^*,\theta_2^*,\theta_3^*)
&=&
\left(\begin{array}{c}
{\cred\xi^{(\ref{201906030126})}_j+\xi^{(\ref{201906030150})}_j}\y
h^{-3/2}\big(\xi_j^{(\ref{201906030334})}+\xi_j^{(\ref{201906030345})}\big)
\end{array}\right),
\eeas

\beas
\Psi_3(\theta_1,\theta_2,\theta_3,\theta_1',\theta_2',\theta_3')
&=&
\sum_{j=1}^nS(Z_\tjm,\theta_1,\theta_3)^{-1}
\bigg[\cald_j(\theta_1',\theta_2',\theta_3')
\\&&\qquad\otimes
\left(\begin{array}{c}0\\
\partial_3 H(Z_\tjm,\theta_3)
+2^{-1}h\partial_3 L_H(Z_\tjm,\theta_1,\theta_2,\theta_3)
\end{array}\right)\bigg]
\eeas
\beas
\widetilde{\Psi}_3(\theta_1,\theta_2,\theta_3,\theta_1',\theta_2',\theta_3')
&=&
\sum_{j=1}^nS(Z_\tjm,\theta_1,\theta_3)^{-1}
\bigg[\widetilde{\cald}_j(\theta_1',\theta_2',\theta_3')
\\&&\qquad\otimes
\left(\begin{array}{c}0\\
\partial_3 H(Z_\tjm,\theta_3)
+2^{-1}h\partial_3 L_H(Z_\tjm,\theta_1,\theta_2,\theta_3)
\end{array}\right)\bigg]
\eeas
\beas
\Psi_{3,1}(\theta_1,\theta_3,\theta_1',\theta_2',\theta_3')
&=&
\sum_{j=1}^nS(Z_\tjm,\theta_1,\theta_3)^{-1}
\bigg[\cald_j(\theta_1',\theta_2',\theta_3'),\>
\left(\begin{array}{c}0\\
\partial_3 H(Z_\tjm,\theta_3)
\end{array}\right)\bigg]
\eeas
\beas
\widetilde{\Psi}_{3,1}(\theta_1,\theta_3,\theta_1',\theta_2',\theta_3')
&=&
\sum_{j=1}^nS(Z_\tjm,\theta_1,\theta_3)^{-1}
\bigg[\widetilde{\cald}_j(\theta_1',\theta_2',\theta_3'),\>
\left(\begin{array}{c}0\\
\partial_3 H(Z_\tjm,\theta_3)
\end{array}\right)\bigg]
\eeas
\beas
\Psi_{3,2}(\theta_1,\theta_2,\theta_3,\theta_1',\theta_2',\theta_3')
&=&
\sum_{j=1}^nS(Z_\tjm,\theta_1,\theta_3)^{-1}
\bigg[\cald_j(\theta_1',\theta_2',\theta_3'),\>
\left(\begin{array}{c}0\\
2^{-1}h\partial_3 L_H(Z_\tjm,\theta_1,\theta_2,\theta_3)
\end{array}\right)\bigg]
\eeas
\beas
\Psi_{3,3}(\theta_1,\theta_3,\theta_1',\theta_2',\theta_3')
&=&
\half\sum_{j=1}^n\big(S^{-1}(\partial_3S)S^{-1}\big)(Z_\tjm,\theta_1,\theta_3)
\big[\cald_j(\theta_1',\theta_2',\theta_3')^{\otimes2}-S(Z_\tjm,\theta_1,\theta_3)\big]
\eeas
\beas
\Psi_{33,1}(\theta_1,\theta_2,\theta_3)
&=&
-\sum_{j=1}^nS(Z_\tjm,\theta_1,\theta_3)^{-1}
\left[
\left(\begin{array}{c}0\\
\partial_3H(Z_\tjm,\theta_3)+2^{-1}h\partial_3L_H(Z_\tjm,\theta_1,\theta_2,\theta_3)
\end{array}\right)^{\otimes2}
\right]
\eeas
\beas
\Psi_{33,2}(\theta_1,\theta_2,\theta_3,\theta_1',\theta_2',\theta_3')
&=&
\sum_{j=1}^nS(Z_\tjm,\theta_1,\theta_3)^{-1}
\bigg[\cald_j(\theta_1',\theta_2',\theta_3')
\\&&\qquad\otimes
\left(\begin{array}{c}0\\
\partial_3^2 H(Z_\tjm,\theta_3)
+2^{-1}h\partial_3^2 L_H(Z_\tjm,\theta_1,\theta_2,\theta_3)
\end{array}\right)\bigg]
\eeas
\beas
\Psi_{33,3}(\theta_1,\theta_3)
&=&
-\half\sum_{j=1}^n
\big\{\big(S^{-1}(\partial_3S)S^{-1}\big)(Z_\tjm,\theta_1,\theta_3)
\big[\partial_3S(Z_\tjm,\theta_1,\theta_3)\big]\big\}
\eeas
\beas
\Psi_{33,4}(\theta_1,\theta_2,\theta_3,\theta_1',\theta_2',\theta_3')
&=&
- {\colorr2}\sum_{j=1}^n
S^{-1}(\partial_3S)S^{-1}(Z_\tjm,\theta_3)
\left[\cald_j(\theta_1',\theta_2',\theta_3')\right.
\\&&\hspace{3cm}\otimes
\left.
\left(\begin{array}{c}0\\
\partial_3H(Z_\tjm,\theta_3)+2^{-1}h\partial_3L_H(Z_\tjm,\theta_1,\theta_2,\theta_3)
\end{array}\right)
\right]
\eeas
\beas
\Psi_{33,5}(\theta_1,\theta_3,\theta_1',\theta_2',\theta_3')
&=&
\half \sum_{j=1}^n\partial_3\big\{
\big(S^{-1}(\partial_3S)S^{-1}\big)(Z_\tjm,\theta_1,\theta_3)\big\}
\big[\cald_j(\theta_1',\theta_2',\theta_3')^{\otimes2}-S(Z_\tjm,\theta_1,\theta_3)\big].
\eeas

\subsection{Section \ref{202001141623}}
\beas 
\bbH_n^{(1)}(\theta_1)
 &=& 
 -\half\sum_{j=1}^n \bigg\{
S(Z_\tjm,\theta_1,\hat{\theta}_3^0)^{-1}\big[\cald_j(\theta_1,\hat{\theta}_2^0,\hat{\theta}_3^0)^{\otimes2}\big]
 +\log\det S(Z_\tjm,\theta_1,\hat{\theta}_3^0)\bigg\}.
\eeas
\beas 
\bbH^{(2,3)}_n(\theta_2,\theta_3)
 &=& 
 -\half\sum_{j=1}^n 
 \hat{S}(Z_\tjm,\hat{\theta}_3^0)^{-1}\big[\cald_j(\hat{\theta}_1^0,\theta_2,\theta_3)^{\otimes2}\big].
\eeas
Recall $\hat{S}(z,\theta_3)=S(z,\hat{\theta}_1^0,\theta_3)$. 

\beas
\Phi_{22,1}(\theta_1,\theta_3,\theta_1',\theta_2',\theta_3')
&=&
-\sum_{j=1}^nS(Z_\tjm,\theta_1,\theta_3)^{-1}
\left[
\left(\begin{array}{c}\partial_2A(Z_\tjm,\theta_2')\\
2^{-1}\partial_2L_H(Z_\tjm,\theta_1',\theta_2',\theta_3')
\end{array}\right)^{\otimes2}
\right]
\eeas
\beas&&
\Phi_{22,2}(\theta_1,\theta_3,\theta_1',\theta_2',\theta_3',
\theta_2'',\theta_3'')
\nn\\&=&
\sum_{j=1}^nS(Z_\tjm,\theta_1,\theta_3)^{-1}
\left[\cald_j(\theta_1',\theta_2',\theta_3'),\>
\left(\begin{array}{c}\partial_2^2A(Z_\tjm,\theta_2'')\\
2^{-1}H_x(Z_\tjm,\theta_3'')\big[\partial_2^2A(Z_\tjm,\theta_2'')\big]
\end{array}\right)
\right]
\eeas
\beas&&
\widetilde{\Phi}_{22,2}(\theta_1,\theta_3,\theta_1',\theta_2',\theta_3',
\theta_2'',\theta_3'')
\nn\\&=&
\sum_{j=1}^nS(Z_\tjm,\theta_1,\theta_3)^{-1}
\left[\widetilde{\cald}_j(\theta_1',\theta_2',\theta_3'),\>
\left(\begin{array}{c}\partial_2^2A(Z_\tjm,\theta_2'')\\
2^{-1}H_x(Z_\tjm,\theta_3'')\big[\partial_2^2A(Z_\tjm,\theta_2'')\big]
\end{array}\right)
\right]
\eeas

\beas &&
\Phi_{23,1}(\theta_1,\theta_3,\theta_1',\theta_2',\theta_3',\theta_2'',\theta_3'')
\\&=&
-\sum_{j=1}^nS(Z_\tjm,\theta_1,\theta_3)^{-1}
\left[\left(\begin{array}{c}0\\
2^{-1}\partial_3L_H(Z_\tjm,\theta_1',\theta_2',\theta_3')
\end{array}\right)\>\right.
\\&&\hspace{5cm}\left.\otimes
\left(\begin{array}{c}\partial_2A(Z_\tjm,\theta_2'')\\
2^{-1}H_x(Z_\tjm,\theta_3'')\big[\partial_2A(Z_\tjm,\theta_2'')\big]
\end{array}\right)
\right]
\eeas
\beas&&
\Phi_{23,2}(\theta_1,\theta_3,\theta_1',\theta_2',\theta_3',\theta_2'',\theta_3'')
\\&=&
\sum_{j=1}^nS(Z_\tjm,\theta_1,\theta_3)^{-1}
\left[\cald_j(\theta_1',\theta_2',\theta_3'),\>
\left(\begin{array}{c}0\\
2^{-1}\partial_3H_x(Z_\tjm,\theta_3'')\big[\partial_2A(Z_\tjm,\theta_2'')\big]
\end{array}\right)
\right]
\eeas
\beas
\Phi_{33,1}(\theta_1,\theta_3,\theta_1',\theta_2',\theta_3')
&=&
-\sum_{j=1}^nS(Z_\tjm,\theta_1,\theta_3)^{-1}
\left[
\left(\begin{array}{c}0\\
\partial_3H(Z_\tjm,\theta_3')+2^{-1}h\partial_3L_H(Z_\tjm,\theta_1',\theta_2',\theta_3')
\end{array}\right)^{\otimes2}
\right]
\eeas
\beas
\Phi_{33,2}(\theta_1,\theta_3,\theta_1',\theta_2',\theta_3',\theta_1'',\theta_2'',\theta_3'')
&=&
\sum_{j=1}^nS(Z_\tjm,\theta_1,\theta_3)^{-1}
\bigg[\cald_j(\theta_1',\theta_2',\theta_3')
\\&&\qquad\otimes
\left(\begin{array}{c}0\\
\partial_3^2 H(Z_\tjm,\theta_3'')
+2^{-1}h\partial_3^2 L_H(Z_\tjm,\theta_1'',\theta_2'',\theta_3'')
\end{array}\right)\bigg]
\eeas

\beas 
\Psi_{1,1}(\theta_1,\theta_2,\theta_3,\theta_1',\theta_2',\theta_3')
&=&
\Psi_1(\theta_1,\theta_2,\theta_3,\theta_1',\theta_2',\theta_3')
\\&=&
\sum_{j=1}^nS(Z_\tjm,\theta_1,\theta_3)^{-1}
\left[\cald_j(\theta_1',\theta_2',\theta_3'),\>
\left(\begin{array}{c}0\\
2^{-1}\partial_1L_H(Z_\tjm,\theta_1,\theta_2,\theta_3)
\end{array}\right)
\right]
\eeas
\beas 
\Psi_{1,2}(\theta_1,\theta_3,\theta_1',\theta_2',\theta_3')
&=&
\half\sum_{j=1}^n
\big(S^{-1}(\partial_1S))S^{-1}\big)(Z_\tjm,\theta_1,\theta_3)
\big[\cald_j(\theta_1',\theta_2',\theta_3')^{\otimes2}-S(Z_\tjm,\theta_1,\theta_3)\big]
\eeas
\beas 
\Psi_{11,1}(\theta_1,\theta_3,\theta_1',\theta_2',\theta_3')
&=&
\sum_{j=1}^nS(Z_\tjm,\theta_1,\theta_3)^{-1}
\left[
\left(\begin{array}{c}0\\
2^{-1}\partial_1L_H(Z_\tjm,\theta_1',\theta_2',\theta_3')
\end{array}\right)^{\otimes2}
\right]
\eeas
\beas 
\Psi_{11,2}(\theta_1,\theta_3,\theta_1',\theta_2',\theta_3',\theta_1'',\theta_2'',\theta_3'')
&=&
\sum_{j=1}^nS(Z_\tjm,\theta_1,\theta_3)^{-1}
\bigg[\cald_j(\theta_1',\theta_2',\theta_3'),\>
\\&&\otimes
\left(\begin{array}{c}0\\
2^{-1}\partial_1^2L_H(Z_\tjm,\theta_1'',\theta_2'',\theta_3'')
\end{array}\right)
\bigg]
\eeas
\beas 
\Psi_{11,3}(\theta_1,\theta_3,\theta_1',\theta_2',\theta_3')
&=&
\sum_{j=1}^n
\partial_1\big\{S^{-1}(\partial_1S)S^{-1}(Z_\tjm,\theta_1,\theta_3)\big\}
\big[\cald_j(\theta_1',\theta_2',\theta_3')^{\otimes2}-S(Z_\tjm,\theta_1',\theta_3')\big]
\eeas
\beas
\Psi_{11,4}(\theta_1,\theta_3)
&=&
\sum_{j=1}^n(S^{-1}(\partial_1S)S^{-1})(Z_\tjm,\theta_1,\theta_3)
\big[\partial_1S(Z_\tjm,\theta_1, \theta_3)\big]
\eeas
\beas
\Psi_{11,5}(\theta_1,\theta_3,\theta_1',\theta_2',\theta_3',\theta_1'',\theta_2'',\theta_3'')
&=&
\sum_{j=1}^n
(S^{-1}(\partial_1S)S^{-1})(Z_\tjm,\theta_1,\theta_3)
\bigg[\cald_j(\theta_1',\theta_2',\theta_3'),\>
\\&&
\hspace{3cm}\otimes
\left(\begin{array}{c}0\\
2^{-1}\partial_1L_H(Z_\tjm,\theta_1'',\theta_2'',\theta_3'')
\end{array}\right)
\bigg].
\eeas

\subsection{Section \ref{202001141632}}
\beas 
\bbH_n(\theta)
 &=& 
 -\half\sum_{j=1}^n \bigg\{
S(Z_\tjm,\theta_1,\theta_3)^{-1}\big[\cald_j(\theta_1,\theta_2,\theta_3)^{\otimes2}\big]
 +\log\det S(Z_\tjm,\theta_1,\theta_3)\bigg\}
\eeas
for $\theta=(\theta_1,\theta_2,\theta_3)$. 

\beas 
\bbD_n(\theta_1,\theta_2,\theta_3,\theta_1',\theta_2',\theta_3')
&=&
\bbH_n(\theta_1,\theta_2,\theta_3)-\bbH_n(\theta_1',\theta_2',\theta_3').
\eeas

\beas &&
\bbD_n^{[1]}(\theta_1,\theta_2,\theta_3,\theta_1',\theta_2',\theta_3')
\\&=&
-\half\sum_{j=1}^nS(Z_\tjm,\theta_1,\theta_3)^{-1}
\left[
\left(\begin{array}{c}
h^{1/2}\big(A(Z_\tjm,\theta_2)-A(Z_\tjm,\theta_2')\big)
\\
\left\{\begin{array}{c}
h^{-1/2}\big(H(Z_\tjm,\theta_3)-H(Z_\tjm,\theta_3')\big)
\\
+2^{-1}h^{1/2}\big(L_H(Z_\tjm,\theta_1,\theta_2,\theta_3)-L_H(Z_\tjm,\theta_1',\theta_2',\theta_3')\big)
\end{array}\right\}
\end{array}\right)^{\otimes2}\right]
\eeas
\beas
\bbD_n^{[2]}(\theta_1,\theta_3,\theta_1',\theta_2',\theta_3')
&=&
h^{-1/2}\sum_{j=1}^nS(Z_\tjm,\theta_1,\theta_3)^{-1}
\left[
\cald_j(\theta_1',\theta_2',\theta_3'),\>
\left(\begin{array}{c}
0\\
H(Z_\tjm,\theta_3)-H(Z_\tjm,\theta_3')
\end{array}\right)\right]
\eeas
\beas&&
\bbD_n^{[3]}(\theta_1,\theta_2,\theta_3,\theta_1',\theta_2',\theta_3')
\\&=&
h^{1/2}\sum_{j=1}^nS(Z_\tjm,\theta_1,\theta_3)^{-1}
\left[
\cald_j(\theta_1',\theta_2',\theta_3'),\>
\left(\begin{array}{c}
A(Z_\tjm,\theta_2)-A(Z_\tjm,\theta_2')\\
2^{-1}\big(L_H(Z_\tjm,\theta_1,\theta_2,\theta_3)-L_H(Z_\tjm,\theta_1',\theta_2',\theta_3')\big)
\end{array}\right)\right]
\eeas

\beas
\bbD_n^{[4]}(\theta_1,\theta_3,\theta_1',\theta_2',\theta_3')
&=&
-\half\sum_{j=1}^n\bigg\{\big(S(Z_\tjm,\theta_1,\theta_3)^{-1}-S(Z_\tjm,\theta_1',\theta_3')^{-1}\big)
\big[\cald_j(\theta_1',\theta_2',\theta_3')^{\otimes2}\big]
\\&&
+\log\frac{\det S(Z_\tjm,\theta_1,\theta_3)}{\det S(Z_\tjm,\theta_1',\theta_3')}\bigg\}
\eeas

\end{document}

%% file: nakamacro291103++.tex
\def\infm{{\infty\text{--}}}
\def\inftym{\infm}
\def\koko{{\coloroy{koko}}}
\def\bd{\begin{description}}
\def\ed{\end{description}}

\def\D2{\bbD_{2,\infty-}}

\def\tj{{t_j}}
\def\tjm{{t_{j-1}}}

\def\D{{\bf D}}

\def\F{{\bf F}}

\def\cald{{\cal D}}

\def\calf{{\cal F}}

\def\calh{{\cal H}}

\def\calx{{\cal X}}

\def\ds{\displaystyle}
\def\yeq{\>=\>}

\def\yl{\><\>}

\def\sfd{{\sf d}}
\def\sfp{{\sf p}}
\def\sfr{{\sf r}}

\def\simleq{\ \raisebox{-.7ex}{$\stackrel{{\textstyle <}}{\sim}$}\ }

\def\ep{\epsilon}
\def\half{\frac{1}{2}}

\def\up{\uparrow}
\def\down{\downarrow}

\def\y{\vspace*{3mm}\\}
\def\halflineskip{\vspace*{3mm}}
\def\nn{\nonumber}
\def\be{\begin{equation}}
\def\ee{\end{equation}}
\def\bea{\begin{eqnarray}}
\def\eea{\end{eqnarray}}
\def\beas{\begin{eqnarray*}}
\def\eeas{\end{eqnarray*}}
\def\bi{\begin{itemize}}
\def\ei{\end{itemize}}
\def\im{\item}
\def\bd{\begin{description}}
\def\ed{\end{description}}
\def\dotc{\stackrel{\circ}{C}}

\newcommand{\bbD}{{\mathbb D}}
\newcommand{\bbE}{{\mathbb E}}
\newcommand{\bbF}{{\mathbb F}}

\newcommand{\bbH}{{\mathbb H}}

\newcommand{\bbN}{{\mathbb N}}

\newcommand{\bbR}{{\mathbb R}}

\newcommand{\bbY}{{\mathbb Y}}
\newcommand{\bbZ}{{\mathbb Z}}